\documentclass[11pt, a4paper, openright]{report}

\usepackage{graphicx} 
\usepackage{amsmath, amsthm,  amsfonts, amscd, amssymb}
\usepackage[twoside, margin={20mm,35mm}, bindingoffset=20mm]{geometry}
\usepackage{url}
\usepackage{mathrsfs}
\usepackage[all]{xy}
\DeclareMathOperator{\hol}{hol}
\DeclareMathOperator{\Map}{Map}
\DeclareMathOperator{\Aut}{Aut}

\DeclareMathOperator{\Ad}{Ad}
\DeclareMathOperator{\ev}{ev}
\DeclareMathOperator{\Spin}{Spin}	
\DeclareMathOperator{\Diff}{Diff}

\theoremstyle{plain}
\newtheorem{theorem}{Theorem}[section]
\newtheorem{lemma}[theorem]{Lemma}
\newtheorem{proposition}[theorem]{Proposition}

\newtheorem{corollary}[theorem]{Corollary}

\theoremstyle{definition}
\newtheorem{definition}[theorem]{Definition}

\theoremstyle{remark}

\numberwithin{equation}{section}
\numberwithin{figure}{section}

\newcommand{\cG}{{\mathcal G}}

\newcommand{\fg}{{\mathfrak g}}

\newcommand{\RR}{{\mathbb R}}
\newcommand{\ZZ}{{\mathbb Z}}

\renewcommand{\a}{\alpha}
\renewcommand{\b}{\beta}

\renewcommand{\d}{\delta}

\newcommand{\LGS}{LG\rtimes S^1}
\newcommand{\Lgs}{L\fg \rtimes i\RR}
\newcommand{\<}{\langle}
\renewcommand{\>}{\rangle}
\newcommand{\ds}{\displaystyle}
\newcommand{\DD}{\textsl{DD}}
\newcommand{\ssv}{\scriptscriptstyle{\vee}}

\renewcommand{\theequation}{\arabic{chapter}.\arabic{section}.\arabic{equation}}

\begin{document}

\pagenumbering{roman}

\begin{titlepage}
%
%
\vspace*{8.0mm}
%
%
\begin{LARGE}
\begin{bf}
\begin{center}
Loop Groups, Higgs Fields and Generalised String Classes
\end{center}
\end{bf}
\end{LARGE}
%
%
\vspace{5.0mm}
%
%
\begin{large}
\begin{center}
Raymond Vozzo
\end{center}
%
%
\vspace*{\fill}
\vspace{2.0mm}
%
%
\begin{center}
Thesis
submitted for the degree
of\\
Doctor of Philosophy\\
in\\
Pure Mathematics\\
at\\
The University of Adelaide\\
\end{center}

%
%
\vspace*{\fill}
\vspace{4.0mm}
%
%
\begin{center}
School of Mathematical Sciences
\end{center}
%
%
\vspace{4.0mm}
\begin{center} 
\includegraphics[height=2cm]{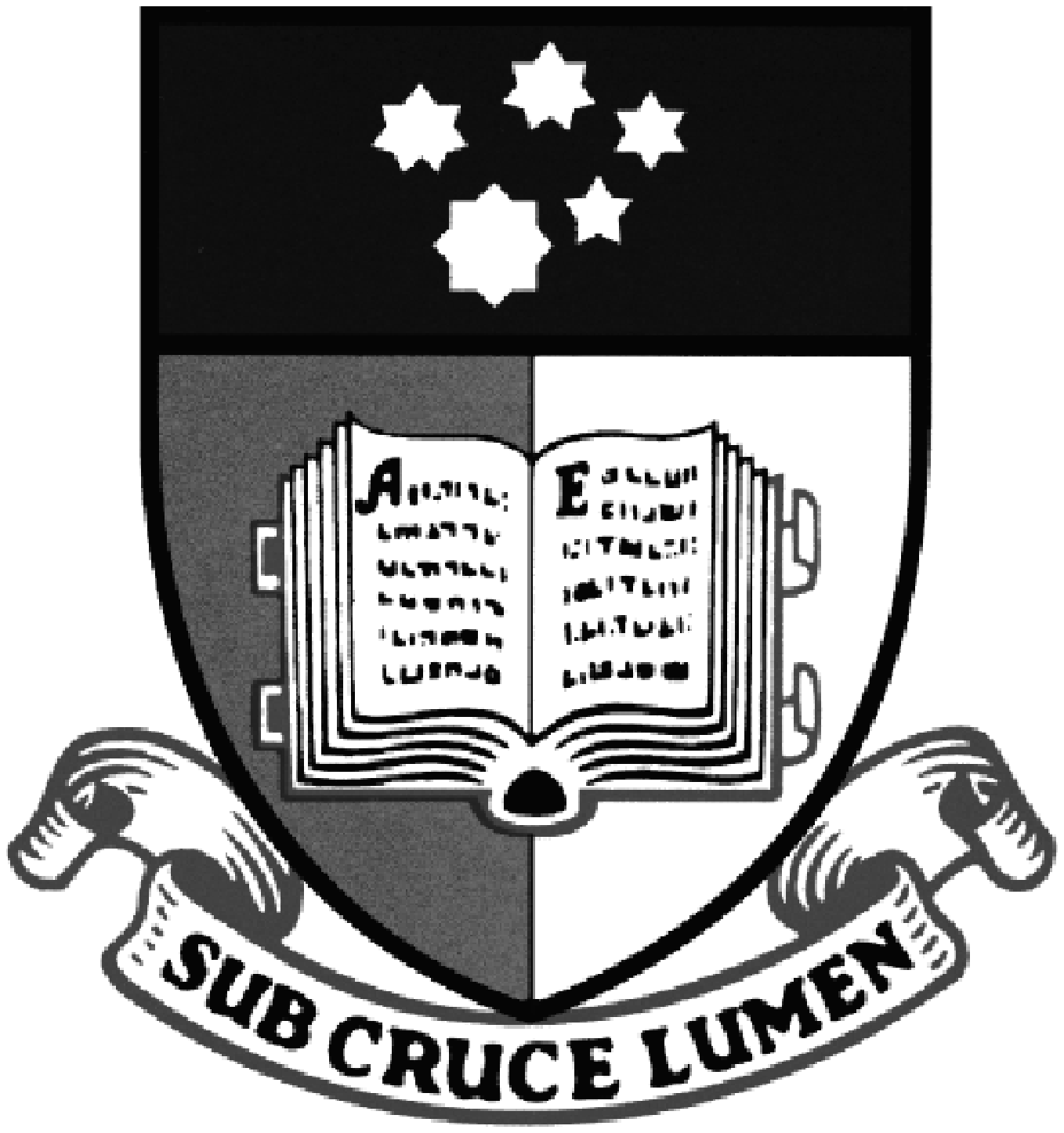}
\end{center}
%
%
%
%
\begin{center}
June 22, 2009
\end{center}
\end{large}
%
%
\end{titlepage}

\tableofcontents

\chapter*{Abstract}
\typeout{Abstract}
\label{ch:abstract}
\addcontentsline{toc}{chapter}{Abstract}

We consider various generalisations of the string class of a loop group bundle. The string class is the obstruction to lifting a bundle whose structure group is the loop group $LG$ to one whose structure group is the Kac-Moody central extension of the loop group.

We develop a notion of higher string classes for bundles whose structure group is the group of based loops, $\Omega G$. In particular, we give a formula for characteristic classes in odd dimensions for such bundles which are associated to characteristic classes for $G$-bundles in the same way that the string class is related to the first Pontrjagyn class of a certain $G$-bundle associated to the loop group bundle in question. This provides us with a theory of characteristic classes for $\Omega G$-bundles analogous to Chern-Weil theory in finite dimensions. This also gives us a geometric interpretation of the well-known transgression map $H^{2k}(BG) \to H^{2k-1}(G).$

We also consider the obstruction to lifting a bundle whose structure group is not the loop group but the semi-direct product of the loop group with the circle, $LG \rtimes S^1$. We review the theory of bundle gerbes and their application to central extensions and lifting problems and use these methods to obtain an explicit expression for the de Rham representative of the obstruction to lifting such a bundle. We also relate this to a generalisation of the so-called `caloron correspondence' (which relates $LG$-bundles over $M$ to $G$-bundles over $M \times S^1$) to a correspondence which relates $LG \rtimes S^1$-bundles over $M$ to $G$-bundles over $S^1$-bundles over $M$.

\chapter*{Signed Statement}
\addcontentsline{toc}{chapter}{Signed Statement}

This work contains no material which has been accepted
for the award of any other degree or diploma in any university or other
tertiary institution and, to the best of my knowledge and belief,
contains no material previously
published or written by another person, except where due reference has
been made in
the text.

\vspace{5mm}

\noindent
I consent to this copy of my thesis, when deposited in the University
Library, being available for loan and photocopying, subject to the provisions of the Copyright Act 1968.

\vspace{20mm}

\noindent
{\sc SIGNED}: {\tt .......................} {\sc DATE}: {\tt .......................}

\chapter*{Acknowledgements}
\typeout{Acknowledgements}
\label{ch:acknowledgements}
\addcontentsline{toc}{chapter}{Acknowledgements}

I would like to thank my supervisors, Michael and Mathai, for their invaluable assistance and guidance, without which this project surely would never have finished (or, at least, would have taken a significantly longer time). I have learned more about how Mathematics is done in the last three years from talking to them than I thought there even was to know. Thanks are also due to Danny Stevenson for many useful discussions (and possibly even more so for his willingness to help me).

Thanks are also due to all the staff in the admin office, whose hard work has made my life much easier over the past few years.

I would also like to thank my fellow postgrads, all of whom have made the last four years of my life more than bearable (to say the least!). I shall not endeavour to name everyone here (since I want to fit this whole thing on one page!), however special thanks are due in particular to Ric Green, David Roberts, David Butler, Jonathan Tuke, Rongmin Lu, Glenis Crane and Jessica Kasza. Their readiness to place their personal health at risk by ingesting dangerous amounts of caffeine (with the possible exception of Rongmin and David of course) just so we can take a break and discuss something other than Mathematics for 15 (or 45 as the case may be) minutes is greatly appreciated. By the same token, special thanks are due to David R and Ric for many, many helpful discussions on not only my research topic but Mathematics and Physics in general. Over the past three or four years we have discovered that other people's (Mathematical) problems are much more interesting than one's own and thanks to David and Ric I believe I have profited immensely from this fact.


Finally, I would like to thank my family for their continual love and support. Without them this PhD would have remained a dream. I am immensely grateful to my parents, Armando and Lucy, for their constant encouragement and to my brother Jonathon and my sister Nicola, who have helped me keep some semblance of my childhood (despite being nearly 25) and have never questioned why their brother doesn't have a job or a home of his own. And to my beautiful fianc\'ee Emily, for being my beautiful fianc\'ee, and for putting up with so much incomprehensible Maths over the years and whose love and patience has made all this possible.

\newpage

\phantom{ }



\pagenumbering{arabic}

\chapter{Introduction}\label{C:intro}



String structures first appeared in Killingback's paper \cite{Killingback:1987} as a string theory version of the well-known spin structures that are important in quantum field theory. The results came out of a study of global anomalies in the worldsheet of a string and the idea was motivated by an observation of Witten \cite{Witten:1985} that the Dirac-Ramond operator in string theory can be considered as Dirac-type operator on the loop space. 

Recall that if one is given a principal $SO(n)$-bundle (for example the frame bundle of a manifold), a spin structure is given by a lifting of the structure group of this bundle to its simply connected double cover $\Spin(n).$ Killingback's idea then, is to replace the bundles which appear in the definition of a spin structure with an infinite-dimensional bundle whose structure group is the loop group of $\Spin(n)$ and consider a lifting of this bundle. More generally, if $G$ is a compact Lie group and $LG$ is its loop group, we could consider lifting any $LG$-bundle $P \to M$ to a bundle whose structure group is the central extension of $LG.$ It turns out that the obstruction $s(P)$ to the existence of such a lift is an element of the degree three cohomology of the base, $H^3 (M, \ZZ).$ Killingback proved that, in the case where the $LG$-bundle $P$ is in fact given by taking loops in a principal $G$-bundle $Q \to X,$ this obstruction class is the transgression of the first Pontrjagyn class of $Q.$ That is,
$$
s(P) = \int_{S^1} \ev^*p_1(Q),
$$
where $\ev \colon S^1 \times LX \to X$ is the evaluation map. The class $s(P) \in H^3(LX, \ZZ)$ is called the string class of $P.$ The link with spin structures and Witten's observation regarding the Dirac-Ramond operator is given by noticing that in quantum field theory the Dirac operator can only be defined if spacetime is spin and correspondingly in string theory the Dirac-Ramond operator can only be defined if spacetime is string (i.e.\! has a string structure).

The present work grew out of an attempt to answer some questions naturally arising from some papers concerning string structures and loop group bundles. In \cite{Murray:2003} Murray and Stevenson use techniques from the theory of bundle gerbes to give an explicit formula for a representative in de Rham cohomology of the string class of a general $LG$-bundle and provide a link with previous work on monopoles. The theory of gerbes was first introduced by Giraud \cite{Giraud:1971} and studied extensively in Brylinski's book \cite{Brylinski:1993}. Gerbes provide a geometric realisation for degree three cohomology in an analogous way to the way in which line bundles (or $U(1)$-bundles) provide a geometric realisation of degree two cohomology. Gerbes are essentially sheaves of groupoids satisfying certain descent conditions but can be tricky to work with in practice. A much more appealing (at least from a differential geometric point of view) approach to the theory of gerbes, called bundle gerbes, was introduced by Murray \cite{Murray:1996}. These have been studied further (see for example \cite{Carey:2000a, Gomi:2003, Meinrenken:2003, Mickelsson:2008, Murray:2000}) and have found applications in physics as well as differential geometry (see for example \cite{Bouwknegt:2002, Carey:1997, Carey:2005, Ekstrand:2000, Schreiber:2007}). Insofar as a gerbe can be considered a sheaf of groupoids, bundle gerbes can be viewed as bundles of groupoids. They have a degree three characteristic class associated with them, called the Dixmier-Douady class, which can be described in terms of cocycles. However, one can also define a notion of connection and curvature (more precisely, 3-curvature) for a bundle gerbe and, using differential geometric methods, obtain a differential form representative for the image in real cohomology of the Dixmier-Douady class in analogy with the way the Chern class of a $U(1)$-bundle is represented in real cohomology by the curvature of the bundle. Bundle gerbes arise very naturally in lifting problems such as the string structure example. This is the approach taken in \cite{Murray:2003} where a de~Rham representative of the string class for a loop group bundle $P$ is given in terms of data on the bundle. Namely, the authors find that the string class is given by
$$
s(P) = -\frac{1}{4\pi^2} \int_{S^1} \<F, \nabla \Phi \> \, d\theta
$$
where $F$ is the curvature of $P,$ $\nabla \Phi$ is the covariant derivative of a Higgs field for $P$ and the bracket is the Killing form suitably normalised. They also extend Killingback's result -- that is, giving the string class in terms of the Pontrjagyn class for some $G$-bundle -- by using the so-called `caloron correspondence' (which first appeared in \cite{Garland:1988}) which relates $LG$-bundles over $M$ to $G$-bundles over $M\times S^1.$ In particular, there is a bijective correspondence between isomorphism classes of principal $LG$-bundles over $M$ and isomorphism classes of principal $G$-bundles over $M \times S^1$ and if $P \to M$ is an $LG$-bundle and $\widetilde{P} \to M \times S^1$ is its corresponding $G$-bundle, then the authors find that the string class of $P$ is given by integrating the first Pontrjagyn class of $\widetilde{P}:$
$$
s(P) = \int_{S^1} p_1(\widetilde{P}).
$$

The first formula above can be used to recover the result from \cite{Carey:1991} in which the authors calculate the string class for the universal $\Omega G$-bundle\footnote{Actually, in \cite{Carey:1991} the authors work with the group of smooth maps from the interval into $G$ whose endpoints agree. In this thesis we extend their work to the group of smooth maps from the circle into $G$.} (where $\Omega G$ is the based loop group) and show that the string class is a characteristic class for loop bundles (that is, $\Omega G$-bundles of the form $\Omega Q \to \Omega X$ for some $G$-bundle $Q \to X$). A model for the classifying space of $\Omega G$ is given by the group $G$ itself and $H^3(G, \ZZ) = \ZZ$ so it is not unreasonable to expect the string class in this case to be the generator of this group. This is in fact true and it is shown that the string class for any loop bundle is given by the pull-back of this class by a classifying map for the bundle.

This thesis deals with two natural questions which arise when one considers these results. The first concerns the relationship between the string class and the Pontrjagyn class and the fact that the string class is a characteristic class for loop bundles. It is natural, firstly, to look for a way to generalise this to $\Omega G$-bundles which are not necessarily loop bundles but, also, it seems possible that there is a more general theory of characteristic classes for loop group bundles which is related to characteristic class theory for $G$-bundles (i.e.\! Chern-Weil theory). In the first part of this thesis we provide  answers to these problems. We give a generalisation of the result from \cite{Carey:1991} to $\Omega G$-bundles which are not loop bundles, that is, we show that the string class is a characteristic class. We then develop a notion of higher string classes for $\Omega G$-bundles which are also characteristic classes and are related to characteristic classes for $G$-bundles. In particular, we develop a kind of Chern-Weil theory for $\Omega G$-bundles which gives characteristic classes from invariant polynomials on the Lie algebra $\fg$ of $G$ and data on the $\Omega G$-bundle. This theory side-steps the complications which arise when trying to define the Chern-Weil map directly for bundles with infinite-dimensional structure group (for example, see \cite{Paycha:2004}). It also provides a geometric interpretation of the well-known transgression map $\tau \colon H^{2k}(BG) \to H^{2k-1}(G).$

The next question which it is natural to ask concerns the caloron correspondence described above (i.e.\! the correspondence between $LG$-bundles over $M$ and $G$-bundles over $M \times S^1$). In trying to find a formula for the string class in terms of the Pontrjagyn class of a $G$-bundle (as in \cite{Murray:2003}) one finds that it is necessary to make use of the caloron correspondence. So it is natural then to ask what kind of correspondence exists in the case where the $G$-bundle is not over $M \times S^1$ but over a non-trivial principal $S^1$-bundle over $M$ and, further, whether the methods of bundle gerbes can be applied to the lifting problem in this case. In fact, the first part of this question has been answered in \cite{Bergman:2005} in connection with the Kaluza-Klein reduction of M-theory to type IIA supergravity. It turns out that there is a bijective correspondence between isomorphism classes of $G$-bundles over $S^1$-bundles and classes of bundles whose structure group is not the loop group, but the semi-direct product $LG \rtimes S^1.$ In the latter part of this thesis we prove that this correspondence also holds on the level of connections (as in the case of a trivial circle bundle) and consider the lifting problem for an $LG \rtimes S^1$-bundle. We use the methods of \cite{Murray:2003} to find a de Rham representative for the image in real cohomology of the class which is the obstruction to the existence of this lift. We also provide a calculation of this class using a different method introduced by Gomi \cite{Gomi:2003}, that of reduced splittings.

The outline of this thesis is as follows: In chapter 2 we describe the necessary background. We recall some important facts about spin structures and give an overview of Killingback's results on string structures. We also review the theory of bundle gerbes and their application to lifting problems. We then present, in some detail, the theory and results from Murray and Stevenson's paper \cite{Murray:2003}, including the calculation of the string class for a general $LG$-bundle and the correspondence between $LG$-bundles over $M$ and $G$-bundles over $M \times S^1.$ We also include the extension of Killingback's result from this paper.

In chapter 3 we show that the string class is a characteristic class for $\Omega G$-bundles (Theorem \ref{T:s(P) is char}) and generalise some of the results from chapter 2 (albeit, only in the case of the based loop group) to higher dimensions. That is, we define cohomology classes in any odd dimension which are related to characteristic classes for $G$-bundles (in the same way that the string class is related to the Pontrjagyn class) and we prove that these are themselves characteristic classes. This gives a method of finding characteristic classes for an $\Omega G$-bundle given a universal characteristic class for $G$-bundles (that is, an element of $H^*(BG)$). This is detailed in Theorem \ref{T:"Chern-Weil"}. We also provide a partial generalisation to the case of the free loop group (although here we work with the group of smooth maps from the interval into $G$ whose endpoints agree). We give a model for the universal bundle and calculate its string class.

In chapter 4 we present the calculation of the string class of an $LG \rtimes S^1$-bundle (that is, the obstruction to lifting the structure group of an $LG \rtimes S^1$-bundle to its central extension). This is given in Theorem \ref{T:LGxS^1string}. We also give the generalisation of the caloron correspondence from \cite{Bergman:2005} which relates $G$-bundles over $S^1$-bundles to $LG \rtimes S^1$-bundles. We show that this correspondence holds on the level of connections as well (Proposition \ref{P:LGxS^1 connection correspondences}). This allows us to prove a generalisation of the result from \cite{Murray:2003} relating the string class to the Pontrjagyn class of the corresponding $G$-bundle (Theorem \ref{T:LGxS^1Pont}). Finally, we briefly outline how these results can be used to gain information about the more general case of lifting a bundle whose structure group is $LG \rtimes \Diff(S^1),$ that is, where the loops in $LG$ are acted upon by general (orientation preserving) diffeomorphisms of the circle.

We make a final comment on terminology and conventions. Throughout this thesis we will work with many variations of the loop group. We give these here for convenience. The group of smooth maps $\Map(S^1, G)$ is denoted by $LG$ and the subgroup of based loops which start at the identity by $\Omega G.$ In chapter 3 we consider slightly more general variants of these groups which consist of smooth maps from the interval $[0, 2\pi]$ into $G$ whose endpoints agree. These are denoted by $L^{\ssv}G$ in the free case and $\Omega^{\ssv}G$ in the based case. Finally, the terms principal $G$-bundle and $G$-bundle are used interchangeably and all bundles are assumed to be principal bundles unless specifically stated otherwise. Also, the circle group is denoted by either $U(1)$ or $S^1$ -- we make no distinction between the two.
\chapter{String structures, bundle gerbes and Higgs fields}\label{C:background}


In this chapter we shall present the relevant background required for the rest of the thesis. Namely, we describe the existing results on string structures and develop the theory of bundle gerbes, which will feature quite heavily in the sequel.


\section{String structures}\label{S:string structures}


The existence of spinors and the Dirac operator is an essential aspect of quantum field theory. It is well known that in order to define these objects the underlying spacetime $M$ must be a spin manifold. In \cite{Witten:1985}, in a study of global anomalies, Witten shows that there occurs a global anomaly in the worldline of a supersymmetric point particle in quantum mechanics unless $M$ admits a spin structure. The analogue of this in string theory, that is, a global anomaly in the worldsheet of a string, was also studied in some detail. Killingback, in \cite{Killingback:1987}, uses these results to determine topological conditions on the spacetime $M.$ These conditions led to the definition of a so-called \emph{string structure} on $M.$ Let us first recall, then, what we mean by a \emph{spin structure} and show how to find the analogue of this in string theory.

\subsection{Spin structures}\label{SS:spin structures}

Let $M$ be an orientable manifold and $F \to M$ its frame bundle. Then $F$ is a principal $SO(n)$-bundle. There is a simply connected double cover of $SO(n),$ called $\Spin(n)$ that fits into the exact sequence
$$
0 \to \ZZ_2 \to \Spin(n) \to SO(n) \to 0.
$$
Thus we can consider lifting the frame bundle of $M$ to a principal $\Spin(n)$-bundle where by a \emph{lift} of $F \to M$ we mean a principal $\Spin(n)$-bundle $\hat{F}\to M$ such that there is a bundle map $\hat{F} \to F$ that commutes with the homomorphism $\Spin(n) \to SO(n).$ If such a lift exists, we say $M$ has a \emph{spin structure}, or simply that $M$ is \emph{spin}. More generally, we can consider any principal $SO(n)$-bundle $P\to M$ and ask for a lift of $P$ to a principal $\Spin(n)$-bundle. If a lift exists in this case we say that $P$ has a spin structure. It can be shown (see for example \cite{Lawson:1989}) that a spin structure exists for $P$ if and only if the second Stiefel-Whitney class, $w_2(P),$ vanishes.

\subsection{String structures}\label{SS:string structures}

As mentioned above, the Dirac operator, an integral element of quantum field theory, cannot be defined unless $M$ is a spin manifold. The analogue of this operator in string theory is the Dirac-Ramond operator. In \cite{Witten:1985} Witten argued that the Dirac-Ramond operator can be considered as a Dirac-like operator on $LM,$ the loop space of $M.$ Thus, in searching for an analogous result for string theory, one is led to study principal bundles over $LM.$ This is the subject of \cite{Killingback:1987}. We shall briefly outline Killingback's argument here. Denote by $LX$ the loop space of $X,$ that is, the set of smooth maps from the circle into $X,$ $\Map (S^1, X).$ Consider a principal $G$-bundle $Q \to M$ (for $G$ a compact, simple, simply-connected Lie group). Then by considering the associated loop spaces, we obtain a principal $LG$-bundle\footnote{For the proof that this in in fact a Fr\'echet principal bundle, see \cite{Carey:1991}} $LQ \to LM,$  We shall call such a bundle a \emph{loop bundle}. In the case that $X = G,$ we have the loop group of $G$ which has been extensively studied (see for example \cite{Pressley-Segal}). There is an extension of this group by the circle $S^1,$
$$
0 \to S^1 \to \widehat{LG} \to LG \to 0.
$$
This extension is central in the sense that the image of $S^1$ in $\widehat{LG}$ is in the centre of $\widehat{LG}.$ We shall look more closely at this central extension later. For now, let us just outline Killingback's result. Killingback considers, as the analogue of a spin structure for string theory, a lifting of the $LG$-bundle $LQ$ to a principal $\widehat{LG}$-bundle $\widehat{LQ}.$ The exact sequence above leads to an exact sequence of sheaves of groups over $LM.$ That is,
$$
\underline{S}^1 \to \underline{\widehat{LG}} \to \underline{LG},
$$
where $\underline{\cG}$ is the sheaf of $\cG$-valued functions over $LM.$ In general, if we have a short exact sequence of sheaves of abelian groups over $X$
$$
\underline{A} \to \underline{B} \to \underline{C},
$$
then this leads to a long exact sequence of sheaf cohomology groups (see \cite{Brylinski:1993})
$$
\cdots \to H^n(X, \underline{A}) \to H^n(X, \underline{B}) \to H^n(X, \underline{C})
 \to H^{n+1}(X, \underline{A}) \to \cdots
$$
The same is not true, however, in the nonabelian case since we cannot define the cohomology groups $H^j (X, \underline{A})$ for $j >1.$ Indeed, if $A, B$ and $C$ are nonabelian, then $H^1(X, \underline{A}),$ $H^1(X, \underline{B})$ and $H^1(X, \underline{C})$ are not groups but pointed sets. In this case, we can write down an exact sequence of pointed sets
$$
0 \to H^0(X, \underline{A}) \to H^0(X, \underline{B}) \to H^0(X, \underline{C}) \to H^1(X, \underline{A}) \to H^1(X, \underline{B}) \to H^1(X, \underline{C}),
$$
where by exactness here we mean the image of any map is exactly the pre-image of the basepoint in the next set in the sequence. There is no connecting homomorphism $H^1(X, \underline{C}) \to H^2(X, \underline{A})$ and so the sequence terminates. If we assume that $A$ is central in $B,$ however, then $H^j(X, \underline{A})$ is an abelian group for all $j$ and it is possible to extend the sequence above one more step to the right (\cite{Brylinski:1993}, Theorem 4.1.4)
\begin{multline*}
0 \to H^0(X, \underline{A}) \to H^0(X, \underline{B}) \to H^0(X, \underline{C})\\
 \to H^1(X, \underline{A}) \to H^1(X, \underline{B}) \to H^1(X, \underline{C}) \to H^2(X, \underline{A}).
\end{multline*}
The short exact sequence above therefore leads to an exact sequence in sheaf cohomology
$$
\ldots \to H^1(LM, \underline{S}^1) \to H^1(LM, \underline{\widehat{LG}}) \to H^1(LM, \underline{LG}) \to H^2(LM, \underline{S}^1),
$$
where, since $\widehat{LG}$ and $LG$ are in general nonabelian, $ H^1(LM, {\widehat{LG}})$ and $H^1(LM, \underline{LG})$ are just pointed sets, whereas $H^1(LM, \underline{S}^1)$ and $H^2(LM, \underline{S}^1)$ are abelian groups. Now, since the set of isomorphism classes of principal $\cG$-bundles over $LM$ is in bijective correspondence with the set $H^1(LM, \underline{\cG})$ we see that the $LG$-bundle $LQ \in H^1(LM, \underline{LG})$ has a lift to an $\widehat{LG}$-bundle exactly when $LQ$ is the image of an element in $H^1(LM, \underline{\widehat{LG}}).$ That is, when the image of $LQ$ in $H^2(LM, \underline{S}^1)$ is zero. Therefore, the obstruction to lifting a loop bundle $LQ \to LM$ is a class in $H^2(LM, \underline{S}^1).$ Now recall that the short exact sequence of groups
$$
0 \to \ZZ \to \RR \to S^1 \to 0,
$$
leads to an exact sequence of sheaves (as above)
$$
\underline{\ZZ} \to \underline{\RR} \to \underline{S}^1,
$$
which in turn leads to a long exact sequence of sheaf cohomology groups
$$
\ldots \to H^2(LM, \underline{\ZZ}) \to H^2(LM, \underline{\RR}) \to H^2(LM, \underline{S}^1) \to H^3(LM, \underline{\ZZ}) \to \ldots
$$
(since $\ZZ, \RR$ and $S^1$ are all abelian). However, because $\underline{\RR}$ is a soft sheaf, $H^*(LM, \underline{\RR}) = 0$ and we have the following well known result (see for example \cite{Brylinski:1993})
$$
H^2(LM, \underline{S}^1) \simeq H^3(LM, \ZZ).
$$
So we see that the obstruction to lifting the $LG$-bundle $LQ$ to an $\widehat{LG}$-bundle is a class in $H^3(LM, \ZZ).$ Since this lifting is the analogue in string theory of a spin structure for $M,$ we call it a \emph{string structure} for $M$ and we call the obstruction class $s(LQ) \in H^3(LM, \ZZ)$ the \emph{string class}. Killingback's main result, then, is a characterisation of this class in terms of the first Pontrjagyn class of the $G$-bundle $Q \to M.$ In particular, if $p_1(Q) \in H^4(M, \ZZ)$ is the first Pontrjagyn class of $Q$, then Killingback shows that the transgression of this is the string class of $LQ.$ That is, the string class is given by pulling-back $p_1(Q)$ by the evaluation map $\ev \colon LM \times S^1 \to M$ to give a class on $LM \times S^1$ and integrating over $S^1:$
$$
s(LQ) = \int_{S^1} \ev^*p_1(Q).
$$
We shall give a proof of this formula later (in section \ref{S:Higgs fields, LG-bundles,...}) following the methods in \cite{Murray:2003}.

\section{Bundle gerbes}\label{S:bundle gerbes}
In order to perform calculations involving the string class and to extend Killingback's result, we shall use the theory of bundle gerbes \cite{Murray:1996}, in particular, the lifting bundle gerbe (see section \ref{S:central extensions}). In this section we briefly outline the theory (developed largely in \cite{Murray:1996} and \cite{Murray:2000}) behind these objects. Bundle gerbes can be considered, in some sense, as `higher' versions of $U(1)$-bundles. Therefore, we start with some basic results on these bundles before describing the theory of bundle gerbes.

\subsection{$U(1)$-bundles}\label{SS:U(1)-bundles}

As mentioned, we shall begin by recalling some facts about $U(1)$-bundles and some constructions involving these bundles. Firstly, note that if $P \to M$ is a $U(1)$-bundle with right action given by $(p, z) \mapsto pz$ (for $p\in P$ and $z \in U(1)$) then there is a \emph{dual} bundle, denoted $P^*,$ which is the same as $P$ but with the action given by $(p, z) \mapsto p z^{-1}.$ Of course this is only a right action because $U(1)$ is abelian. Further, if $Q$ is another $U(1)$-bundle over $M,$ we can form the fibre product over $M,$ $P\times_M Q,$ which is a principal $U(1) \times U(1)$-bundle over $M$ whose fibres are the product of the fibres of $P$ and $Q$ (i.e.\! $(P\times_M Q)_m = P_m \times Q_m$). By factoring out by the `anti-diagonal' inside $U(1) \times U(1),$ that is, the set $\{(z, z^{-1})\},$ we obtain a principal $U(1)$-bundle called the \emph{contracted product} of $P$ and $Q$ and denoted $P\otimes Q.$ It is easy to see that $P\otimes P^*$ is canonically trivialised by the section $s \colon m \mapsto [p, p^*],$ where $p$ is any point in the fibre of $P$ above $m$ and $p^*$ is the same point considered as an element of $P^*.$ For if $s_\a$ and $s_\b$ are two such local sections then suppose $s_\a (m) = [p, p^*]$ and $ s_\b (m) = [q,q^*],$ then we have that $[q, q^*] = [pz, p^*z^{-1}] $ for some $z \in U(1)$ and so $s_\a = s_\b.$

Note that if instead of considering $U(1)$-bundles we equivalently considered complex hermitian line bundles then the dual would correspond to the linear dual of a line bundle (i.e.\! the bundle whose fibres are the dual of those of the original bundle) and the contracted product would correspond to the tensor product of line bundles (the bundle whose fibres are the tensor product of the fibres of the original two bundles). Note also that if $P$ and $Q$ have transition functions $g_{\a \b}$ and $h_{\a \b}$ respectively relative to some open cover of $M$ then $P^*$ has transition functions $g_{\a \b}^{-1}$ and $P\otimes Q$ has transition functions $g_{\a \b} h_{\a \b}.$

Another important property of $U(1)$-bundles on $M$ is the way in which they relate to $H^2(M, \ZZ).$ If a $U(1)$-bundle $P$ has transition functions $g_{\a \b}$ then on triple overlaps these satisfy the cocycle condition $g_{\b \gamma}^{\phantom{-1}} g_{\a \gamma}^{-1} g_{\a \b}^{\phantom{-1}} = 1$ and thus form a class in $H^1 (M,\underline{U(1)}).$ Thus, from the argument in the previous section we have that a $U(1)$-bundle defines a class in $H^2(M, \ZZ).$ This class is called the \emph{Chern class} of the bundle $P.$ It is a standard result (see for example \cite{Brylinski:1993}) that the Chern class classifies $U(1)$-bundles up to isomorphism and, further, that given any class in $H^2(M, \ZZ)$ one can construct a $U(1)$-bundle. So we see that isomorphism classes of $U(1)$-bundles are in bijective correspondence with $H^2(M, \ZZ).$ The Chern class is additive in the sense that if $c(P)$ and $c(Q)$ are the Chern classes of $P$ and $Q$ respectively, then $c(P\otimes Q) = c(P) + c(Q)$ and $c(P^*) = -c(P).$ It is natural in the sense that if we pull-back the bundle $P \to M$ by a map $f \colon N \to M$ to give a $U(1)$-bundle $f^*P \to N$ then $c(f^*P) = f^*c(P).$

We can actually represent the image of the Chern class in real cohomology using differential forms quite easily. If $A$ is a connection on $P$ whose curvature is $F,$ then $F/ 2\pi i $ is a closed integral form and its class in the de Rham cohomology group $H^2(M)$ is the image in real cohomology of the Chern class of $P.$

\subsection{Bundle gerbes}\label{SS:bundle gerbes}

\subsubsection{Definitions and basic constructions}

Having reviewed some of the basic properties of $U(1)$-bundles in the previous section, we would now like to present another object, first introduced in \cite{Murray:1996} and studied further in \cite{Murray:2000}, which is in some sense a higher dimensional version of a $U(1)$-bundle as we shall see shortly.

Consider a surjective submersion $Y \xrightarrow{\pi} M.$ We can form the fibre product of $Y$ with itself, which we denote $Y^{[2]},$ and we have (as before)
$$
Y^{[2]} = \{ (y_1, y_2) \in Y\times Y \mid \pi(y_1) = \pi(y_2)\}.
$$
Note that since $\pi$ is a submersion $Y^{[2]}$ is a submanifold of $Y^2$. In general we have the \emph{p-fold fibre product} $Y^{[p]}$ defined similarly. We define the maps $\pi_i  \colon Y^{[p+1]} \to Y^{[p]} (i=1, \ldots, p+1)$ to be omission of the $i^\text{th}$ factor,
$$
\pi_i (y_1, \ldots, y_{p+1}) = (y_1, \ldots, y_{i-1}, y_{i+1}, \ldots, y_{p+1}).
$$
We have, then, the following definition:
\begin{definition}[\cite{Murray:1996}]
A \emph{bundle gerbe} over a manifold $M$ is a pair $(P, Y)$ where $Y \to M$ is a surjective submersion and $P \to Y^{[2]}$ is a $U(1)$-bundle and such that there is a \emph{bundle gerbe multiplication}, which is a smooth isomorphism
$$
m \colon \pi_3^* P \otimes \pi_1^* P \xrightarrow{\sim} \pi_2^* P
$$
of $U(1)$-bundles over $Y^{[3]}.$ Further, this multiplication is required to be associative whenever triple products are defined. That is, if $P_{(y_1, y_2)}$ denotes the fibre of $P$ over $(y_1, y_2) \in Y^{[2]}$ then the following diagram commutes for all $(y_1, y_2, y_3, y_4) \in Y^{[4]}$:
$$
\xymatrix@C=6.5ex{P_{(y_1, y_2)} \otimes P_{(y_2, y_3)} \otimes P_{(y_3, y_4)} \ar[d]_{\text{id} \otimes m} \ar[r]^-{m \otimes \text{id}} & P_{(y_1, y_3)}\otimes P_{(y_3, y_4)} \ar[d]^{m}\\
		P_{(y_1, y_2)} \otimes P_{(y_2, y_4)} \ar[r]_-{m} & P_{(y_1, y_4)}}
$$
We sometimes denote a bundle gerbe simply by $P.$
\end{definition}

We typically depict a bundle gerbe thusly:
$$
\xymatrix{P\ar[d]\\
		Y^{[2]} \ar@<0.3ex>[r] \ar@<-0.7ex>[r] &Y^{\vphantom{[2]}}\ar[d]\\
		&M}
$$

We can characterise the bundle gerbe multiplication and its associativity in a different way using sections of bundles related to $P$ as follows. If $Q\to Y^{[p-1]}$ is a $U(1)$-bundle, define the bundle $\d Q \to Y^{[p]}$ as
$$
\d Q = \pi_1^*Q \otimes (\pi_2^* Q )^* \otimes \pi_3^* Q \otimes \ldots
$$
Then it is easy to show that $\d\d Q$ is canonically trivial. One can show that the bundle gerbe multiplication is equivalent to a section $s$ of $\d P \to Y^{[3]}$ and that the associativity condition is equivalent to the condition that $\d s = 1$ as a section of $\d\d P$ (where $1$ denoted the canonical section of $\d\d P$). Indeed if $p$ and $q$ are elements of $P_{(y_1, y_2)}$ and $P_{(y_2, y_3)}$ respectively, we can define a section $s$ of $\d P$ by
$$
s(y_1, y_2, y_3) = p\otimes m(p,q)^*\otimes q,
$$
then the associativity of $m$ forces the condition $\d s =1.$ Note that these conditions reflect the definition of a \emph{simplicial line bundle} from \cite{Brylinski:1994}. So we see that a bundle gerbe is the same as a simplicial line bundle over the simplicial space defined by the fibre products $Y^{[p]}.$ We shall discuss simplicial spaces and this relationship more in section \ref{S:central extensions}.

In \cite{Murray:1996} Murray claimed that bundle gerbes were essentially bundles of groupoids. Although it is not essential for our purposes let us briefly explain what is meant by this. Recall (see \cite{MacLane:1971}) that a groupoid is a small category with all arrows invertible. Consider then a bundle gerbe $(P,Y)$ over $M.$ If we consider the elements of the fibre over $m,$ $Y_m,$ as the objects of a category, then the elements of the fibre $P_{(y_1, y_2)}$ are the morphisms from $y_1$ to $y_2$ and the bundle gerbe multiplication gives a way of composing these morphisms. Since $P _{(y_1, y_2)}^{\vphantom{*}} \simeq P_{(y_2, y_1)}^*$ and $P_{(y,y)} \simeq Y^{[2]} \times U(1)$ (which can be shown using the bundle gerbe multiplication), this category is a groupoid. In \cite{Murray:1996} the theory of $U(1)$-groupoids is presented in more detail as a prelude to the introduction of bundle gerbes.

Just as for $U(1)$-bundles, various constructions are possible with bundle gerbes \cite{Murray:1996}. Consider a map $f \colon N \to M.$ We can pull-back the submersion $Y\to M$ to a submersion $f^*Y \to N.$ This gives a map $\hat{f} \colon f^*Y \to Y$ covering $f$ which induces a map (also called $\hat{f}$) $(f^*Y)^{[2]} \to Y^{[2]}.$ Thus we can pull-back the $U(1)$-bundle $P \to Y^{[2]}$ by $\hat{f}$ to give a bundle $\hat{f}^*P \to (f^*Y)^{[2]}.$ So we have a bundle gerbe over $N$ called the \emph{pull-back} and which we will denote $f^*P.$ We can also define the \emph{dual} of $(P, Y)$ by taking the dual of the $U(1)$-bundle $P$ over $Y^{[2]}.$ We denote this by $P^*.$ We can form the \emph{product} of two bundle gerbes $(P, Y)$ and $(Q, X)$ over $M,$ denoted $P\otimes Q,$ by taking the fibre product $Y\times_M X$ over $M$ and the $U(1)$-bundle $P\otimes Q$ over $(Y\times_M X)^{[2]}.$

We say two bundle gerbes $(P, Y)$ and $(Q, X)$ over $M$ are \emph{isomorphic} if there is an isomorphism $Y\to X$ covering the identity on $M$ and a bundle isomorphism $P \to Q$ covering the induced map $Y^{[2]} \to X^{[2]}$ and which commutes with the bundle gerbe multiplication.

A particular example of a bundle gerbe is given by taking a $U(1)$-bundle $P$ over $Y$ and defining $\d P $ over $Y^{[2]}$ as above. That is, $\d P = \pi_1^*P \otimes (\pi_2^* P)^*.$ Since $\d\d P$ is canonically trivial over $Y^{[3]},$ it has a canonical section $s$ which defines the bundle gerbe multiplication. This is called the \emph{trivial bundle gerbe} and in general we say a bundle gerbe is \emph{trivial} if it is isomorphic to one of this form.

As was pointed out in \cite{Murray:2000} there is another notion of equivalence, in addition to isomorphism, for bundle gerbes. This is the notion of \emph{stable isomorphism}, first introduced in \cite{Carey:1997} and studied in detail in \cite{Murray:2000}. Two bundle gerbes $(P, Y)$ and $(Q, X)$ are called \emph{stably isomorphic} if there are trivial bundle gerbes $T_1$ and $T_2$ such that $P\otimes T_1 \simeq Q\otimes T_2$ or, equivalently, if $P\otimes Q^*$ is trivial. It turns out that stable isomorphism is in some sense the correct notion of equivalence for bundle gerbes because, as we shall see next, all bundle gerbes have a characteristic class associated to them and this class classifies them up to stable isomorphism. That is, two bundle gerbes have the same associated class exactly when they are stably isomorphic. This class is called the \emph{Dixmier-Douady class} and it is to this which we now turn our attention.

\subsubsection{Bundle gerbes and degree three cohomology}

As mentioned earlier, bundle gerbes can be considered as higher dimensional $U(1)$-bundles. We now explain why this is the case and describe how to construct a characteristic class for bundle gerbes which is analogous to the Chern class for $U(1)$-bundles.

Let $(P, Y)$ be a bundle gerbe over $M$ and choose a good cover $\{U_\a \}$ of $M$ over which $Y \to M$ admits local sections. This is always possible (see \cite{Bott-Tu}). Suppose that $s_\a \colon U_\a \to Y$ is a local section. We have a section of $Y^{[2]}$ over double overlaps given by
$$
(s_\a, s_\b) \colon U_{\a \b} \to Y^{[2]},
$$
where $U_{\a \b} = U_\a \cap U_\b.$ As $U_{\a\b}$ is contractible, the pull-back $P_{\a \b} = (s_\a, s_\b)^*P$ of $P$ by this section is trivial. The fibres of $P_{\a\b}$ are given by $(P_{\a\b})_m = P_{(s_\a (m), s_\b (m))}.$ Choose a section $\sigma_{\a \b}$ of this bundle. That is, a map
$$
\sigma_{\a\b} \colon U_{\a\b} \to P
$$
such that $\sigma_{\a\b}(m) \in P_{(s_\a (m), s_\b (m))}.$ On triple overlaps $U_{\a\b\gamma}$ the bundle gerbe multiplication gives
$$
m(\sigma_{\a\b}, \sigma_{\b\gamma}) = g_{\a\b\gamma} \sigma_{\a\gamma}
$$
for some $g_{\a\b\gamma} \colon U_{\a\b\gamma} \to U(1).$ On overlaps $U_{\a \b \gamma \d}$ the associativity of this multiplication gives the cocycle condition
$$
g_{\b \gamma \d}^{\vphantom{-1}} g_{\a \gamma \d}^{-1} g_{\a \b \d}^{\vphantom{-1}} g_{\a \b \gamma}^{-1} = 1.
$$
Thus the functions $g_{\a\b\gamma}$ define a class in $H^2(M, \underline{U(1)}) \simeq H^3(M, \ZZ).$ This class is independent of any choices and is called the \emph{Dixmier-Douady class} of $P$ and denoted $\DD(P).$ In \cite{Murray:1996} it is proven that this class is precisely the obstruction to the bundle gerbe being trivial. We also have the following results regarding the Dixmier-Douady class for the constructions presented earlier: If $(P, Y)$ and $(Q, X)$ are bundle gerbes over $M$ then $\DD(P\otimes Q) = \DD(P) + \DD(Q)$ and $\DD(P^*) = -\DD(P).$ The Dixmier-Douady class is natural with respect to pull-backs, that is, $\DD(f^*P) = f^*\DD(P).$

As mentioned at the end of the previous section, the Dixmier-Douady class classifies bundle gerbes up to stable isomorphism. This is clear because $P$ and $Q$ are stably isomorphic exactly when $P\otimes Q^*$ is trivial and so the result follows from the fact that $\DD(P\otimes Q^*) = \DD(P) - \DD(Q)$ and that trivial bundle gerbes have zero Dixmier-Douady class.

In \cite{Murray:1996} it is also shown that every class in $H^3(M, \ZZ)$ is the Dixmier-Douady class of some bundle gerbe. This means that there is a bijection between $H^3(M, \ZZ)$ and stable isomorphism classes of bundle gerbes. Thus bundle gerbes provide a geometric realisation of elements in $H^3(M, \ZZ)$ in an analogous way to that of $U(1)$-bundles and $H^2(M, \ZZ).$

\subsubsection{Connective structures on bundle gerbes}

We have seen now the way in which bundle gerbes play a role for degree three cohomology analogous to that of $U(1)$-bundles and degree two cohomology. As we saw in section \ref{SS:U(1)-bundles} $U(1)$-bundles have the nice property that the image of their Chern class in real cohomology is represented by the form $F/2\pi i,$ where $F$ is the curvature of the bundle. We would now like to study connective structures on bundle gerbes and, as we shall see, a similar result is true in this case.

Consider first the $p$-fold fibre product $Y^{[p]}$ as before. Let $\Omega^q(Y^{[p]})$ denote the space of differential $q$-forms on $Y^{[p]}.$ Then we can define a map $\d \colon \Omega^q(Y^{[p]}) \to \Omega^q(Y^{[p+1]})$ as the alternating sum of pull-backs by the projections $\pi_i:$
$$
\d = \sum_{i=1}^{p+1} (-1)^{i-1} \pi_i^*.
$$
Then $\d^2 = 0$ and so we have a complex
$$
0 \to \Omega^q(M) \xrightarrow{\pi^*} \Omega^q(Y) \xrightarrow{\,\d\,} \Omega^q(Y^{[2]}) \xrightarrow{\,\d\,} \Omega^q(Y^{[3]}) \xrightarrow{\,\d\,} \ldots
$$
In \cite{Murray:1996} it is proven that this complex has no cohomology. That is, the above sequence is exact for all $q \geq 0.$ We shall use this result shortly.

A bundle gerbe connection is a connection $A$ for the $U(1)$-bundle $P$ that respects the bundle gerbe product in the sense that the induced connection on $\pi_2^* P$ is the same as the image of the induced connection on $\pi_3^*P \otimes \pi_1^*P$ under the bundle gerbe multiplication. Note that if $s \colon Y^{[3]} \to \d P$ is the section defining this multiplication, then this means that a bundle gerbe connection satisfies $s^*(\d A) = 0.$ That is, $\d A$ is flat with respect to $s.$ Using this observation, it is easy to see that bundle gerbe connections always exist. For consider a connection $A$ on $P$ that does not necessarily commute with the product. We cannot say that $s^*(\d A) = 0$ but note that $\d ( s^*(\d A)) = (\d s)^* (\d\d A),$ which is zero since $\d s =1$ as a section of $\d\d P$ and $\d\d A$ is flat with respect to the canonical trivialisation of $\d\d P.$ Therefore, by the exact sequence above there is some $a \in \Omega^1(Y^{[2]})$ such that $\d a = s^*(\d A)$ and so $s^*(\d (A - \pi^*a)) = 0$ (where $\pi \colon P \to Y^{[2]}$ is the projection). Therefore, $A - \pi^*a$ is a bundle gerbe connection.

If $F$ is the curvature of a bundle gerbe connection $A$ viewed as a 2-form on $Y^{[2]},$ then $\d F = s^* (\d dA) = d( s^* (\d A)) = 0.$ This means that there is some $B \in \Omega^2(Y)$ satisfying $F = \d B.$ A choice of such a $B$ is called a \emph{curving} for $P.$ Note that if $B'$ is another choice of curving then $B$ and $B'$ differ by a $\d$-closed (and hence $\d$-exact) 2-form on $Y$. As $\d$ and $d$ commute, we have that $\d (dB) = d(\d B) = dF = 0.$ Therefore there is a 3-form $H$ on $M$ such that $dB = \pi^* H$ (for $\pi$ the projection $Y \to M$). $H$ is called the \emph{3-curvature} of $P.$ It is closed and a different choice of $B$ or $H$ would result in a difference of an exact form. So $H$ defines a cohomology class in $H^3(M).$ It turns out that the 3-form $H/2\pi i$ is integral and that $H/ 2\pi i$ is a representative of the Dixmier-Douady class of $P$ in real cohomology.

\section{Central extensions and the lifting bundle gerbe}\label{S:central extensions}


In this thesis, we wish to apply the theory of bundle gerbes to the study of central extensions of Lie groups and, in particular, to lifting problems as in section \ref{S:string structures}. For this purpose we use a particular bundle gerbe called the \emph{lifting bundle gerbe} and in this section we review the basic definitions and results required to develop the theory. We shall start by outlining the theory of central extensions, following \cite{Brylinski:1994}.

\subsection{Simplicial line bundles and central extensions}\label{SS:slb's and ce's}

We begin by recalling some simplicial techniques. Recall (see \cite{Dupont:1978}) that a \emph{simplicial space} is a collection of spaces $\{X_p\}\, (p= 0, 1, 2,\ldots)$ together with maps $d_i \colon X_p \to X_{p-1}$ and $s_j\colon X_p \to X_{p+1}$ for $i,j = 0, \ldots, p$, called \emph{face} and \emph{degeneracy} maps respectively, which satisfy the simplicial identities
\begin{align*}
	d_i d_j	&= d_{j-1}d_i,	\qquad i < j,\\
	s_i s_j	&= s_{j+1} s_i,	\qquad i\leq j,\\
	d_i s_j	&= \begin{cases}
				s_{j-1} d_i, 	& i < j\\
				\text{id}, 		& i=j, \, j+1\\
				s_j d_{i-1},	& i > j+1.
			\end{cases}
\end{align*}
If we are working in the category of manifolds and smooth maps we say that $\{X_p\}$ is a \emph{simplicial manifold}. For example, consider the collection\footnote{Note that here $X_0 = Y, X_1 = Y^{[2]}, X_2 = Y^{[3]}, \ldots$ and so on} $\{Y^{[p+1]}\}$ of fibre products as in the previous section. These form a simplicial manifold with the obvious face and degeneracy maps. Note that for a general simplicial manifold $\{ X_p\}$ we can define a complex similar to the one described in section \ref{S:bundle gerbes} by using the pull-backs of the face maps $d_i.$ That is, we define $\d\colon \Omega^q(X_p) \to \Omega^q(X_{p+1})$ by
$$
\d = \sum_{i=0}^p (-1)^i d_i^*.
$$
Also, as before, if $Q$ is a $U(1)$-bundle (or an hermitian line bundle) over $X_p$ then we can define a bundle over $X_{p+1}$ by
$$
\d Q = d_0^*Q \otimes (d_1^* Q)^* \otimes d_2^*Q \otimes \ldots
$$
The particular example of interest to us is a certain simplicial manifold associated to a Lie group which we describe presently. Let $\cG$ be a Lie group. There is a simplicial manifold called $N\cG = \{ N\cG_p \}$ given by the manifolds $\{ \cG^p\}$ and face and degeneracy maps $d_i$ and $s_j$ where
$$
d_i (g_1, \ldots, g_{p+1}) = \begin{cases}
					(g_2,\ldots, g_{p+1}), 	& i=0\\
					(g_1, \ldots, g_{i-1}g_i, g_{i+1}, \ldots, g_{p+1}), &1\leq i \leq p-1\\
					(g_1, \ldots, g_p),		& i=p
					\end{cases}
$$
and
$$
s_j(g_1, \ldots, g_{p+1}) = (g_1, \ldots, g_{j-1}, 1, g_{j}, \ldots, g_{p+1}).
$$
We would like to consider central extensions of $\cG$ by the circle and show how they are related to $N\cG.$ For this, we follow Brylinski and McLaughlin \cite{Brylinski:1994} where the result is phrased in terms of simplicial line bundles. We have the following definition
\begin{definition}[\cite{Brylinski:1994}]
Let $\{X_p\}$ be a simplicial manifold. A \emph{simplicial line bundle} over $\{X_p\}$ is a line bundle $L$ on $X_1$ together with a section $s$ of the bundle $\d L \to X_2$ such that $\d s = 1$ as a section of $\d\d L.$
\end{definition}
Notice the similarity with the definition of a bundle gerbe. In fact, instead of using $U(1)$-bundles, we can rephrase everything about bundle gerbes in terms of line bundles and we see that a bundle gerbe is the same thing as a simplicial line bundle over the simplicial space $\{Y^{[p]}\}.$

Now consider a central extension of $\cG$ by the circle
$$
U(1) \to \widehat{\cG} \xrightarrow{p} \cG.
$$
If we think of this as a $U(1)$-bundle $\widehat{\cG} \to \cG$ then we must have a multiplication $M \colon \widehat{\cG} \times\widehat{\cG} \to \widehat{\cG}$ which covers the multiplication on $\cG,$ that is, $m = d_1 \colon \cG\times \cG \to \cG.$ Because $\widehat{\cG}$ is a central extension we must have $M(\hat{g}z, \hat{h}w) = M(\hat{g}, \hat{h})(zw)$ for any $\hat{g}, \hat{h} \in \widehat{\cG}$ and $z,w \in U(1).$ In a similar way to that in which the bundle gerbe multiplication on a bundle gerbe $P$ gave rise to a section of $\d P,$ this gives a section of $\d \widehat{\cG},$
$$
s(g, h) = \hat{g} \otimes M(\hat{g}, \hat{h})^* \otimes \hat{h},
$$
where $\hat{g}$ and $\hat{h}$ are points in the fibres over $g$ and $h$ respectively. The associativity of this multiplication is equivalent to the condition $\d s = 1$ as before and hence a central extension gives rise to a simplicial line bundle. In fact it can be shown that they are equivalent and we have the result from \cite{Brylinski:1994}:
\begin{theorem}[\cite{Brylinski:1994}]
A simplicial line bundle over the simplicial manifold $N\cG$ is a central extension of $\cG$ by the circle.
\end{theorem}
We wish to perform explicit calculations using differential forms so, following \cite{Murray:2001} and \cite{Murray:2003}, we shall rephrase this result in terms of differential forms on $\cG^p$ and give a method of constructing central extensions using these forms. Consider then, a connection $\nu$ for $\widehat{\cG}$ thought of as a $U(1)$-bundle over $\cG.$ As in the treatment of bundle gerbe connections in section \ref{S:bundle gerbes} we can consider the induced connection $\d \nu$ on the bundle $\d \widehat{\cG} \to \cG \times \cG$ and then, as this bundle is trivial, we can pull-back $\d \nu$ by the section $s.$ Let $\a = s^*(\d \nu).$ In general $\a$ is non-zero. However, we have that $\d \a = \d (s^*(\d \nu)) = (\d s)^*(\d \d \nu) = 0.$ Furthermore, we also have $d \a = s^*(d \d \nu) = \d R,$ where $R$ is the curvature of $\nu$ viewed as a form on $\cG.$ Therefore we have constructed from the central extension a pair of forms $(R, \a),$ where $R \in \Omega^2(\cG)$ is closed and integral and $\a \in \Omega^1(\cG \times \cG)$ is such that $\d R = d \a$ and $\d \a = 0.$ In fact, as we shall now show, this pair is sufficient to reconstruct the central extension. Recall (see for example \cite{Brylinski:1993}) that given an integral 2-form $R\in \Omega^2(\cG)$ there exists a principal $U(1)$-bundle $P\to \cG$ with a connection $a$ whose curvature is $R.$ 
Also, it is a standard result (see \cite{Kobayashi:1963}) that if $Q$ is a bundle over a simply connected base which admits a flat connection $A,$ then $Q$ is trivial and there is a section $s$ of $Q$ such that $s^*A = 0.$ In terms of the construction here, this means that we can find a bundle $P \to \cG$ with curvature $R$ and because $d \a = \d R,$ we have that $\d a - \pi^*\a$ is a flat connection on $\d P \to \cG \times \cG.$ Therefore, there is a section $s$ of $\d P$ satisfying $s^*(\d a) = \a.$ As before, this section defines a multiplication and we can calculate $\d s$ which we want to be equal to $1.$ Now, $(\d s)^*(\d \d a) = \d (s^* (\d a)) = \d \a = 0$ and for the canonical section $1$ we also have $1^*(\d\d a) = 0.$ This means that they differ by an element of $U(1)$ and so rather than associativity of the multiplication $M$ defined by $s$ we have
$$
M(M(\hat{g}, \hat{h}), \hat{k}) = z M(\hat{g}, M(\hat{h}, \hat{k}))
$$ 
for some $z \in U(1).$ However, if we choose some $\hat{g}$ in the fibre above the identity $e$ in $\cG$ then $M(\hat{g}, \hat{g})$ is also in the fibre above $e$ and so $\hat{g}$ and $M(\hat{g}, \hat{g})$ differ by some $w \in U(1).$ That is, $M(\hat{g}, \hat{g}) = \hat{g} w.$ Let $\hat{h}$ and $\hat{k}$ both be equal to $\hat{g} \in \pi^{-1}(e).$ Then the formula above reads
$$
M(M(\hat{g}, \hat{g}), \hat{g}) = z M(\hat{g}, M(\hat{g}, \hat{g}))
$$
and so $\hat{g} w^2 = \hat{g} w^2 z$ and we see that in fact $z=1.$

Thus we have constructed a central extension from the pair $(R,\a)$ and this construction recovers the original extension (which follows from the fact that $P$ has curvature $R$ and the definition of $\a$ above). Note that isomorphic central extensions (where by isomorphic, we mean isomorphic as $U(1)$-bundles and as groups) give rise to the same $R$ and $\a$ and that in constructing the pair $(R, \a)$ if we had chosen a different connection, by adding on the pull-back of a 1-form $\eta$ on $\cG,$ then we would have the pair $(R + d\eta, \a + \d \eta).$ Also, note that the section constructed above from the flat connection is not unique but changing this by multiplying by a constant $z$ in $U(1)$ would change $M$ to $Mz$ and, as the extension is central, this would give an isomorphic central extension. So, as in \cite{Murray:2003}, we have a bijection between isomorphism classes of central extensions with connection and pairs of forms satisfying the conditions above.

\subsection{The lifting bundle gerbe}

Having reviewed a method for constructing central extensions, we would like now to link the theory of central extensions with that presented earlier on bundle gerbes. We present a particular example of a bundle gerbe related to central extensions, first introduced in \cite{Murray:1996}, called the \emph{lifting bundle gerbe} whose Dixmier-Douady class is precisely the obstruction to lifting a $\cG$-bundle $P$ to a $\widehat{\cG}$-bundle $\widehat{P}.$

Consider then a principal $\cG$-bundle $P \to M.$ Choose a good cover of $M$ and consider the transition functions $g_{\a\b}$ of $P$ relative to this cover. We can choose lifts of these functions $\hat{g}_{\a\b}$ which take values in $\widehat{\cG}$ and these are candidates for the transition functions of the lift $\widehat{P}.$ However, transition functions are required to satisfy the cocycle condition $g_{\a\b}  g_{\b\gamma} = g_{\a\gamma}$ on triple overlaps but the lifts $\hat{g}_{\a\b}$ only satisfy
$$
\hat{g}_{\a\b}  \hat{g}_{\b\gamma} = \epsilon_{\a\b\gamma}\hat{g}_{\a\gamma}
$$
for some $U(1)$-valued function $\epsilon_{\a\b\gamma}.$ This means that the $\hat{g}_{\a\b}$'s are not necessarily transition functions. However, due to the fact that $\widehat{\cG}$ is a central extension, it can be shown that the functions $\epsilon_{\a\b\gamma}$ satisfy the cocycle condition
$$
\epsilon_{\b\gamma\d}^{\vphantom{-1}} \epsilon_{\a\gamma\d}^{{-1}} \epsilon_{\a\b\d}^{\vphantom{-1}} \epsilon_{\a\b\gamma}^{{-1}}  = 1.
$$
Therefore, $\epsilon_{\a\b\gamma}$ defines a class in $H^2(M, \underline{U(1)}) \simeq H^3(M,\ZZ).$ As per the discussion in section \ref{S:string structures}, this class is the obstruction to lifting the transition functions $g_{\a\b}$ to transition functions $\hat{g}_{\a\b}$ and hence the obstruction to lifting $P$ to $\widehat{P}.$

If we take the principal $\cG$-bundle $P \to M$ and consider the fibre product $P^{[2]} \rightrightarrows P$ then there is a natural map $\tau \colon P^{[2]} \to \cG,$ called the \emph{difference map}, given by $p_1 \tau(p_1, p_2) = p_2.$ If we view $\widehat{\cG}$ as a $U(1)$-bundle over $\cG$ then we can pull-back $\widehat{\cG}$ by this map to obtain a $U(1)$-bundle over $P^{[2]}:$
$$
\xymatrix{ \tau^* \widehat{\cG} \ar[r] \ar[d]	& \widehat{\cG} \ar[d]\\
			P^{[2]} \ar[r]^\tau			& \cG^{\vphantom{[2]}}}
$$
where
$$
\tau^* \widehat{\cG} = \left\{ (p_1, p_2, \hat{g}) \mid p(\hat{g}) = \tau(p_1, p_2) \right\}.
$$
Note that $\tau(p_1, p_2) \tau(p_2, p_3) = \tau(p_1, p_3)$ and so, because the multiplication in $\widehat{\cG}$ covers that in $\cG,$ we have an induced map
$$
\tau^*\widehat{\cG}_{(p_1, p_2)} \otimes \tau^*\widehat{\cG}_{(p_2, p_3)} \to \tau^*\widehat{\cG}_{(p_1, p_3)}
$$
which serves as a bundle gerbe multiplication for the bundle gerbe $(\tau^*\widehat{\cG}, P)$ over $M.$ This bundle gerbe is called the \emph{lifting bundle gerbe}. We would now like to examine its Dixmier-Douady class. Recall from section \ref{S:bundle gerbes} the construction of the Dixmier-Douady class of a bundle gerbe. This involves taking sections $s_\a$ and $s_\b$ of $P$ to give a section $(s_\a, s_\b)$ of $P^{[2]}$ over $U_{\a\b}.$ We then pull-back the bundle $\tau^*\widehat{\cG}$ by $(s_\a, s_\b)$ to give a bundle $(s_\a, s_\b)^* (\tau^*\widehat{\cG}) \to U_{\a\b}.$ The Dixmier-Douady class of $\tau^*\widehat{\cG}$ is related to sections of this bundle, that is, maps $\sigma_{\a\b} \colon U_{\a\b} \to \tau^*\widehat{\cG}$ such that $\sigma(m) \in \tau^*\widehat{\cG}_{(s_\a(m), s_\b (m))}$. The bundle gerbe multiplication (which in this case is given by the multiplication in $\widehat{\cG}$) gives $\sigma_{\a\b} \sigma_{\b\gamma} = g_{\a\b\gamma} \sigma_{\a\gamma}$ for some $U(1)$-valued function $g_{\a\b\gamma}$ and the image of this in $H^3(M, \ZZ)$ is a representative for the Dimier-Douady class of $\tau^*\widehat{\cG}.$ Note at this point, however, that as $P$ is a principal $\cG$-bundle, the sections $s_\a$ and $s_\b$ are related by the transition functions $g_{\a\b}.$ That is, $s_\b = s_\a g_{\a\b}.$ This means that $(s_\a, s_\b)^* (\tau^*\widehat{\cG})$ is given by triples $(s_\a, s_\b, \hat{g})$ where $p(\hat{g}) = g_{\a\b}.$ So in fact a section $\sigma_{\a\b}$ is given by the candidate transition functions $\hat{g}_{\a\b}.$ Therefore, the sections $\sigma_{\a\b}$ satisfy
$$
\hat{g}_{\a\b} \hat{g}_{\b\gamma} = \epsilon_{\a\b\gamma}\hat{g}_{\a\gamma},
$$
or
$$
\hat{g}_{\b\gamma}^{\vphantom{-1}} \hat{g}_{\a\gamma}^{-1} \hat{g}_{\a\b}^{\vphantom{-1}} = \epsilon_{\a\b\gamma},
$$
which is precisely the relation above for the obstruction to the existence of a lift. Thus the Dixmier-Douady class of the lifting bundle gerbe $(\tau^*\widehat{\cG}, P)$ measures the obstruction to lifting the $\cG$-bundle $P$ to a $\widehat{\cG}$-bundle $\widehat{P}.$ So the lifting bundle gerbe is trivial exactly when $P$ lifts to a $\widehat{\cG}$ bundle.

In the next section we shall demonstrate how to find a representative for the obstruction class of a particular lifting problem using the methods outlined already from the theory of bundle gerbes.

\section{The string class of an $LG$-bundle}\label{S:LG string class}


Having outlined the theory of central extensions and bundle gerbes we are now in a position to extend Killingback's result to general $LG$-bundles. In this section we will review the calculations from \cite{Murray:2003} which give an explicit expression for (the image in real cohomology of) the string class of an $LG$-bundle $P \to M,$ where here we do not require $P$ to be a loop bundle as in section \ref{S:string structures}.

\subsubsection{The central extension of the loop group}

In the previous section we showed how to classify isomorphism classes of central extensions of a Lie group $\cG$ using a 2-form $R$ on $\cG$ and a 1-form $\a$ on $\cG \times \cG.$ Now suppose that $\cG = LG,$ the loop group of a compact, simple, simply connected Lie group. In this case we can give these forms explicitly, thus making it possible to perform calculations involving the central extension $\widehat{LG}$ of $LG.$

In \cite{Pressley-Segal} Pressley and Segal give a well known expression for the curvature of a connection on the central extension $\widehat{LG}.$ Namely,
$$
R = \frac{i}{4 \pi} \int_{S^1} \<\Theta, \partial \Theta \> \, d\theta,
$$
where $\Theta$ is the (left-invariant) Maurer-Cartan form on $LG,$ which is defined pointwise, $\partial$ denotes the derivative in the loop direction, that is, the derivative with respect to $\theta$ and $\<\,\, ,\, \>$ is 
an invariant inner product\footnote{We shall refer to this as the Killing form since all invariant, bilinear, symmetric forms on $\fg$ are proportional and so this is just the Killing form with a suitable normalisation.}
 on $L\fg$ (defined pointwise) normalised so the longest root has length squared equal to 2. To construct the central extension we also need a 1-form $\a$ satisfying $\d R = d \a$ and $\d \a = 0.$ In this case it is easy to find such an $\a.$ First note that $\d R = \pi_1^*R - m^* R + \pi_2^*R$ where $m$ is the multiplication in $LG$ and $\pi_i$ is the projection $LG \times LG \to LG$ which omits the $i^\text{th}$ factor. Then $\pi_i^*R$  is given by
$$
\frac{i}{4\pi} \int_{S^1} \< \pi_i^*\Theta, \partial \pi_i^* \Theta \> \, d\theta.
$$
and using the identities
$$
\partial \Theta = ad(\gamma^{-1})d(\partial\gamma^{\vphantom{-1}} \gamma^{-1}),
$$
at the point $\gamma \in LG,$ and
$$
\partial \left(ad(\gamma^{-1}) X \right) = ad(\gamma^{-1})[X, \partial \gamma \gamma^{-1}] + ad(\gamma^{-1}) \partial X,
$$
for a vector $X \in L\fg,$ we can calculate $m^*R$ to be
\begin{multline*}
\frac{i}{4\pi} \int_{S^1} \< \Theta_1, \partial \Theta_1\> + \< [\Theta_1, \Theta_1], \partial \gamma_2^{\vphantom{-1}} \gamma_2^{-1} \> + \< \Theta_1, d(\partial \gamma_2^{\vphantom{-1}} \gamma_2^{-1})\>\\
 + \< \Theta_2,\partial (ad(\gamma_2^{-1}) \Theta_1 ) \> + \< \Theta_2, \partial \Theta_2\> \, d\theta,
\end{multline*}
where we have written $\Theta_1$ for $\pi_2^*\Theta$ and so on. So
$$
\d R = -\frac{i}{4\pi} \int_{S^1} \< [\Theta_1, \Theta_1], \partial \gamma_2^{\vphantom{-1}} \gamma_2^{-1} \> + \< \Theta_1, d(\partial\gamma_2^{\vphantom{-1}} \gamma_2^{-1})\> 
+ \< \Theta_2,\partial (ad(\gamma_2^{-1}) \Theta_1 ) \> \, d\theta,
$$
and using the identities above and integration by parts, we have
$$
\d R = \frac{i}{2 \pi} \int_{S^1} \< d \Theta_1, \partial \gamma_2^{\vphantom{-1}} \gamma_2^{-1}\> - \< \Theta_1, d(\partial\gamma_2^{\vphantom{-1}} \gamma_2^{-1}) \>\, d\theta.
$$
Therefore, if we define
$$
\a = \frac{i}{2 \pi} \int_{S^1} \< \pi_2^*\Theta, \pi_1^* Z\> \, d\theta,
$$
for $Z \colon LG \to L\fg$ the function $\gamma \mapsto \partial \gamma \gamma^{-1},$ then we see that $d\a = \d R.$ Also, one can check that  $\d\a = 0.$

\subsubsection{A connection for the lifting bundle gerbe}

Now that we have a construction of $\widehat{LG}$ in terms of the differential forms $R$ and $\a,$ we can consider the problem of lifting the $LG$-bundle $P \to M$ to an $\widehat{LG}$-bundle $\widehat{P} \to M.$ We can write down the lifting bundle gerbe for this problem, that is, the bundle gerbe $(\tau^*\widehat{LG}, P)$ over $M,$ and we would like a connection on this bundle gerbe so we can calculate its Dixmier-Douady class. 

Consider, then, the map $\tau \colon P^{[2]} \to LG$ above. We can extend this to a map $\tau \colon P^{[k+1]} \to LG^k$ by defining
$$
\tau (p_1, \ldots, p_{k+1}) = (\tau(p_1, p_2), \ldots, \tau(p_{k}, p_{k+1})).
$$
This is a \emph{simplicial map}. That is, it commutes with the face and degeneracy maps for the simplicial manifolds $\{P^{[k]}\}$ and $\{LG^k\}.$ This means that for differential forms on these manifolds, $\d$ commutes with pull-back by $\tau.$ Now consider the connection $\nu$ on $\widehat{LG}$ (whose curvature is the form $R$). The natural choice for a bundle gerbe connection would be the pull-back, $\tau^*\nu,$ of this form to $\tau^*\widehat{LG}.$ However, $\tau^*\nu$ is not a bundle gerbe connection because it does not respect the product. That is, $s^*(\d\tau^*\nu)$ is non-zero. We know from the discussion on bundle gerbe connections in section \ref{S:bundle gerbes} that $\d(s^*(\d\tau^*\nu)) = 0$ and so there is some form $\epsilon$ on $P^{[2]}$ such that $\d \epsilon = s^*(\d\tau^*\nu).$ Then $\tau^*\nu - \epsilon$ will be a bundle gerbe connection on $\tau^*\widehat{LG}.$ In fact, in this case, since $\a = s^*(\d\nu)$ by definition, we have $s^*(\d \tau^*\nu) = \tau^*\a.$ So $\d(s^*(\d \tau^*\nu)) = \d \tau^*\a = \tau^* \d\a = 0$ as $\d\a = 0$ and so $\epsilon$ satisfies $\d\epsilon = \tau^* \a.$ Thus it suffices to find a 1-form $\epsilon$ on $P^{[2]}$ satisfying $\d\epsilon = \tau^* \a.$

The form $\tau^*\a$ is given by
$$
\frac{i}{2\pi} \int_{S^1} \< \tau_{12}^* \Theta , \tau_{23}^*Z\>\, d\theta
$$
where we have written $\tau_{ij}$ for $\tau(p_i, p_j).$ In order to solve for $\epsilon,$ we need to choose a connection $A$ on $P.$ Then using the equation $p_1 \tau(p_1, p_2) = p_2$ and the Leibnitz rule (see \cite{Kobayashi:1963}), we find the identity
$$
\pi_1^*A = ad(\tau_{12}^{-1})\pi_2^* A + \tau_{12}^* \Theta.
$$
Therefore we have
$$
\tau^*\a = \frac{i}{2\pi} \int_{S^1} \< \pi_{13}^* A - ad(\tau_{12}^{-1}) \pi_{23}^* A , \partial \tau_{23}^{\vphantom{-1}} \tau_{23}^{-1} \> \, d\theta,
$$
where $\pi_{23}(p_1, p_2, p_3) = p_1,$ etc. Now define
$$
\epsilon = \frac{i}{2\pi} \int_{S^1} \< \pi_2^* A , \tau^*Z\>\, d\theta.
$$
Then, using the simplicial identities and the fact that $\tau_{ij}\tau_{jk} = \tau_{ik},$ we have
\begin{align*}
\d \epsilon		&= \pi_1^* \epsilon - \pi_2^*\epsilon + \pi_3^*\epsilon\\
			&= \frac{i}{2\pi}\int_{S^1} \<\pi_{13}^* A, \tau_{23}^*Z\> - \<\pi_{23}^* A, \tau_{13}^*Z \> + \< \pi_{23}^* A, \tau_{12}^*Z\> \, d\theta\\
			&= \frac{i}{2\pi}\int_{S^1} \<\pi_{13}^* A, \tau_{23}^*Z\> - \<\pi_{23}^* A, ad(\tau_{12}^{\vphantom{^*}}) \tau_{23}^* Z  \> \, d\theta\\
			& = \frac{i}{2\pi} \int_{S^1} \< \pi_{13}^* A - ad(\tau_{12}^{-1}) \pi_{23}^* A , \partial \tau_{23}^{\vphantom{-1}} \tau_{23}^{-1} \> \, d\theta.
\end{align*}
It turns out \cite{Stevenson:comm} that in general, $\epsilon$ can be written in terms of $\a$ and $A.$ We shall demonstrate in section \ref{S:LGxS^1 string class} how to find $\epsilon$ in general.

Since we want to calculate the 3-curvature of the lifting bundle gerbe, we are really interested in the curvature of the connection $\tau^*\nu - \epsilon.$ This is given by $\tau^*R - d\epsilon.$ Using the identities given above, we have
\begin{align*}
\tau^*R	&= \frac{i}{4\pi} \int_{S^1} \< \tau^*\Theta, \partial \tau^*\Theta\>\, d\theta\\
		&=  \frac{i}{4\pi} \int_{S^1} \< A_2 - ad(\tau^{-1}) A_1, \partial ( A_2 - ad(\tau^{-1}) A_1)\>\, d\theta\\
		&=  \frac{i}{4\pi} \int_{S^1} \< A_2, \partial A_2\> + \< A_1, \partial A_1\> + \< [A_1, A_1], \tau^*Z\> - 2\< ad(\tau^{-1}) A_1, \partial A_2\> \, d\theta,\\
\intertext{and}
d\epsilon	&= \frac{i}{2\pi} \int_{S^1} \< dA_1, \tau^*Z\> - \< A_1, d(\tau^*Z)\> \, d\theta\\
		&= \frac{i}{2\pi} \int_{S^1} \< dA_1, \tau^*Z\> - \<A_1, \partial A_1\>  + \< [A_1, A_1], \tau^*Z\> - \< ad(\tau^{-1})A_1, \partial A_2\> \, d\theta.
\end{align*}
Therefore
$$
\tau^*R - d\epsilon = \frac{i}{4\pi} \int_{S^1} \< \pi_1^*A, \partial \pi_1^* A\> - \< \pi_2^*A, \partial \pi_2^* A\> - 2\< \pi_2^*F, \tau^*Z\> \, d\theta,
$$
where $F = dA + \frac{1}{2} [A, A]$ is the curvature of $A.$

\subsubsection{A curving for the lifting bundle gerbe}

The next step is to find a curving for $\tau^*\widehat{LG}.$ That is, we wish to find some 2-form $B$ on $P$ such that $\d B = \tau^*R - d\epsilon.$ Note that $\d \colon \Omega^2(P) \to \Omega^2(P^{[2]})$ is given by $\d = \pi_1^* - \pi_2^*,$ so we can write $\tau^*R - d\epsilon$ as
$$
\d \left( \frac{i}{4\pi} \int_{S^1} \< A, \partial A\>\, d\theta \right) - \frac{i}{2\pi} \int_{S^1} \< \pi_2^*F, \tau^*Z\>\, d\theta.
$$
Thus we just need to find some $B_2 \in \Omega^2(P)$ such that
$$
\d B_2 = \frac{i}{2\pi} \int_{S^1} \< \pi_2^*F, \tau^*Z\>\, d\theta.
$$
To solve this equation, we use a \emph{Higgs field} for the bundle $P.$ A Higgs field is a map $\Phi \colon P \to L\fg$ satisfying
$$
\Phi (p \gamma) = ad(\gamma^{-1}) \Phi(p) + \gamma^{-1} \partial \gamma.
$$
It is clear that Higgs fields exist. Since they exist when $P$ is trivial and convex combinations of Higgs fields are also Higgs fields, we can use a partition of unity to construct a Higgs field in general. We shall explain the geometric significance of this map in the next section. For now, note that if we pull back $\Phi$ to $P^{[2]}$ it satisfies
$$
ad(\tau) \pi_1^* \Phi = \pi_2^* \Phi + \tau^*Z.
$$
This just comes from the condition above and the definition of $\tau.$ Then we see that
\begin{align*}
\<\pi_2^*F, \tau^*Z\>	&= \<\pi_2^*F, ad(\tau)\pi_1^*\Phi\> - \< \pi_2^*F, \pi_2^*\Phi\>\\
				&= \<ad(\tau^{-1})\pi_2^*F, \pi_1^*\Phi\> - \< \pi_2^*F, \pi_2^*\Phi\>.
\end{align*}
But one can demonstrate (in a similar manner to the proof of the equation above relating $\pi_1^*A$ and $\pi_2^*A$) that the curvature $F$ satisfies
$$
\pi_1^*F = ad(\tau^{-1})\pi_2^*F
$$
and so we have
$$
\<\pi_2^*F, \tau^*Z\> = \<\pi_1^*F, \pi_1^*\Phi\> - \< \pi_2^*F, \pi_2^*\Phi\>.
$$
Therefore, a curving is given by
$$
B = \frac{i}{2\pi} \int_{S^1} \tfrac{1}{2} \< A, \partial A\> - \< F, \Phi\> \,d \theta.
$$

\subsubsection{The string class of an $LG$-bundle}

Now that we have a curving for the lifting bundle gerbe we can find a representative for the string class $s(P)$ by calculating the 3-curvature $H = dB.$ We have
$$
dB =  \frac{i}{2\pi} \int_{S^1} \tfrac{1}{2} \< dA, \partial A\> - \tfrac{1}{2} \< A, \partial dA\> - \< dF, \Phi\> - \< F, d\Phi\> \,d \theta.
$$
Integration by parts and the Bianchi identity $dF = [F, A]$ yields
$$
dB =  \frac{i}{2\pi} \int_{S^1}  \< dA, \partial A\>  - \< F, [A,\Phi]\> - \< F, d\Phi\> \,d \theta
$$
and since the integral over the circle of $\< [A, A], \partial A\>$ vanishes, we find
$$
dB =  \frac{i}{2\pi} \int_{S^1}  \< F, \partial A\>  - \< F, [A,\Phi]\> - \< F, d\Phi\> \,d \theta.
$$
This descends to a form on $M$ and so
$$
H = -\frac{i}{2\pi} \int_{S^1} \< F, \nabla \Phi\> \, d\theta,
$$
where
$$
\nabla \Phi = d\Phi + [A, \Phi] - \partial A.
$$
Thus we have the result from \cite{Murray:2003}
\begin{theorem}[\cite{Murray:2003}]\label{T:LGstring}
Let $P \to M$ be a principal $LG$-bundle. Let $A$ be a connection on $P$ with curvature $F$ and let $\Phi$ be a Higgs field for $P.$ Then the string class of $P$ is represented in de Rham cohomology by the form
$$
-\frac{1}{4\pi^2} \int_{S^1} \< F, \nabla \Phi\> \, d\theta,
$$
where $\nabla \Phi$ is the covariant derivative above.
\end{theorem}

\section{Higgs fields, $LG$-bundles and the string class}\label{S:Higgs fields, LG-bundles,...}

Recall Killingback's result from section \ref{S:string structures} regarding string structures of a loop bundle. That is, if $Q \to M$ is a principal $G$-bundle and $LQ \to LM$ is the $LG$-bundle obtained by taking loops, then the string class of $LQ$ is the transgression of the first Pontrjagyn class of $Q,$ i.e.
$$
s(LQ) = \int_{S^1} \ev^* p_1 (Q).
$$
In the last section we obtained, following the methods of \cite{Murray:2003}, a general expression for the string class of a principal $LG$-bundle $P \to M$ which is not necessarily a loop bundle. In this case we can prove a result analogous to Killingback's by using a correspondence between $LG$-bundles and certain $G$-bundles. This will also enable us to provide an easy proof of Killingback's result.

\subsection{Higgs fields and $LG$-bundles}\label{SS:LG correspondences}

The following correspondence first appeared in \cite{Garland:1988} in a study of calorons (monopoles for the loop group) and, in the context in which we are interested, in \cite{Murray:2003}. We shall present the construction here in some detail since we will generalise this result in section \ref{S:Higgs fields, LGxS^1-bundles,...} to $LG\rtimes S^1$-bundles and it will be instructive to see the introductory case in depth.

We wish to set up a bijective correspondence between $LG$-bundles over $M$ and $G$-bundles over $M\times S^1.$ Consider the $LG$-bundle $P\times S^1 \to M\times S^1$ where the $LG$ action is trivial on the $S^1$ factor. Then use the evaluation map $\ev \colon LG \times S^1 \to G$ to form the associated $G$-bundle $\widetilde{P} \to M\times S^1$. That is, define $\widetilde{P}$ by
$$
\widetilde{P} = (P\times G\times S^1)/LG
$$
where $LG$ acts on $P\times G \times S^1$ by $(p,g,\theta)\gamma = (p\gamma, \gamma(\theta)^{-1} g, \theta).$ Then there is a right $G$ action on $\widetilde{P}$ given by $[p, g, \theta]h = [p, gh, \theta]$ (where square brackets denote equivalence classes) and a projection $\tilde{\pi} \colon \widetilde{P} \to M\times S^1$ given by $\tilde{\pi}([p, g, \theta]) = (\pi(p), \theta).$ This action is free and transitive on the fibres (which are the orbits of the $G$ action) and hence $\widetilde{P} \to M \times S^1$ is a principal $G$-bundle.

Conversely, given a $G$-bundle $\widetilde{P} \to M\times S^1$ we can define fibrewise an $LG$-bundle $P\to M$ by taking sections of $\widetilde{P}$ restricted to a point in $M.$ That is, the fibre of $P$ over $m$ is
$$
P_m = \Gamma (\widetilde{P}_{|\{m\}\times S^1})
$$
or
$$
P_m= \{f \colon S^1 \to \widetilde{P} \,|\, \tilde{\pi}(f(\theta)) = (m,\theta)\}.
$$
The $LG$ action here is the obvious one derived from the $G$ action on $\widetilde{P}.$ The transition functions of this bundle are simply the transition functions of $\widetilde{P}$ considered as functions from an open set of $M$ to $LG,$ for if $\{U_\a \times S^1\}$ is an open cover of $M\times S^1$ and $\tilde{s}_\a$ is a section of $\widetilde{P}$ then since elements of $P$ are loops in $\widetilde{P},$ a section of $P$ is given by $s_\a (m)(\theta) = \tilde{s}_\a (m,\theta).$ If $s_\b$ is another such section, then the transition functions of $P,$ $g_{\a\b} \colon U_\a\cap U_\b \to LG,$ are given by
$$
s_\b(m) = s_\a(m) g_{\a\b}(m).
$$
Evaluating at $\theta$ gives
$$
s_\b(m)(\theta) = s_\a(m)(\theta) g_{\a\b}(m)(\theta).
$$
But $s_\b(m)(\theta) = \tilde{s}_\b(m,\theta)$ (and similarly for $\a$), so we have
$$
g_{\a\b}(m)(\theta) = \tilde{g}_{\a\b}(m,\theta)
$$
where $\tilde{g}_{\a\b}$ are the transition functions for $\widetilde{P}.$ We can actually give a global description of this bundle quite easily by considering the map
$$
\eta \colon M\to L(M\times S^1); \quad m\mapsto (\theta \mapsto (m, \theta)).
$$
That is, $\eta(m)(\theta) = (m, \theta).$ Then the bundle $P$ is the pullback of the $LG$-bundle $L\widetilde{P} \to L(M\times S^1):$
$$
\xymatrix{ **[l] \eta^*L\widetilde{P} =P\ar[r] \ar[d] &L \widetilde{P}\ar[d]	\\
			M\ar^-\eta[r]	& L(M\times S^1).}
$$
Thus we have a way of constructing a $G$-bundle given an $LG$-bundle and vice versa. It remains to be shown that this is a bijection on the set of isomorphism classes of these bundles. That is, if we start with a $G$-bundle $\widetilde{P}$ and construct $P$ and then form the $G$-bundle corresponding to that bundle, say $\widetilde{P}',$ we have that $\widetilde{P}'$ is isomorphic to $\widetilde{P}.$ And similarly, if we start with $P$ and construct $\widetilde{P}$ and then construct the $LG$-bundle corresponding to that, say $P',$ then these are isomorphic. To see this, first consider a $G$-bundle $\widetilde{P}$ and construct $P$ as above. Then $\widetilde{P}'$ is given by
$$
\widetilde{P}' = (P \times G \times S^1)/LG
$$
where for $[p, g, \theta] \in (P \times G \times S^1)/LG,\, p$ is a map $S^1 \to \widetilde{P}$ as above. Define a bundle map by
$$
f \colon \widetilde{P}' \to \widetilde{P}; \quad [p,g,\theta] \mapsto p(\theta)g.
$$
This is well-defined, since $[p\gamma, \gamma(\theta)^{-1}g, \theta] \stackrel{f}{\mapsto} (p\gamma)(\theta) \gamma(\theta)^{-1}g = p(\theta)g$ and commutes with the $G$ action, since $[p,g,\theta]h = [p, gh, \theta] \stackrel{f}{\mapsto} p(\theta)gh = (p(\theta)g)h.$ Hence $f$ is a bundle isomorphism. On the other hand, if we consider an $LG$-bundle $P$ and construct $\widetilde{P} = (P\times G\times S^1)/LG$ then $P'$  is given by the pull-back above. Notice that if we define the map $\hat{\eta} \colon P \to L\widetilde{P}$ by
$$
\hat{\eta}(p)(\theta) = [p, 1, \theta]
$$
then $\hat{\eta}$ covers $\eta \colon M\to L(M\times S^1),$ that is,
$$
\xymatrix{ P\ar^{\hat{\eta}}[r] \ar[d] &L \widetilde{P}\ar[d]	\\
			M\ar^-\eta[r]	& L(M\times S^1)}
$$
commutes, and so $P$ is isomorphic to the pull-back $P'.$ Thus we have proven

\begin{proposition}[\cite{Garland:1988,Murray:2003}]\label{P:LG correspondences}
There is a bijective correspondence between isomorphism classes of principal $G$-bundles over $M\times S^1$ and isomorphism classes of principal $LG$-bundles over $M.$
\end{proposition}

Importantly for our purposes, this correspondence holds on the level of connections as well. More specifically, if we have a $G$-bundle with connection we can construct an $LG$-bundle with connection and Higgs field and, conversely, given an $LG$-bundle with connection and Higgs field we can construct a $G$-bundle with connection. We shall see that the Higgs field is essentially the $S^1$ component of the connection on $\widetilde{P}.$

Suppose we have a connection $\tilde{A}$ on $\widetilde{P}.$ We can define a connection on $P$ (which is an $L\fg$-valued 1-form) by $A_p(X)(\theta) = \tilde{A}_{p(\theta)}(X_\theta),$ where $X$ is a vector in $T_p P$ (i.e. a vector field along $p$ in $\widetilde{P}$), and so $X_\theta \in T_{p(\theta)} \widetilde{P}.$ This is a connection by virtue of the fact that $\tilde{A}$ is. If we view $\tilde{A}$ as a splitting of the tangent space at each point in $\widetilde{P},$ then we can easily see that $A$ is given by essentially the same splitting since for each $\theta \in S^1,$ $T_pP$ splits as
$$
T_{p(\theta)}\widetilde{P} \simeq V_{p(\theta)} \widetilde{P} \oplus H_{p(\theta)}\widetilde{P}
$$
where $V_{p(\theta)} \widetilde{P}$ is the vertical subspace at $p(\theta)$ and $H_{p(\theta)}\widetilde{P}$ is the horizontal subspace.

Suppose instead we are given an $LG$-bundle $P$ with connection $A$ and Higgs field $\Phi.$ Then we can define a form on $P\times G \times S^1$ by
$$
\tilde{A} = ad(g^{-1}) A(\theta) + \Theta + ad(g^{-1}) \Phi \, d\theta.
$$
This form descends to a form on $\widetilde{P}$ and the connection (also called $\tilde{A}$) is given by this equation considered as a form on $(P\times G\times S^1)/LG.$ To show that this is well defined, we need to check that it is independent of the lift of a vector  in $\widetilde{P}.$ That is, if $\hat{X}$ and $\hat{X}'$ are two lifts of the vector $X\in T_{[p,g , \theta]}\widetilde{P}$ to the fibre in $P\times G\times S^1$ above $[p,g,\theta],$ then $\tilde{A}(\hat{X}) = \tilde{A}(\hat{X}').$ Suppose then, that $\hat{X} \in T_{(p,g,\theta)}(P\times G\times S^1)$ and $\hat{X}' \in T_{(p,g,\theta)\gamma}(P\times G\times S^1).$ Then $\hat{X}\gamma \in T_{(p,g,\theta)\gamma}(P\times G\times S^1),$ and $\hat{X}'$ and $\hat{X}\gamma$ differ by a vertical vector (with respect to the $LG$ action) at $(p,g,\theta)\gamma = (p\gamma,\gamma(\theta)^{-1}g, \theta)$ and so it is sufficient to show that $\tilde{A}$ is zero on vertical vectors and invariant under the $LG$ action (since then $\tilde{A}(\hat{X}') = \tilde{A}(\hat{X}\gamma + \text{vertical}) = \tilde{A}(\hat{X})$). Because any compact Lie group has a faithful representation as matrix group \cite{Pressley-Segal}, we can expand the exponential map as $\exp(t\xi) = 1 + t\xi +\ldots.$ Therefore, the vertical vector at $(p,g,\theta)$ generated by $\xi \in L\fg$ is
\begin{align*}
V 	&= \frac{d}{dt}{\bigg|_0}(p,g,\theta)\exp(t\xi)\\
	&= \frac{d}{dt}{\bigg|_0}(p\exp(t\xi),\exp(-t\xi(\theta))g,\theta)\\
	&= (\iota_p(\xi), -\xi(\theta)g, 0),
\end{align*}
(where we have written $\left.\frac{d}{dt}\right|_0$ for the derivative evaluated at $t=0$), and so
\begin{align*}
\tilde{A}(V) 	&= ad(g^{-1}) A(\iota_p(\xi))(\theta) - g^{-1}\xi(\theta)g\\
			&= g^{-1} \xi(\theta)g - g^{-1}\xi(\theta)g\\
			&= 0.
\end{align*}
So $\tilde{A}$ is zero on vertical vectors. Now, suppose $\hat{X} = (X, g\zeta, x_\theta)$ is given by
$$
\frac{d}{dt}{\bigg|_0} ( \gamma_X (t), g\exp(t\zeta), \theta + tx),
$$
where $\gamma_X (t)$ is a path in $P$ whose tangent vector at $0$ is $X$ and where $\zeta$ and $x$ are elements of the Lie algebras of $G$ and $S^1$ respectively. Then
\begin{align*}
\hat{X}\gamma	& = \frac{d}{dt}{\bigg|_0} ( \gamma_X (t)\gamma, \gamma(\theta + tx)g\exp(t\zeta), \theta + tx)\\
			&= \frac{d}{dt}{\bigg|_0} ( \gamma_X (t)\gamma, \gamma(\theta)gt\zeta + tx\partial\gamma(\theta)g, \theta + tx)\\
			&= ( X\gamma, \gamma(\theta)g(\zeta + x ad(g^{-1})\gamma(\theta)^{-1}\partial\gamma(\theta)), x).
\end{align*}
So
\begin{multline*}
\tilde{A}_{(p\gamma, \gamma(\theta)^{-1}g, \theta)} (\hat{X}\gamma) = \tilde{A}_{(p\gamma, \gamma(\theta)^{-1}g, \theta)} (X\gamma, \gamma(\theta)g(\zeta + x ad(g^{-1})\gamma(\theta)^{-1}\partial\gamma(\theta)), x)\\
\shoveleft{\phantom{\tilde{A}_{(p\gamma, \gamma(\theta)^{-1}g, \theta)} (\hat{X}\gamma)} = ad((\gamma(\theta)^{-1}g)^{-1}) A( X\gamma) + \zeta + x ad(g^{-1})\gamma(\theta)^{-1}\partial\gamma(\theta)}\\
	\shoveright{ + ad((\gamma(\theta)^{-1}g)^{-1})x\Phi(p\gamma)}\\
\shoveleft{\phantom{\tilde{A}_{(p\gamma, \gamma(\theta)^{-1}g, \theta)} (\hat{X}\gamma)} = ad(g^{-1})ad(\gamma)ad(\gamma^{-1}) A(X)(\theta) + \zeta + x ad(g^{-1})\gamma(\theta)^{-1}\partial\gamma(\theta)}\\
	\shoveright{ + ad(g^{-1})x ad(\gamma)(ad(\gamma^{-1})\Phi(p) + \gamma^{-1} \partial \gamma)}\\
\shoveleft{\phantom{\tilde{A}_{(p\gamma, \gamma(\theta)^{-1}g, \theta)} (\hat{X}\gamma)} = ad(g^{-1}) A(X)(\theta) + \zeta + ad(g^{-1})x \Phi(p).}\\
\end{multline*}
Therefore $\tilde{A}$ is invariant under the $LG$ action and so defines a form on $\widetilde{P}.$ This form is a connection form since if $[X, g\zeta, x_\theta]$ is a vector at $[p,g,\theta],$ then $[X, g\zeta, x_\theta]h = [X, gh\, ad(h^{-1})\zeta , x_\theta]$ and so
\begin{align*}
\tilde{A} ([X, g\zeta, x_\theta]h)	&= ad(h^{-1}g^{-1}) A(X)(\theta) + ad(h^{-1})\zeta  + ad(h^{-1}g^{-1})x\Phi(p)\\
						&= ad(h^{-1}) \tilde{A}([X, g\zeta ,x_\theta])
\end{align*}
and further, the vertical vector at $[p,g,\theta]$ generated by $\zeta \in \fg$ is given by
\begin{align*}
V_\zeta 	&= \frac{d}{dt}_{|_0}[p,g \exp(t\zeta),\theta]\\
		&= [0, g\zeta, 0]
\end{align*}
and so $\tilde{A}(V_\zeta) = \zeta.$

We have shown already that the correspondence outlined above is a bijection between isomorphism classes of bundles. Now we will show that in fact it is a bijection between isomorphism classes of bundles with connection. So given a $G$-bundle $\widetilde{P}$ with connection $\tilde{A},$ we construct the $LG$-bundle $P$ with the connection $A$ as above. Then construct the $G$-bundle $\widetilde{P}'$ (which is isomorphic to $\widetilde{P}$) and give it the connection $\tilde{A}'$ which we just outlined. Of course, to do this we'll need a Higgs field for $P.$ Recalling that elements of $P$ are essentially loops in $\widetilde{P},$ we can define a Higgs field by
$$
\Phi (p) = \tilde{A}(\partial p).
$$
This is a Higgs field since if we calculate $\Phi(p\gamma)$ we get
\begin{align*}
\tilde{A}(\partial (p\gamma)) 	&= \tilde{A}( (p\gamma)_* \frac{\partial}{\partial \theta})\\
					&= \tilde{A}((\partial p)\gamma + \iota_{p\gamma} (\gamma^{-1} \partial\gamma))\\
					&= ad(\gamma^{-1})\tilde{A}(\partial p) + \gamma^{-1}\partial \gamma.
\end{align*}
(Note that this is essentially the $S^1$ part of $\tilde{A}.$ That is, if we take a section $\tilde{s}$ of $\widetilde{P} \to M \times S^1$ we can get a section $s$ of $P\to M$ by $s(m)(\theta):= \tilde{s}(m,\theta).$ Then if we pull-back $\Phi$ by $s$ we get
\begin{align*}
(s^*\Phi)(m)(\theta) &= (\tilde{s}^*\tilde{A})(m,\theta)\left(\frac{\partial}{\partial\theta}\right)\\
				&=(\tilde{s}^*\tilde{A})_\theta (m, \theta)
\end{align*}
where $(\tilde{s}^*\tilde{A})_\theta$ is the $S^1$ part of $(\tilde{s}^*\tilde{A})$ -- i.e. the coefficient of $d\theta$ -- and since the $\frac{\partial}{\partial\theta}$ kills all but the $d\theta$ part.)

Therefore, the connection $\tilde{A}'$ is given in terms of $\tilde{A}$ as
$$
\tilde{A}'_{[p,g,\theta]} = ad(g^{-1}) \tilde{A}_{p(\theta)} + \Theta + ad(g^{-1}) \tilde{A}(\partial p) d\theta.
$$
Recall that $\widetilde{P}'$ is isomorphic to $\widetilde{P}$ via the map
$$
f \colon \widetilde{P}' \to \widetilde{P}; \quad [p,g,\theta] \mapsto p(\theta)g,
$$
so we would like to have $ f^* \tilde{A} = \tilde{A}'.$ Now,  $ f^* \tilde{A}([X, g\zeta, x_\theta]) = \tilde{A}(f_* [X, g\zeta, x_\theta])$ and, as before, if $\gamma_X (t)$ is a path in $P$ whose tangent vector at $0$ is $X$ and if $\zeta$ and $x$ are elements of the Lie algebras of $G$ and $S^1$ respectively, then
\begin{multline*}
f_* [X, g\zeta, x_\theta] = \frac{d}{dt}{\bigg|_0} (\gamma_X(t) (\theta + tx) g\exp(t\zeta))\\
\shoveleft{\phantom{f_* [X, g\zeta, x_\theta]} =\left(\frac{d}{dt}(\gamma_X (t)) (\theta + tx) g\exp(t\zeta)+ \gamma_X (t)(\theta + tx) g \frac{d}{dt} \exp(t\zeta)\right.}\\
\shoveright{ + \partial \gamma_X (t)(\theta+ tx) x g \exp(t\zeta)\left.\vphantom{\frac{d}{dt}}\right){\bigg|_0}}\\
\shoveleft{\phantom{f_* [X, g\zeta, x_\theta]} = X(\theta)g + \iota_{p(\theta)g}(\zeta) + \partial p(\theta) x g}\\
\end{multline*}
and so
\begin{align*}
f^* \tilde{A}([X, g\zeta, x_\theta]) &= ad(g^{-1})A(X) + \zeta + ad(g^{-1}) A(\partial p(\theta))x\\
						&= \tilde{A}'([X, g\zeta, x_\theta]).
\end{align*}
If, on the other hand, we had started with the $LG$-bundle $P$ with connection $A$ (and Higgs field $\Phi$), then $A'$ would be given by
$$
A'_{p}(X)(\theta) = \tilde{A}_{p(\theta)} (X_\theta)
$$
and recalling that the isomorphism between $P$ and $P'$ is essentially given by $f(p) = (\theta \mapsto [p, 1, \theta]),$ we have
$$
f^*A'_{p}(X)(\theta) = A_p(X)(\theta).
$$
Hence, we have

\begin{proposition}[\cite{Murray:2003}]\label{P:LG connection correspondences}
The correspondence from Proposition \ref{P:LG correspondences} extends to a bijection between $G$-bundles on $M \times S^1$ with connection and $LG$-bundles on $M$ with connection and Higgs field.
\end{proposition}


\subsection{The string class and the first Pontrjagyn class}\label{SS:Killingback}

As mentioned previously, the correspondence above provides us with a result analogous to Killingback's. We have

\begin{theorem}[\cite{Murray:2003}]\label{T:LGPont}
Let $P \to M$ be an $LG$-bundle and $\widetilde{P} \to M\times S^1$ the corresponding $G$-bundle. Then the string class of $P$ is given by integrating over the circle the first Pontrjagyn class of $\widetilde{P}.$ That is,
$$
s(P) = \int_{S^1} p_1(\widetilde{P}).
$$
\end{theorem}

\begin{proof}
If $\tilde{F}$ is the curvature of a connection on $\widetilde{P}$ then the Pontrjagyn form is given by
$$
p_1(\widetilde{P}) = -\frac{1}{8\pi^2} \< \tilde{F}, \tilde{F}\>.
$$
In this case we know that $\tilde{A}$ is given as in the previous section. That is,
$$
\tilde{A} = ad(g^{-1}) A + \Theta + ad(g^{-1}) \Phi\, d\theta,
$$
so we can calculate its curvature using $\tilde{F} = d\tilde{A} + \frac{1}{2} [\tilde{A}, \tilde{A}].$ Now,
\begin{multline*}
\tfrac{1}{2}[\tilde{A}, \tilde{A}] = \tfrac{1}{2}[ad(g^{-1}) A + \Theta + ad(g^{-1}) \Phi\, d\theta, ad(g^{-1}) A + \Theta + ad(g^{-1}) \Phi\, d\theta]\\
\shoveleft{\phantom{\tfrac{1}{2}[\tilde{A}, \tilde{A}]} = \tfrac{1}{2}ad(g^{-1})[A , A] + \tfrac{1}{2}[\Theta, \Theta]  + [\Theta, ad(g^{-1})A]}\\
 + ad(g^{-1})[A, \Phi]d\theta + [\Theta, ad(g^{-1}) \Phi]d\theta.
\end{multline*}
So we just need to calculate $d \tilde{A} = d ( ad(g^{-1})A) + d\Theta + d(ad(g^{-1}) \Phi) d\theta.$ Now, if $\omega$ is a 1-form then for tangent vectors $X$ and $Y$ we have
$$
d\omega (X, Y) = \tfrac{1}{2} \left\{ X(\omega (Y)) - Y(\omega (X)) - \omega ([X, Y]) \right\},
$$
so let $(X, g\xi, x_\theta)$ and $(Y, g\zeta, y_\theta)$ be two tangent vectors to $\widetilde{P}$ at the point $[p, g, \theta].$ Then for $d (ad(g^{-1}) A),$ first calculate
\begin{align*}
(X, g\xi, x_\theta)(ad(g^{-1})  & A_p(Y)_\theta)\\
		&= \frac{d}{dt}\bigg|_0 (1-t\xi)g^{-1} A_{\gamma_X(t)} (Y)_{(\theta + tx)} g (1+ t\xi)\\
		&= \frac{d}{dt}\bigg|_0 \left( ad(g^{-1}) A_{\gamma_X(t)} (Y)_\theta \right) + ad(g^{-1})\partial A_p(Y) x - [\xi, ad(g^{-1})A_p (Y)_\theta].
\end{align*}
This yields
$$
d(ad(g^{-1})A) = ad(g^{-1}) dA - ad(g^{-1}) \partial A \wedge d\theta - [\Theta, ad(g^{-1})A].
$$
Similarly, for $d(ad(g^{-1})\Phi)d\theta$ we have
\begin{align*}
(X, g\xi, x_\theta)(ad(g^{-1})  & \Phi(p)_\theta )\\
		&= \frac{d}{dt}\bigg|_0 (1-t\xi)g^{-1} \Phi(\gamma_X(t))_{(\theta + tx)} g (1+ t\xi)\\
		&= \frac{d}{dt}\bigg|_0 \left( ad(g^{-1}) \Phi(\gamma_X(t))\right) + ad(g^{-1})\partial \Phi x - [\xi, ad(g^{-1}) \Phi(p)_\theta],
\end{align*}
and so
$$
d(ad (g^{-1})\Phi) d\theta = ad(g^{-1}) d\Phi \wedge d\theta- [\Theta, ad(g^{-1}) \Phi ] d\theta.
$$
Putting these together gives
\begin{multline*}
\tilde{F} = ad(g^{-1}) dA - ad(g^{-1}) \partial A \wedge d\theta - [\Theta, ad(g^{-1})A] + d\Theta\\
 + ad(g^{-1}) d\Phi \wedge d\theta - [\Theta, ad(g^{-1}) \Phi ] d\theta +  \tfrac{1}{2}ad(g^{-1})[A , A]\\
  + \tfrac{1}{2}[\Theta, \Theta]  + [\Theta, ad(g^{-1})A]
\shoveright{ + ad(g^{-1})[A, \Phi]d\theta + [\Theta, ad(g^{-1}) \Phi]d\theta}\\
\shoveleft{\phantom{\tilde{F}} = ad(g^{-1})\left( dA + \tfrac{1}{2}[A, A] + d\Phi\wedge d\theta + [A, \Phi]d\theta - \partial A \wedge d\theta\right)}\\
\end{multline*}
That is,
$$
\tilde{F} = ad(g^{-1}) \left( F + \nabla \Phi\, d\theta \right).
$$
Then the Pontrjagyn form is given by
$$
p_1(\widetilde{P}) = -\frac{1}{8\pi^2} \left( \< F, F\> + 2\< F, \nabla \Phi\>\, d\theta \right),
$$
and integrating over the circle gives the required result.
\end{proof}

\subsubsection{A proof of Killingback's result}

We now have a result which is more general than Killingback's result since it can be applied to a general $LG$-bundle, not just a loop bundle. We now show how Theorem \ref{T:LGPont} gives a method for proving Killingback's result.

\begin{corollary}\label{C:Killingback}
Let $LQ \to LM$ be a loop bundle, that is, a principal $LG$-bundle obtained by taking loops in a $G$-bundle $Q \to M.$ Then
$$
s(LQ) = \int_{S^1} \ev^* p_1(Q).
$$
\end{corollary}

\begin{proof}
We know that the string class of $LQ$ is given by the integral over the circle of the first Pontrjagyn class of the corresponding $G$-bundle over $LM \times S^1.$ We show that this bundle is isomorphic to the pull-back of $Q$ by the evaluation map, then the result follows. The $G$-bundle $\widetilde{LQ}$ is given by $(LQ \times G \times S^1)/LG.$ Define the map $\widetilde{LQ} \to Q$ by $[q, g, \theta] \mapsto q(\theta)g.$ As in section \ref{SS:LG correspondences} above, this map is well-defined and commutes with the $G$-action. Furthermore, it covers the evaluation map $LM \times S^1 \to M$ and so $\widetilde{LQ}$ is isomorphic to $\ev^*Q$ and hence $p_1(\widetilde{LQ}) = \ev^*p_1(Q).$
\end{proof}

\chapter{Higgs fields and characteristic classes for $\Omega G$-bundles}\label{C:characteristic classes}


In our discussion of string structures in chapter \ref{C:background} we were concerned mainly with the loop group $LG$ and its central extension $\widehat{LG}.$ In this chapter we shall, for the most part, be considering the subgroup of $LG$ given by those loops which begin at the identity in $G,$ that is, the \emph{based} loop group, which we shall denote $\Omega G.$ We will return to the discussion of free loops in section \ref{S:universal string class}.

\section{String structures and the path fibration}\label{S:CM'91}

In this section we will outline the result from \cite{Carey:1991} concerning string structures for certain $\Omega G$-bundles.\footnote{Actually, in \cite{Carey:1991} Carey and Murray work with the group of smooth maps from the interval $[0, 2\pi]$ into the group $G$ whose endpoints agree. We shall look more closely at this group in section \ref{S:universal string class}. Here we will be extending their results to the subgroup of based smooth maps $S^1 \to G.$} In particular, we shall see that if $Q \to M$ is a principal $G$-bundle, then the string class for the $\Omega G$-bundle $\Omega Q \to \Omega M$ is a characteristic class for such bundles. To be precise, what we mean here is that we have chosen a base point $m_0$ in $M$ and a base point $q_0$ in the fibre above $m_0$ and then $\Omega Q \to \Omega M$ is an $\Omega G$-bundle. By `string class' we mean the obstruction to lifting $\Omega Q$ to an $\widehat{\Omega G}$-bundle, where $\widehat{\Omega G}$ is the central extension of $\Omega G.$ (Actually, since we are working with differential forms, we are really concerned with the image in real cohomology of the string class -- however, we make no distinction between the terms here.) We will also generalise this to the case of a general $\Omega G$-bundle, that is, one which is not necessarily a loop bundle.

\subsection{Classifying maps and characteristic classes}

In the interests of being self-contained we shall begin by giving a short overview of the theory of classifying maps and characteristic classes before moving on to the specific case we are interested in. Recall that $\cG$-bundles over $M$ are classified by (homotopy classes of) maps to the classifying space $B\cG.$ A $\cG$-bundle is then (isomorphic to) the pull-back by this map of the universal bundle $E\cG \to B\cG.$ This bundle is characterised by the fact that it is a principal $\cG$-bundle and that $E\cG$ is a contractible space. If $P\to M$ is a $\cG$-bundle, a map $f \colon M \to B\cG$ such that $P$ is isomorphic to the pull-back $f^*E\cG$ is called a classifying map for $P.$

A characteristic class associates to a $\cG$-bundle $P\to M$ a class $c(P)$ in $H^*(M).$ It must be natural with respect to pull-backs in the sense that if $g \colon N\to M$ is a smooth map then $c$ must associate to the pull-back bundle $g^*P \to N$ the class given by the pull-back of $c(P).$ That is,
$$
c(g^*P) = g^* c(P).
$$
Note that since all $\cG$-bundles are pulled-back from the universal bundle, then if $P \to M$ is a $\cG$-bundle with classifying map $f,$ all its characteristic classes are of the form $f^*c(E\cG)$ for some characteristic class $c.$ That is, the set of characteristic classes for $\cG$-bundles is in bijective correspondence with the cohomology group $H^*(B\cG).$

\subsection{String structures and the path fibration}

In general, both the classifying space and the universal bundle for a group can be difficult to describe. For the based loop group $\Omega G,$ however, we have the following construction \cite{Carey:2000}: Let $PG$ be the space of paths in $G,$ $p \colon \RR \to G$ such that $p(0)$ is the identity and $p^{-1} \partial p$ is periodic. Then this is acted on by $\Omega G$ and
$$
\xymatrix@C=4ex{\Omega G\ar[r]	& PG\ar[d]\\
					&G}
$$
is an $\Omega G$-bundle called the \emph{path fibration}, where the projection $\pi$ sends a path $p$ to its value at $2\pi.$ $PG$ is contractible and so the path fibration is a model for the universal $\Omega G$-bundle and we have $B\Omega G = G.$ (See Appendix \ref{A:infinite-dimensional} for details.)

Since we are assuming that $G$ is compact, simple and simply connected, we know that $H^3(G, \ZZ) = \ZZ$ and there is an expression for the generator of this group. Namely, the 3-form on $G$ given by
$$
\omega = \frac{1}{48\pi^2} \< \Theta, [\Theta, \Theta] \>.
$$
In \cite{Carey:1991} Carey and Murray show the string class of the path fibration (for the case of loops which are smooth on $(0, 2\pi)$) is given by the 3-form $\omega$ by giving an explicit construction of the lift of $PG$ which exists precisely when this class vanishes. We will use Theorem \ref{T:LGstring} to calculate the string class of the path fibration. Firstly we need a connection on $PG.$ This is given in \cite{Carey:2002}: Let $\a$ be a smooth real-valued function on $[0, 2\pi]$ such that $\a(0) = 0, \a(2\pi) = 1$ and all the derivatives of $\a$ vanish at the endpoints. Then $\a$ can be extended to a function on $\RR$ and a connection in $PG$ is given by
$$
A = \Theta - \a \, ad(p^{-1}) \pi^* \widehat{\Theta},
$$
where $\widehat{\Theta}$ is the \emph{right} invariant Maurer-Cartan form. The horizontal projection of a tangent vector $X$ using this connection is
$$
hX = \a \, X(2\pi) p(2\pi)^{-1} p.
$$
We can calculate the curvature of $A$ using the covariant derivative $F = DA.$ For tangent vectors $X$ and $Y,$ we have
\begin{align*}
F(X, Y)	&= \frac{1}{2} A([hX, hY])\\
		&= \frac{1}{2} A\left(\a^2 \left[X(2\pi) p(2\pi)^{-1}, Y(2\pi) p(2\pi)^{-1}\right] p \right)\\
		&= \frac{1}{2}\left(\Theta -\a \, ad(p^{-1}) \pi^* \widehat{\Theta}\right) \left(\a^2 \left[X(2\pi) p(2\pi)^{-1}, Y(2\pi) p(2\pi)^{-1}\right] p \right)\\
		&= \frac{1}{2}\left(\a^2 -\a \right) ad(p^{-1}) \left[X(2\pi) p(2\pi)^{-1}, Y(2\pi) p(2\pi)^{-1}\right].
\end{align*}
So
$$
F =  \frac{1}{2}\left(\a^2 -\a \right) ad(p^{-1})[ \pi^* \widehat{\Theta}, \pi^* \widehat{\Theta}].
$$
In order to use Theorem \ref{T:LGstring} we also need a Higgs field for $PG.$ Define the map $\Phi \colon PG \to L \fg$ by
$$
\Phi(p) = p^{-1} \partial p.
$$
Then $\Phi$ is a Higgs field, since for $\gamma \in \Omega G$ we have
\begin{align*}
\Phi(p\gamma)	&= (p\gamma)^{-1} \partial (p\gamma)\\
			&= ad(\gamma^{-1}) p^{-1} \partial p + \gamma^{-1} \partial \gamma.
\end{align*}
The formula for the string class uses $\nabla \Phi = d\Phi + [A, \Phi] - \partial A.$ We can calculate
$$
d\Phi = \partial \Theta + [\Phi, \Theta],
$$
$$
[A, \Phi] = [\Theta, \Phi] - \a \, [ad(p^{-1}) \pi^*\widehat{\Theta}, \Phi]
$$
and
$$
\partial A = \partial \Theta -\partial \a \, ad(p^{-1}) \pi^* \widehat{\Theta} - \a \, [ad(p^{-1}) \pi^*\widehat{\Theta}, \Phi].
$$
So we have
$$
\nabla \Phi = \partial \a \, ad(p^{-1}) \pi^* \widehat{\Theta}.
$$
Therefore, by Theorem \ref{T:LGstring} we have
\begin{align*}
s(PG)	&= -\frac{1}{8\pi^2} \int_{S^1} \left\< \left( \a^2 - \a \right) ad(p^{-1})[ \pi^* \widehat{\Theta}, \pi^* \widehat{\Theta}], \partial \a \, ad(p^{-1}) \pi^* \widehat{\Theta} \right\> \, d\theta\\
			&= -\frac{1}{8\pi^2}  \< [\widehat{\Theta}, \widehat{\Theta}], \widehat{\Theta}\> \int_{S^1}\left( \a^2 - \a \right) \partial \a \, d\theta\\
			&= \frac{1}{48\pi^2} \<\Theta, [\Theta, \Theta]\>,
\end{align*}
where the last line follows from the $ad$-invariance of the Killing form. Thus we see that the string class of the path fibration is the generator of the degree three cohomology of $G.$

Now, consider again a based loop bundle $\Omega Q \xrightarrow{\,\Omega G\,} \Omega M.$ In \cite{Carey:1991} Carey and Murray write down the classifying map for such bundles and then show, by explicitly calculating the integral of the (pull-back by the evaluation map of the) first Pontrjagyn class of $Q,$ that the string class is the pull-back by this map of the 3-form $\omega.$ To write down the classifying map of the bundle $\Omega Q \to \Omega M$ choose a connection for it. Then take a loop $\gamma \in \Omega Q$ and project it down to $\pi\circ \gamma \in \Omega M.$ Lift this back up horizontally to $\gamma_h \in \Omega Q,$ so that $\pi\circ \gamma = \pi \circ \gamma_h.$ Then the \emph{holonomy}, $\hol(\gamma) \in PG$ is given by $\gamma = \gamma_h \hol(\gamma).$ This covers the usual holonomy\footnote{Note that we can define the holonomy since we have chosen basepoints in $M$ and $Q.$}  $\hol \colon\Omega M\to G,$ so we have:
$$
\xymatrix{\Omega Q \ar^{\hol}[r] \ar[d] & PG \ar[d]\\
		\Omega M \ar^{\hol}[r]	& G}
$$
Thus $\hol$ is a classifying map for the bundle $\Omega Q \to \Omega M.$ Now, using Corollary \ref{C:Killingback} and by calculating explicitly $\int_{S^1} \ev^*p_1(Q),$  we can show that
$$
s(\Omega Q) = \hol^*\omega + \text{\emph{exact}}.
$$
We shall show this in more detail in the next section when we generalise this result to the case of higher classes for general $\Omega G$-bundles (that is, an $\Omega G$-bundle which is not necessarily a loop bundle). For now let us assume this result and show how it leads us to a more general statement.





To generalise the result above to a general $\Omega G$-bundle $P\xrightarrow{\,\Omega G\,} M,$ we need a classifying map for such bundles. Consider the $\Omega G$-bundle $P\to M.$ Choose a Higgs field $\Phi \colon P \to L \fg$ for $P.$ It is possible to solve the equation $\Phi(p) = g^{-1}\partial g$ for $g \in PG.$ We define the \emph{Higgs field holonomy}, $\hol_\Phi,$ to be the solution to this equation satisfying the initial condition $g(0) =1$. Note that if $\hol_\Phi(p) = g$ then since 
$$
\Phi(p h) = ad(h^{-1})\Phi(p) + h^{-1}\partial h
$$
and 
$$
(gh)^{-1}\partial(gh) = ad(h^{-1}) g^{-1}\partial g + h^{-1}\partial h,
$$
we see that $\hol_\Phi (p\cdot h) = \hol_\Phi(p)h$ and hence $\hol_\Phi$ descends to a map (also called $\hol_\Phi)\, M\to G$ and is a classifying map for $P\to M.$ 

A natural question arises at this point: If $Q \to M$ is a $G$-bundle with connection $A$ then we can define the holonomy of a loop $\gamma \in \Omega Q.$ However, since the loop bundle $\Omega Q \to \Omega M$ is an $\Omega G$-bundle, we can also choose a Higgs field for it and define the Higgs field holonomy of a loop $\gamma$ in this bundle. Can we find the Higgs field $\Phi$ such that $\hol_\Phi = \hol$? Define $\Phi$ in terms of $A$ as in section \ref{S:Higgs fields, LG-bundles,...}, that is,
$$
\Phi(\gamma) = A(\partial \gamma).
$$
Then using $\gamma = \gamma_h \hol(\gamma),$ we find
$$
\partial \gamma =\partial \gamma_h \cdot \hol(\gamma) + \iota_{\gamma_{h}}(\hol(\gamma)^{-1}\partial \hol(\gamma)).
$$
Since $\gamma_h$ is horizontal (in the sense that all its tangent vectors are horizontal), applying the connection form $A$ gives
$$
A(\partial \gamma) = \hol(\gamma)^{-1} \partial \hol(\gamma).
$$
Therefore, $\hol_\Phi = \hol.$

We can extend the result from \cite{Carey:1991} by finding a relationship between $\hol_\Phi$ and $\hol$ in general:

We can modify the correspondence in section \ref{S:Higgs fields, LG-bundles,...}, which relates $LG$-bundles over $M$ and $G$-bundles over $M\times S^1,$ to one which applies to $\Omega G$-bundles. We say a $G$-bundle over $M \times S^1$  is \emph{framed} over $M \times \{0\}$ if it is trivial over $M \times \{0\}$. A particular trivialisation is called a \emph{framing}. Given this, then, $\Omega G$-bundles correspond to $G$-bundles over $M\times S^1$ which are framed over $M\times \{0\}.$ This means we take a $G$-bundle $\widetilde{P} \to M\times S^1$ and a section (i.e.\! a framing) $s \colon M\times \{0\} \to \widetilde{P}$ and the fibre of $P$ over $m$ has a base point given by $s(m, 0).$ Using this correspondence, define a bundle map
$$
\xymatrix{ P \ar^\eta[r] \ar[d]	& \Omega \widetilde{P}\ar[d]\\
		M\ar^(.3)\eta[r]		&\Omega(M\times S^1)}
$$
by $\eta(m) = \theta\mapsto (m,\theta),$ or, on the total space, $\eta(p) = \theta\mapsto [p,1,\theta].$ Then we have:
\begin{lemma}\label{L:holphi=holeta}
Let $P \to M$ be an $\Omega G$-bundle with connection and Higgs field $\Phi,$ $\widetilde{P} \to M \times S^1$ its corresponding $G$-bundle and $\eta$ as above. Then $\ds\hol_\Phi = \hol \circ \eta.$
\end{lemma}

\begin{proof}
If $\tilde{A}$ is the connection form on $\widetilde{P}$ then $\tilde{\Phi}\colon \Omega \widetilde{P} \to L\fg$ defined by
$$
\tilde{\Phi}(\gamma) = \tilde{A}(\partial \gamma)
$$
gives us that
$$
\hol_{\tilde{\Phi}} = \hol
$$
as above. Therefore we need only show that $\hol_{\Phi} = \hol_{\tilde{\Phi}} \circ \eta.$

Let $p\in P.$ Consider the unique horizontal path $\eta(p)_{h}$ such that
$$
\tilde{\pi}(\eta(p)) = \tilde{\pi}(\eta(p)_h)
$$
given by projecting $\eta(p)$ to $\Omega(M\times S^1)$ and lifting horizontally back to $\Omega \widetilde{P}.$ The tangent vector to the loop $\eta(p)$ at the point $\theta$ is given by the derivative $\partial \eta(p)_\theta$ and since $\eta(p)_h$ is horizontal we have that
$$
\tilde{A} (\eta(p)_{h,\theta}) = 0.
$$
Now, $\eta(p)_\theta = [p,1,\theta]$, so we can explicitly calculate $\partial \eta(p)_\theta:$
$$
\frac{\partial}{\partial \theta} \eta(p)_\theta = [0,0,1].
$$
Recall that the connection $\tilde{A}$ is given in terms of the connection $A$ and Higgs field $\Phi$ for $P$ as
$$
\tilde{A} = ad(g^{-1})A + \Theta + ad(g^{-1}) \Phi \, d\theta.
$$
Therefore, we have $\tilde{A}(\partial \eta(p)) = \Phi(p).$ Or, in terms of the Higgs field for $\Omega \widetilde{P},$
$$
\Phi = \tilde{\Phi} \circ \eta.
$$
As above, we have
$$
\tilde{\Phi}(\eta(p)) = \hol(\eta(p))^{-1} \partial \hol(\eta(p)),
$$
and therefore $\hol_{\Phi} = \hol_{\tilde{\Phi}} \circ \eta.$

\end{proof}

We see that $\hol_\Phi$ factors through $\hol.$ In order to use this we need the following result:

\begin{lemma}\label{L:eta* int ev* = int}
In the situation of Lemma \ref{L:holphi=holeta}, for degree 4 differential forms on $M\times S^1$ we have
$$\displaystyle \eta^* \int_{S^1} \ev^* = \int_{S^1}.$$
\end{lemma}

\begin{proof}
Note first that we have
\begin{alignat*}{5}
M&\times S^1 \,&&\xrightarrow{\,\eta \times 1\,} \,\,&&\Omega(M\times S^1) \times S^1 \,&&\xrightarrow{\,\ev\,} \,\,&M\times S^1\\
(m&,\phi) &&\longmapsto && (\theta\mapsto (m,\theta),\phi) &&\longmapsto &(m,\phi)
\end{alignat*}
so, $\ev \circ (\eta \times 1)$ is the identity. Therefore, we have
$$
\int_{S^1} = \int_{S^1} (\eta \times 1)^*  \ev^*,
$$
so it suffices to show that
$$
\int_{S^1} (\eta \times 1)^* = \eta^* \int_{S^1}.
$$
That is, that the following diagram commutes
$$
\xymatrix@C=7ex@R=6ex{
\Omega^4(\Omega (M \times S^1) \times S^1) \ar[r]^-{(\eta \times 1)^*} \ar[d]_{\int_{S^1}} & \Omega^4(M \times S^1) \ar[d]^{\int_{S^1}}\\
\Omega^3(\Omega(M\times S^1)) \ar[r]^-{\eta^*} & \Omega^3(M)
}
$$
Consider $\omega \in \Omega^4(\Omega (M \times S^1) \times S^1).$ Then if $X_1, X_2$ and $X_3$ are tangent vectors to $M$ we have
\begin{align*}
\left( \int_{S^1} (\eta \times 1)^* \omega \right) ( X_1, X_2, X_3) &= \int_{S^1} (\eta \times 1)^* \omega ( \widehat{X}_1, \widehat{X }_2 , \widehat{X}_3)\\
	&= \int_{S^1} \omega ( (\eta\times 1)_* \widehat{X}_1, (\eta\times 1)_* \widehat{X}_2, (\eta\times 1)_* \widehat{X}_3 ),
\end{align*}
where $\widehat{X}_i$ ($i=1,2,3$) is a lift of $X_i$ to $M \times S^1$. On the other hand, if $\widehat{\eta_* X}_i$ is a lift of $\eta_* X_i$ to $\Omega( M \times S^1) \times S^1$, then
\begin{align*}
\eta^* \left( \int_{S^1} \omega \right) (X_1, X_2, X_3) &= \left( \int_{S^1} \omega \right) (\eta_* X_1, \eta_* X_2, \eta_* X_3 )\\
	&= \int_{S^1} \omega ( \widehat{\eta_* X}_1, \widehat{\eta_* X}_2, \widehat{\eta_* X}_3 ).
\end{align*}
Since the expressions above are independent of the lift chosen, we can use the natural splitting of the tangent bundles to $M \times S^1$ and $\Omega (M \times S^1) \times S^1$ to define $\widehat{X}_i = (X_i, 0)$ and $\widehat{\eta_* X}_i = (\eta_* X_i, 0)$ and so we have
\begin{align*}
\left( \int_{S^1} (\eta \times 1)^* \omega \right) ( X_1, X_2, X_3) &= \int_{S^1} \omega ( (\eta_* \widehat{X}_1, 0), (\eta_* \widehat{X}_2, 0) , (\eta_* \widehat{X}_3, 0) )\\
	&= \eta^* \left( \int_{S^1} \omega \right) (X_1, X_2, X_3).
\end{align*}

\end{proof}

Combining Lemmas \ref{L:holphi=holeta} and \ref{L:eta* int ev* = int}, we have:

\begin{theorem}\label{T:s(P) is char}
The string class of an $\Omega G$-bundle $P \to M$ is the characteristic class corresponding to $\omega \in H^3(G)$.
\end{theorem}

\begin{proof}
On the level of cohomology we have
\begin{align*}
s(P) 	&= \int_{S^1} p_1(\widetilde{P})\\
	&= \eta^*\int_{S^1} \ev^*p_1(\widetilde{P})\\
	&= \eta^* s(\Omega \widetilde{P})\\
	&= \eta^*\hol^*\omega\\
	&=\hol_\Phi^*\omega.
\end{align*}

\end{proof}


\section{Higher string classes for $\Omega G$-bundles}\label{S:Omega G characteristic classes}


We have seen in the last section that the string class is a characteristic class for $\Omega G$-bundles and we know from section \ref{S:Higgs fields, LG-bundles,...} (Theorem \ref{T:LGPont}) that it is naturally associated to the first Pontrjagyn class of the corresponding $G$-bundle. Indeed, the fact that the string class is given by integrating the first Pontrjagyn class was used to show that it is natural. In this section we will generalise these ideas to higher degree classes for $\Omega G$-bundles. These classes will be naturally associated to a characteristic class for $G$-bundles in the same way the string class is related to the Pontrjagyn class.

We can summarise the results from the previous section with the following diagram
$$
\xymatrix{ H^4(BG) \ar[r]^-{C\text{-}W_{\widetilde{P}}} \ar_\tau[d] & H^4(M\times S^1) \ar^{\int_{S^1}}[d]\\
		H^3(G) \ar[r]^{\hol_\Phi^*} & H^3(M)}
$$
The top arrow here is the usual Chern-Weil map (see below). The map $\tau$ is the \emph{transgression} (see for example \cite{Chern:1974} or \cite{Heitsch:1974}) which we shall describe presently. As long as $G$ is compact and connected, $H^{2k}(BG)$ is isomorphic to the set of multilinear, symmetric, $ad$-invariant functions on $\fg\times \ldots\times\fg$ ($k$ times). Let $f$ be such a function and let $Q \to M$ be a $G$-bundle with connection. Then the Chern-Weil map, $C\text{-}W_{Q},$ takes $f$ to the class on $M$ given by $f(F, \ldots, F),$ where $F$ is the curvature of the connection on $Q.$ This is well-defined and independent of choice of connection. (For details we refer the reader to \cite{Kobayashi:1969}.) In this case the transgression map $\tau$ is given by
$$
\tau(f) = \left( -\frac{1}{2} \right)^{k-1} \frac{k! (k - 1)!}{(2k-1)!} \, f(\Theta, [\Theta, \Theta], \ldots, [\Theta, \Theta]),
$$
where, as usual, $\Theta$ is the Maurer-Cartan form on $G.$ In terms of the result above, we have seen that in the case where the polynomial $f$ is given by $f(X,Y) = -\frac{1}{8\pi^2}\<X,Y\>$ and the $G$-bundle is $\widetilde{P} \to M\times S^1,$ then the Chern-Weil map gives the Pontrjagyn class of $\widetilde{P}$ and the diagram commutes. Furthermore, the element that fits in the bottom right hand corner is the string class of the corresponding $\Omega G$-bundle $P\to M.$ That is,
\begin{align*}
\int_{S^1} p_1(\widetilde{P}) &= s(P)\\
 &= -\frac{1}{4\pi^2}\int_{S^1} \<F, \nabla \Phi\> d\theta\\
 &=\frac{1}{48\pi^2}\hol_\Phi^*\<\Theta,[\Theta,\Theta]\>.
 \end{align*}

%

It is natural to ask now whether there is a similar theory for general and higher degree characteristic classes. That is, whether we can set up the following diagram
$$
\xymatrix{ H^{2k}(BG) \ar[r]^-{C\text{-}W_{\widetilde{P}}} \ar_\tau[d] & H^{2k}(M\times S^1) \ar^{\int_{S^1}}[d]\\
		H^{2k-1}(G) \ar[r]^{\hol_{\Phi}} & H^{2k-1}(M)}
$$
and give a formula for the element that ends up in the bottom right-hand corner given a general polynomial in the top left.

As above, the usual Chern-Weil theory tells us that if we start with an invariant polynomial $f \in H^{2k}(BG)$ then the element in $H^{2k}(M\times S^1)$ that we end up with is $f(\tilde{F},\ldots,\tilde{F})$ where $\tilde{F}$ is the curvature of the $G$-bundle $\widetilde{P}$ on $M\times S^1.$ Note that if we write out $f(\tilde{F},\ldots,\tilde{F})$ in terms of the curvature and Higgs field on the corresponding $\Omega G$-bundle $P \to M,$ we get
\begin{align*}
f(\tilde{F},\ldots,\tilde{F}) &= f(F+\nabla\Phi \, d\theta,\ldots,F+\nabla\Phi \, d\theta)\\
				&= f(F,\ldots, F) + k f(\nabla\Phi \, d\theta, F,\ldots,F)
\end{align*}
since $f$ is multilinear and symmetric and all terms with more than one $d\theta$ will vanish. From now on we will adopt the convention that whenever $f$ has repeated entries they will be ordered at the end and we will write them only once. That is, whatever appears as the last entry in $f$ is repeated however many times required to fill the remaining slots. (For example, $f(F) = f(F,\ldots, F)$ and $f(\nabla\Phi, F)d\theta = f(\nabla\Phi, F,\ldots, F)d\theta.$) So integrating this over the circle gives
$$
\int_{S^1} f(\tilde{F}) = k \int_{S^1} f(\nabla\Phi,F)\,d\theta.
$$
So $ k \int_{S^1} f(\nabla\Phi,F)d\theta$ is our candidate for the element in $H^{2k-1}(M)$ which corresponds to $f \in H^{2k}(BG)$ and makes the diagram commute. 

Note that if we evaluate this expression for the path fibration we have
\begin{align*}
k \int_{S^1} f( \nabla \Phi, F) \, d\theta &= f(\Theta, [\Theta, \Theta]) \left( \frac{1}{2} \right)^{k-1} k \int_{S^1} \left( \a^2 - \a \right)^{k-1} \partial \a \, d\theta\\
	&=  f(\Theta, [\Theta, \Theta]) \left( \frac{1}{2} \right)^{k-1} k \int_{S^1} \sum_{i=0}^{k-1} \binom{k-1}{i} (-1)^{k-1-i} \a^{2i} \a^{k-1-i} \partial \a \, d\theta\\
	&=  f(\Theta, [\Theta, \Theta]) \left( -\frac{1}{2} \right)^{k-1} k \sum_{i=0}^{k-1} \binom{k-1}{i} (-1)^{i} \frac{1}{k+i}.
\end{align*}
It turns out \cite{Sury:2004} that the coefficient above is equal to the coefficient in the definition of the transgression map $\tau$. That is,
$$
k \sum_{i=0}^{k-1} \binom{k-1}{i}  \frac{(-1)^i}{k+i} =  \frac{k! (k - 1)!}{(2k-1)!}.
$$
Therefore, we have for the path fibration
$$
k \int_{S^1} f( \nabla \Phi, F) \, d\theta = \tau (f).
$$
So what we are really asking for is a theory which associates to any characteristic class for $G$-bundles (that is, any polynomial in $H^{2k}(BG)$) a characteristic class for an $\Omega G$-bundle over $M.$ That is a map $H^{2k}(BG) \to H^{2k-1}(M)$ which gives characteristic classes for $\Omega G$-bundles over $M.$ Thus we need to show firstly that $k \int_{S^1} f(\nabla\Phi,F)d\theta$ is closed and independent of choice of connection and Higgs field. Also, we need to show that it is cohomologous to the pull-back by the classifying map $\hol_\Phi$ of the $(2k-1)$-form $\tau(f)$ defined above. We shall call  $ k \int_{S^1} f(\nabla\Phi,F)d\theta$ the \emph{string $(2k-1)$-class associated to $f$} and write $s_{2k-1}^P(f).$ To be more precise

\begin{definition}\label{D:string (2k-1)-class}
Let $\widetilde{P}$ be a framed $G$-bundle over $M\times S^1$ and $P$ its corresponding $\Omega G$-bundle over $M.$ Suppose that $f \in H^{2k}(BG)$ is an invariant polynomial representing the characteristic class $f(\tilde{F}) \in H^{2k}(M\times S^1).$ Then its associated \emph{string $(2k-1)$-class} is the class in $H^{2k-1}(M)$ given by
$$
s_{2k-1}^P(f) = k \int_{S^1} f(\nabla \Phi, F)\, d\theta,
$$
where $\Phi$ is a Higgs field for $P$ and $F$ is the curvature of a connection on $P.$
\end{definition}

Note that we still have to show that $s_{2k-1}^P(f)$ is closed and well-defined. We have

\begin{proposition}\label{P:closed}
The string $(2k-1)$-class is closed.
\end{proposition}

\begin{proof}

To show that  $s_{2k-1}^P(f)$ is closed we use the following result which follows from Lemmas 1 and 2 on pages 294--295 of \cite{Kobayashi:1969}:

\begin{lemma}\label{L:dpsi=Dpsi}
Let $\psi$ be an $ad$-invariant, vertical form on the total space of a principal bundle. Then $\psi$ projects to a form on the base space. For such a form, the exterior derivative is equal to the covariant exterior derivative. That is, $d\psi = D\psi.$
\end{lemma}

Thus we only need to show that $Ds_{2k-1}^P(f)=0.$ Now,
\begin{align*}
Dk\int_{S^1}f(\nabla\Phi, F) \, d\theta &= k\int_{S^1} f(D(\nabla\Phi), F) \, d\theta + k(k-1)\int_{S^1} f(\nabla\Phi, DF, F) \, d\theta\\
	&= k\int_{S^1} f(D(\nabla\Phi), F) \, d\theta
\end{align*}
using the Bianchi identity. We can calculate $D(\nabla\Phi).$ For tangent vectors $X$ and $Y,$
\begin{align*}
D(\nabla\Phi)(X,Y) &= d(\nabla\Phi)(hX,hY)\\
		&= (d^2\Phi + [dA,\Phi] - [A,d\Phi] - \partial (dA))(hX,hY)
\end{align*}
where $(hX,hY)$ is the projection of $(X,Y)$ onto the horizontal subspace at that point. Using the fact that $dA(hX,hY) = F(X,Y)$ and $A(hX)=A(hY)=0,$ we have:
$$
D(\nabla\Phi)(X,Y) = [F(X,Y),\Phi] - \partial F(X,Y)
$$
That is,
$$
D(\nabla\Phi) = [F,\Phi] - \partial F.
$$
So we have,
\begin{align*}
Dk\int_{S^1}f(\nabla\Phi, F) \, d\theta &= k\int_{S^1} f([F,\Phi],F) \, d\theta - k\int_{S^1} f(\partial F,F) \, d\theta\\
\end{align*}
and $ad$-invariance of $f$ (which we will discuss in more detail later) implies the first term on the right hand side vanishes while integration by parts implies the second term vanishes. Therefore, $s_{2k-1}^P(f)$ is closed.

\end{proof}

We also have

\begin{proposition}\label{P:independence}
The string $(2k-1)$-class is independent of choice of connection and Higgs field.
\end{proposition}

\begin{proof}

In order to see that $s_{2k-1}^P(f)$ is independent of choice of connection and Higgs field consider 2 different connection forms, $A_0$ and $A_1,$ on $P$ and 2 different Higgs fields, $\Phi_0$ and $\Phi_1.$ Since the space of connections is an affine space and the same is true for Higgs fields, we can consider lines joining the 2 connections and Higgs fields respectively. Define:
$$
\a := A_1 - A_0, \qquad  \qquad \varphi := \Phi_1 - \Phi_0
$$
and
$$
A_t := A_0 + t\a, \qquad \qquad \Phi_t := \Phi_0 + t\varphi
$$
for $t \in [0,1].$ Now consider the corresponding connection form on $\widetilde{P}$
\begin{multline*}
\tilde{A}_t = \tilde{A}_0 + t(\tilde{A}_1 - \tilde{A}_0)\\
\shoveleft{\phantom{\tilde{A}_t} = ad(g^{-1}) A_0 + \Theta + ad(g^{-1})\Phi_0 d\theta + t(ad(g^{-1}) A_1 + ad(g^{-1})\Phi_1 d\theta}\\
\shoveright{ - ad(g^{-1}) A_0 - ad(g^{-1})\Phi_0 d\theta)}\\
\shoveleft{\phantom{\tilde{A}_t}= ad(g^{-1}) A_0 + \Theta + ad(g^{-1})\Phi_0 d\theta + t\tilde{\a}}\\
\end{multline*}
where
$$
\tilde{\a} = ad(g^{-1})\a + ad(g^{-1})\varphi d\theta.
$$
Note that
$$
\tilde{A}_t = ad(g^{-1})A_t + \Theta + ad(g^{-1}) \Phi_t d\theta.
$$
Recall that $f(\tilde{F}) = f(F) + kf(\nabla \Phi, F)d\theta.$ We shall show $f(\tilde{F}_0)$ and $f(\tilde{F}_1)$ differ by an exact form, (where $\tilde{F}_0$ and $\tilde{F}_1$ are the curvature forms of $\tilde{A}_0$ and $\tilde{A}_1$ respectively) so that the class defined by $k\int_{S^1} f(\nabla \Phi_1, F_1)d\theta = \int_{S^1} f(\tilde{F})$ is independent of $A$ and $\Phi.$ For this we will need the following lemma:
\begin{lemma}\label{L:Da = d/dt(F)}
$\ds D_t \tilde{\a} = \frac{d}{dt}\tilde{F}_t.$
\end{lemma}

\begin{proof}
Firstly, we calculate $\tilde{F}_t$:
\begin{align*}
\tilde{F}_t &= d\tilde{A}_t + \tfrac{1}{2}[\tilde{A}_t,\tilde{A}_t]\\
		&= ad(g^{-1})\left( F_t + \nabla \Phi_t \wedge d\theta \right)\\
		&= ad(g^{-1})\left( dA_t  + \tfrac{1}{2}[A_t,A_t] + (d\Phi_t + [A_t,\Phi_t] - \partial A_t)\wedge d\theta \right)\\
		&= ad(g^{-1})\left( dA_0  + td\a + \tfrac{1}{2}[A_t,A_t] + (d\Phi_0 + td\varphi + [A_t,\Phi_t] - \partial A_0 - t\partial \a)\wedge d\theta \right).
\end{align*}
Therefore $\frac{d}{dt}\tilde{F}_t$ is given by
\begin{align*}
\frac{d}{dt}\tilde{F}_t &= ad(g^{-1})\left( d\a +\frac{1}{2}\frac{d}{dt}[A_t,A_t] +(d\varphi + \frac{d}{dt}[A_t,\Phi_t] - \partial\a)\wedge d\theta \right)\\
				&=  ad(g^{-1})\left( d\a +\frac{1}{2}[\a,A_t] + \frac{1}{2}[A_t, \a] +(d\varphi + [\a,\Phi_t] + [A_t,\varphi]- \partial\a)\wedge d\theta \right)\\
				&=  ad(g^{-1})\left( d\a +[\a,A_t]  +(d\varphi + [\a,\Phi_t] + [A_t,\varphi]- \partial\a)\wedge d\theta \right),
\end{align*}
since $\ds\frac{d}{dt}A_t = \a$ and $\ds\frac{d}{dt} \Phi_t =\varphi.$ Next we calculate $D_t\tilde{\a}$ by calculating $d\tilde{\a}$ and evaluating it on horizontal (with respect to $\tilde{A}_t$) vectors. At a point $(p,g, \theta)$ in $\widetilde{P}$ and for vectors $(X, g\xi, x_\theta)$ and $(Y, g\zeta, y_\theta)$ at $(p,g,\theta)$ we have:
\begin{multline*}
d\tilde{\a}_{(p,g,\theta)} (X, g\xi , x_\theta, Y, g\zeta, y_\theta)\\
=\tfrac{1}{2} \left\{ (X, g\xi, x_\theta)(\tilde{\a}_{(p,g,\theta)} (Y, g\zeta, y_\theta)) - (Y, g\zeta, y)(\tilde{\a}_{(p,g,\theta)} (X, g\xi, x_\theta))\right.\\
\left. -\tilde{\a}_{(p,g,\theta)} ([(X, g\xi, x_\theta), (Y, g\zeta, y_\theta)])\right\}.
\end{multline*}
So we need to calculate
\begin{enumerate}
\item $(X, g\xi, x_\theta)(\tilde{\a}_{(p,g,\theta)} (Y, g\zeta, y_\theta))$, and\\
\item $\tilde{\a}_{(p,g,\theta)} ([(X, g\xi, x_\theta), (Y, g\zeta, y_\theta)]).$
\end{enumerate}
If $\gamma_X (t)$ is a curve whose tangent vector is $X,$ we have:
\begin{multline*}
(X, g\xi, x_\theta)(\tilde{\a}_{(p,g,\theta)} (Y, g\zeta, y_\theta))\\
\shoveleft{=\frac{d}{dt}\bigg|_0\left\{ (1-t\xi)g^{-1}\a_{\gamma_X (t)}(Y)_{(\theta+tx)} g(1+t\xi) + (1-t\xi)g^{-1}\varphi_{\gamma_X (t),(\theta+tx)} g(1+t\xi)y\right\}}\\
\shoveleft{=\frac{d}{dt}\bigg|_0\left\{-t\xi g^{-1}\a_{\gamma_X (t)}(Y)_{(\theta+tx)} g +  g^{-1}\a_{\gamma_p (t)}(Y)_{(\theta+tx)} g t\xi + g^{-1}\a_{\gamma_X (t)}(Y)_{\theta} g\right.}\\
 + g^{-1} \partial \a_{\gamma_X (0)}(Y)_\theta x t g + -t\xi g^{-1}\varphi_{\gamma_X (t),(\theta+tx)} g y + g^{-1}\varphi_{\gamma_X (t),(\theta+tx)} g t\xi y\\
\shoveright{ \left.+ g^{-1} \varphi_{\gamma_X (t),\theta} g y + g^{-1}\partial \varphi_{\gamma_X (0), \theta}g t x y\right\}}\\
\shoveleft{= -\xi g^{-1} \a_p (Y)_\theta g + g^{-1} \a_p (Y)_\theta g\xi + g^{-1}\frac{d}{dt}\bigg|_0\a_{\gamma_X (t)} (Y)_\theta g + g^{-1}\partial \a_p(Y)_\theta g x}\\
-\xi g^{-1}\varphi_{p,\theta} gy + g^{-1}\varphi_{p,\theta} g\xi y + g^{-1}\frac{d}{dt}\bigg|_0\varphi_{\gamma_X (t), \theta} g y + g^{-1}\partial \varphi_{p, \theta} g xy.
\end{multline*}
Also,
\begin{align*}
\tilde{\a}_{(p,g,\theta)} ( [(X,g\xi, x_\theta), &(Y, g\zeta, y_\theta)]\\
	&= ad(g^{-1})\a_p([X,Y]) + ad(g^{-1})\varphi_p d\theta([x,y])\\
	&= ad(g^{-1})\a_p([X,Y])
\end{align*}
Therefore,
\begin{multline*}
d\tilde{\a}_{(p,g,\theta)} ( X, g\xi , x_\theta, Y, g\zeta, y_\theta)\\
= \frac{1}{2} \left\{ [ad(g^{-1})\a_p(Y), \xi] + ad(g^{-1})\left(\frac{d}{dt}\bigg|_0 \a_{\gamma_X (t)}(Y)_\theta\right) + ad(g^{-1})\partial\a_p(Y)x\right.\\
+ [ad(g^{-1})\varphi_{p,\theta}, \xi]y + ad(g^{-1}) \left(\frac{d}{dt}\bigg|_0 \varphi_{\gamma_p (t), \theta}\right)y\\
-[ad(g^{-1})\a_p(X), \zeta] - ad(g^{-1})\left(\frac{d}{dt}\bigg|_0 \a_{\gamma_X (t)}(X)_\theta\right) - ad(g^{-1})\partial\a_p(X)y\\
-[ad(g^{-1})\varphi_{p,\theta}, \zeta]x - ad(g^{-1}) \left(\frac{d}{dt}\bigg|_0 \varphi_{\gamma_X (t), \theta}\right)x\\
\left.- ad(g^{-1}) \a_p([X,Y])\vphantom{\frac{d}{dt}\bigg|_0} \right\}
\end{multline*}
That is,
\begin{multline*}
d\tilde{\a} = - [ad(g^{-1})\a,\Theta] + ad(g^{-1}) d\a - ad(g^{-1})\partial \a \wedge d\theta\\
+ [ad(g^{-1})\varphi, \Theta]\wedge d\theta + ad(g^{-1}) d\varphi \wedge d\theta
\end{multline*}
\begin{multline*}
\phantom{d\tilde{\a}} = ad(g^{-1}) \left( d\a + d\varphi \wedge d\theta - \partial \a \wedge d\theta\right) - [ad(g^{-1})\a + ad(g^{-1})\varphi d\theta, \Theta]\\ {}
\end{multline*}
\vspace{-5ex}

To calculate $D_t \tilde{\a}$ we need to know what the horizontal projection (with respect to $\tilde{A}_t$) of a vector looks like. If $X$ is a tangent vector at $p$ we can calculate its horizontal projection as $hX = X - \iota_p ( A(X)),$ where $\iota_p ( A(X))$ is the vector at $p$ generated by the Lie algebra element $A(X).$ So for the vector $(X, g\xi, x_\theta)$ we have
$$
h(X, g\xi, x_\theta) = (X, g\xi, x_\theta) - \iota_{(p, g, \theta)} ( \tilde{A}_t (X, g\xi, x_\theta)).
$$
Now,
\begin{align*}
\iota_{(p, g, \theta)} ( \tilde{A}_t (X, g\xi, x_\theta)) &= \frac{d}{ds}\bigg|_0 (p,g(1+s  \tilde{A}_t (X, g\xi, \theta + x)), \theta)\\
	&=  \frac{d}{ds}\bigg|_0 (p,gs  \tilde{A}_t (X, g\xi, \theta + x), \theta)\\
	&= (0, g \tilde{A}_t (X, g\xi, x_\theta), 0),
\end{align*}
and therefore,
$$
h(X, g\xi, x_\theta) = (X, g(\xi - \tilde{A}_t (X, g\xi, x_\theta)), x_\theta).
$$
Putting this into the formula above for $d\tilde{\a},$ we obtain
$$
D_t\tilde{\a} = ad(g^{-1}) \left( d\a + d\varphi \wedge d\theta - \partial \a \wedge d\theta\right) - [ad(g^{-1})\a + ad(g^{-1})\varphi d\theta, \Theta - \tilde{A}_t]
$$
and inserting the formula for $\tilde{A}_t$ in terms of $A_t$ and $\Phi_t$ in the second term we obtain
\begin{align*}
 - [ad(g^{-1})\a + &ad(g^{-1})\varphi d\theta, \Theta - \tilde{A}_t]\\
 	&= - [ad(g^{-1})\a + ad(g^{-1})\varphi d\theta, \Theta - ad(g^{-1})A_t - \Theta - ad(g^{-1}) \Phi_t d\theta]\\
 	&= - [ad(g^{-1})\a + ad(g^{-1})\varphi d\theta, - ad(g^{-1})A_t  - ad(g^{-1}) \Phi_t d\theta]\\
	&= ad(g^{-1})[\a + \varphi d\theta,  A_t  + \Phi_t d\theta]
\end{align*}
and therefore
$$
D_t\tilde{\a} = ad(g^{-1}) \left( d\a + d\varphi \wedge d\theta - \partial \a \wedge d\theta + [\a, A_t] + [\a, \Phi_t]d\theta + [A_t,\varphi]d\theta\right)
$$
which is equal to $\dfrac{d}{dt}\tilde{F}_t.$ This completes the proof of Lemma \ref{L:Da = d/dt(F)}.
\renewcommand{\qedsymbol}{}
\end{proof}

Now, if we set
$$
\psi = k\int_0^1 f(\tilde{\a},\tilde{F}_t)dt
$$
then
\begin{align*}
d\psi & = D\psi \qquad \text{ (by Lemma \ref{L:dpsi=Dpsi}) }\\
	& = k\int_0^1 f(D_t\tilde{\a},\tilde{F}_t)dt\\
	& = k\int_0^1 f(\frac{d}{dt}\tilde{F}_t,\tilde{F}_t)dt\\
	& = \int_0^1 \frac{d}{dt} f(\tilde{F}_t)dt\\
	& = f(\tilde{F}_1) - f(\tilde{F}_0).
\end{align*}
So $s_{2k-1}^P(f)$ is independent of choice of connection and Higgs field.

\end{proof}

It remains only to prove that $s_{2k-1}^P(f)$ is the pull-back of $\tau(f)$ by $\hol_\Phi.$ For this we follow the argument in \cite{Carey:1991} that will give us a formula for $f(\tilde{F})$ that we can use to calculate $s_{2k-1}^{\Omega \widetilde{P}}(f)$ for a loop bundle $\Omega \widetilde{P} \xrightarrow{\Omega G} \Omega (M\times S^1)$ and then we can use Lemma \ref{L:holphi=holeta} to generalise to a general $\Omega G$-bundle.

If we start with the $G$-bundle $\widetilde{P} \to M\times S^1$ we can pull-back by the evaluation map $\ev \colon [0,1] \times \Omega (M\times S^1) \to  (M\times S^1)$ to get a trivial bundle $\ev^*\widetilde{P}$ over $[0,1] \times \Omega (M\times S^1).$ A section is given by
$$
h \colon [0,1]\times \Omega  (M\times S^1) \to \ev^*\widetilde{P}; \quad (t,\gamma) \mapsto  \hat{\gamma}(t),
$$
where $\hat{\gamma}$ is the horizontal lift of $\gamma.$ If $\tilde{A}$ is the connection in $\widetilde{P}$ we can pull it back to $\ev^*\widetilde{P}$ and then back to $[0,1] \times \Omega (M\times S^1)$ to obtain
$$
\tilde{A}' := h^*\ev^*\tilde{A}.
$$
We can calculate the curvature $\tilde{F}$ of $\tilde{A}$ and pull it back by $\ev$ to $[0,1] \times \Omega (M\times S^1)$ and because this is a product manifold we can decompose it into parts with a $dt$ and parts without a $dt.$ Under this decomposition, we have
$$
\ev^*\tilde{F} = -\frac{\partial}{\partial t} \tilde{A}'\wedge dt + \tilde{F}',
$$
where we call the component without a $dt$ $\tilde{F}'$ since if we view the form $\tilde{A}'$ for fixed $t_0$ as a connection form on $\Omega (M\times S^1)$ then its curvature is $\tilde{F}'$ evaluated at $t_0.$ 

Now, we want to calculate $\int_{S^1} f(\tilde{F})$ and using Lemma \ref{L:eta* int ev* = int} we have for a general $\Omega G$-bundle $P\to M,$
\begin{align*}
\int_{S^1} f(\tilde{F}) &= \eta^*\int_{S^1}\ev^*f(\tilde{F})\\
				&= \eta^*\int_{S^1}f(\ev^*\tilde{F}).				
\end{align*}
So we wish to calculate explicitly $\int_{S^1}f(\ev^*\tilde{F}).$ If we view the circle as the interval $[0,1]$ with endpoints identified, then we can write
$$
\int_{S^1} f(\ev^*\tilde{F}) = \int_{[0,1]} f(\ev^*\tilde{F})
$$
and so we have
\begin{align*}
k\int_{S^1} f(\nabla\Phi, F)d\theta & = \eta^*\int_{S^1} f(\ev^*\tilde{F})\\
						&= \eta^*\int_{[0,1]} f(- \frac{\partial}{\partial t} \tilde{A}'\wedge dt + \tilde{F}')\\
						&= \eta^*\int_{[0,1]} f(\tilde{F}') -k \eta^*\int_{[0,1]} f(- \frac{\partial}{\partial t} \tilde{A}', \tilde{F}') dt\\
						&= -k \eta^*\int_{[0,1]} f(- \frac{\partial}{\partial t} \tilde{A}', \tilde{F}') dt.
\end{align*}
Using the formula $\tilde{F}' = d\tilde{A}' + \frac{1}{2} [\tilde{A}', \tilde{A}'],$ we can write this as:
\begin{multline*}
-k \eta^*\left\{ \int_{[0,1]} f(\partial \tilde{A}',d\tilde{A}')dt\right. \\
\left.+ (k-1)\frac{1}{2} \int_{[0,1]} f(\partial\tilde{A}', d\tilde{A}', \ldots, d\tilde{A}',[\tilde{A}',\tilde{A}'])dt +\ldots \right.\\
...+\binom{k-1}{k-2} \left(\frac{1}{2}\right)^{k-2}\int_{[0,1]} f(\partial\tilde{A}', d\tilde{A}', [\tilde{A}',\tilde{A}'])dt\\
\left.+ \left(\frac{1}{2}\right)^{k-1} \int_{[0,1]} f(\partial\tilde{A}', [\tilde{A}',\tilde{A}'])dt \right\}
\end{multline*}
where we have written $\partial \tilde{A}'$ for $\partial\tilde{A}' / \partial t.$ Thus we need to work with the general term
$$
\binom{k-1}{i} \left(\frac{1}{2}\right)^{i}\int_{[0,1]} f(\partial \tilde{A}', \underbrace{d\tilde{A}', \ldots, d\tilde{A}'}_{k-i-1},  \underbrace{[\tilde{A}',\tilde{A}'], \ldots,  [\tilde{A}',\tilde{A}']}_{i})dt.
$$
To deal with these terms we shall use integration by parts  and the $ad$-invariance of $f.$ Thus we need to know in detail how $ad$-invariance works.

\begin{lemma}\label{L:ad-invariance}
Let $\varphi_1, \ldots, \varphi_k$  be $\fg$-valued forms of degree $q_1, \ldots, q_k$ respectively. Then if $A$ is a $\fg$-valued $p$-form, we have
\begin{multline*}
f([\varphi_1, A], \varphi_2, \ldots, \varphi_k)\\
= f(\varphi_1, [A, \varphi_2], \ldots, \varphi_k) + (-1)^{p q_2} f(\varphi_1, \varphi_2, [A, \varphi_3],\ldots, \varphi_k) +\ldots\\
 \ldots +  (-1)^{p(q_2 + \ldots q_{k-1})} f(\varphi_1, \ldots,\varphi_{k-1}, [A,\varphi_k]).
\end{multline*}
\end{lemma}

\begin{proof}

We can expand $\varphi_i$ as $\varphi_i = \varphi_{i,j} \omega_i^j$ for $\varphi_{i,j} \in \fg$ and $\omega_i^j$ a $q_i$-form. Then we have
$$
f(\varphi_1, \ldots, \varphi_k) = f(\varphi_{1,j_1}, \ldots, \varphi_{k, j_k}) \omega_1^{j_1}\wedge \ldots \wedge \omega_k^{j_k}.
$$
Now if $A$ is a $\fg$ valued $p$-form and we write $A = A_i \a^i$ as above, then 
\begin{multline*}
f([A,\varphi_1], \varphi_2,\ldots,\varphi_k)\\
\shoveleft{\phantom{f}=f([A_i,\varphi_{1,j_1}], \varphi_{2,j_2},\ldots, \varphi_{k,j_k}) \a^i \wedge \omega_1^{j_1} \wedge \ldots \wedge \omega_k^{j_k}}\\
\shoveleft{\phantom{f} = f(\varphi_{1,j_1}, [\varphi_{2,j_2},A_i], \ldots, \varphi_{k,j_k})(-1)^{p(q_1 + q_2)} \omega_1^{j_1} \wedge  \omega_2^{j_2} \a^i \wedge\ldots \wedge \omega_k^{j_k}}\\
+ f(\varphi_{1,j_1},\varphi_{2,j_2}, [\varphi_{3,j_3},A_i], \ldots, \varphi_{k,j_k})(-1)^{p(q_1 + q_2 + q_3)} \omega_1^{j_1} \wedge  \omega_2^{j_2} \wedge  \omega_3^{j_3}\wedge\a^i \wedge\ldots \wedge \omega_k^{j_k}\\
\ldots + f(\varphi_{1,j_1},\varphi_{2,j_2}, \ldots, [\varphi_{k,j_k},A_i])(-1)^{p(q_1 + q_2 +\ldots + q_k)} \omega_1^{j_1} \wedge  \omega_2^{j_2} \wedge\ldots \wedge \omega_k^{j_k}\wedge\a^i
\end{multline*}
That is,
\begin{multline*}
f([A,\varphi_1], \varphi_2,\ldots,\varphi_k)\\
= (-1)^{pq_1} f(\varphi_1, [\varphi_2,A], \ldots,\varphi_k) + (-1)^{p(q_1 + q_2)} f(\varphi_1, \varphi_2, [\varphi_3,A], \ldots,\varphi_k) + \ldots\\
\ldots + (-1)^{p(q_1 + \ldots + q_k)} f(\varphi_1, \ldots,\varphi_{k-1},[\varphi_k, A]),
 \end{multline*}
which we can write as:
\begin{multline*}
f([\varphi_1, A], \varphi_2, \ldots, \varphi_k)\\
= f(\varphi_1, [A, \varphi_2], \ldots, \varphi_k) + (-1)^{p q_2} f(\varphi_1, \varphi_2, [A, \varphi_3],\ldots, \varphi_k) +\ldots\\
 \ldots +  (-1)^{p(q_2 + \ldots q_{k-1})} f(\varphi_1, \ldots,\varphi_{k-1}, [A,\varphi_k]).
\end{multline*}
 
\end{proof}

We are now in a position to prove

\begin{proposition}\label{P:s(P)=hol^*f}
$$
s_{2k-1}^P(f)= \hol_\Phi ^*\tau(f).
$$
\end{proposition}

\begin{proof}

To calculate the general term given above, we integrate by parts in the $\Omega(M\times ~S^1)$ and $t$ directions giving
\begin{multline*}
\int_{[0,1]}f_i dt = \int_{[0,1]} f(d\partial \tilde{A}', \tilde{A}', d\tilde{A}', \ldots, d\tilde{A}', [\tilde{A}', \tilde{A}'])dt\\
 + i \int_{[0,1]} f(\partial \tilde{A}', \tilde{A}', d\tilde{A}', \ldots, d\tilde{A}', d[\tilde{A}', \tilde{A}'], [\tilde{A}', \tilde{A}']) dt\\
 - d\int_{[0,1]} f(\partial \tilde{A}', \tilde{A}', d\tilde{A}', \ldots, d\tilde{A}', [\tilde{A}', \tilde{A}'])dt
\end{multline*}
and
\begin{multline*}
\int_{[0,1]}f_i dt = f(\tilde{A}'_1, d\tilde{A}'_1, \ldots, d\tilde{A}'_1, [\tilde{A}'_1, \tilde{A}'_1]) - f(\tilde{A}'_0, d\tilde{A}'_0, \ldots, d\tilde{A}'_0, [\tilde{A}'_0, \tilde{A}'_0])\\
 - (k-1-i)\int_{[0,1]} f(\tilde{A}', \partial d\tilde{A}', d\tilde{A}', \ldots, d\tilde{A}', [\tilde{A}', \tilde{A}'])dt\\
 - i \int_{[0,1]} f(\tilde{A}', d\tilde{A}', \ldots, d\tilde{A}', \partial [\tilde{A}', \tilde{A}'], [\tilde{A}', \tilde{A}']) dt
\end{multline*}
where we have written $f_i$ for the integrand of the general term given earlier. Combining these gives
\begin{multline*}
(k-i)\int_{[0,1]} f_i dt = f_{i,1} - f_{i,0} - i\int_{[0,1]} f(\tilde{A}', d\tilde{A}', \ldots, d\tilde{A}', \partial[\tilde{A}', \tilde{A}'], [\tilde{A}', \tilde{A}']) dt\\
+ i(k-1-i) \int_{[0,1]} f(\partial \tilde{A}', \tilde{A}', d\tilde{A}', \ldots, d\tilde{A}', d[\tilde{A}', \tilde{A}'], [\tilde{A}', \tilde{A}'])dt\\
- (k-1-i) d\int_{[0,1]} f(\partial \tilde{A}', \tilde{A}', d\tilde{A}', \ldots, d\tilde{A}', [\tilde{A}', \tilde{A}'])dt
\end{multline*}
where we have written $f_{i,1}$ and $f_{i,0}$ for $f_i$ evaluated at $t=1$ and $0$ respectively. Using $ad$-invariance, the term on the middle line simplifies as follows:
\begin{multline*}
\int_{[0,1]} f(\partial \tilde{A}', \tilde{A}', d\tilde{A}', \ldots, d\tilde{A}', d[\tilde{A}', \tilde{A}'], [\tilde{A}', \tilde{A}'])dt\\
\shoveleft{= 2\int_{[0,1]} f([d\tilde{A}', \tilde{A}'], \partial \tilde{A}', \tilde{A}', d\tilde{A}', \ldots, d\tilde{A}', [\tilde{A}', \tilde{A}'])dt}\\
\shoveleft{= 2 \int_{[0,1]} f(d\tilde{A}', [\tilde{A}', \partial\tilde{A}'], \tilde{A}', d\tilde{A}', \ldots, d\tilde{A}', [\tilde{A}', \tilde{A}'])dt}\\
 - 2 \int_{[0,1]} f( d\tilde{A}', \partial \tilde{A}',[\tilde{A}', \tilde{A}'], d\tilde{A}', \ldots, d\tilde{A}', [\tilde{A}', \tilde{A}'])dt\\
\shoveright{ + 2(k-2-i)  \int_{[0,1]} f( d\tilde{A}', \partial \tilde{A}', \tilde{A}', [\tilde{A}', d\tilde{A}'],d\tilde{A}',  \ldots, d\tilde{A}', [\tilde{A}', \tilde{A}'])dt}\\
\shoveleft{ = \int_{[0,1]} f(d\tilde{A}', \partial[\tilde{A}',\tilde{A}'], \tilde{A}', d\tilde{A}', \ldots, d\tilde{A}', [\tilde{A}', \tilde{A}'])dt}\\
- 2 \int_{[0,1]} f(\partial \tilde{A}', d\tilde{A}', \ldots, d\tilde{A}', [\tilde{A}', \tilde{A}'])dt\\
- (k-2-i)  \int_{[0,1]} f(\partial \tilde{A}', \tilde{A}',d\tilde{A}',  \ldots, d\tilde{A}', d[\tilde{A}', \tilde{A}'], [\tilde{A}', \tilde{A}'])dt
\end{multline*}
and so
\begin{multline*}
(k-1-i)  \int_{[0,1]} f(\partial \tilde{A}', \tilde{A}',d\tilde{A}',  \ldots, d\tilde{A}', d[\tilde{A}', \tilde{A}'], [\tilde{A}', \tilde{A}'])dt\\
= \int_{[0,1]} f(\tilde{A}', d\tilde{A}', \ldots, d\tilde{A}', \partial[\tilde{A}',\tilde{A}'],[\tilde{A}', \tilde{A}'])dt\\
- 2 \int_{[0,1]} f(\partial \tilde{A}', d\tilde{A}', \ldots, d\tilde{A}', [\tilde{A}', \tilde{A}'])dt.
\end{multline*}
Inserting this into the formula for $\int f_i dt$ gives
\begin{multline*}
(k-i)\int_{[0,1]} f_i dt = f_{i,1} - f_{i,0} -2i\int_{[0,1]} f_i dt\\
- (k-1-i) d\int_{[0,1]} f(\partial \tilde{A}', \tilde{A}', d\tilde{A}', \ldots, d\tilde{A}', [\tilde{A}', \tilde{A}'])dt
\end{multline*}
and hence
\begin{multline*}
(k+i)\int_{[0,1]} f_i dt\\
 = f_{i,1} - f_{i,0} - (k-1-i) d\int_{[0,1]} f(\partial \tilde{A}', \tilde{A}', d\tilde{A}', \ldots, d\tilde{A}', [\tilde{A}', \tilde{A}'])dt.
\end{multline*}

So we have the following expression for $s_{2k-1}^P(f):$
\begin{multline*}
k\int_{S^1} f(\nabla\Phi, F)d\theta\\
= -k \eta^*\left\{ \sum_{i=0}^{k-1} \binom{k-1}{i} \left(\frac{1}{2}\right)^{i} \frac{1}{k+i} \left( f_{i, 1} - f_{i, 0} -  (k-i-1) d c_i \vphantom{\tilde{f}} \right) \right\}
\end{multline*}
where $c_i$ is the last integral in the equation above (with $i$ $[\tilde{A}', \tilde{A}']$'s).

Now since $\tilde{A}'_0 = 0$ and $h(0, \gamma) = h(1, \gamma)\hol(\gamma)$ (where $h$ is the section from earlier), we have that 
$$
\tilde{A}'_0 = ad(\hol^{-1}) \tilde{A}'_1 + \hol^{-1} d\hol
$$
and so
$$
\tilde{A}'_1 = - d\hol \hol^{-1}.
$$
Therefore we have that $f_{i,0} = 0$ and we can calculate $f_{i,1}$ in terms of $f_{0,1}$ as follows:
\begin{align*}
f_{0,1} 	&= f(\tilde{A}'_1, d\tilde{A}'_1)\\
		&= f(- d\hol \hol^{-1}, d(- d\hol \hol^{-1}))\\
		&= (-1)^{k} \left( \frac{1}{2}\right)^{k-1} \hol^* f(\Theta, [\Theta, \Theta])
\end{align*}
and in general,
\begin{align*}
f_{i,1}	&= f(\tilde{A}'_1, d\tilde{A}'_1, \ldots, d\tilde{A}'_1, [\tilde{A}'_1,\tilde{A}'_1])\\
		&= (-1)^{k-i} \left(\frac{1}{2}\right)^{k-1-i} \hol^* f(\Theta, [\Theta, \Theta])\\
		&= (-1)^{i} 2^{i} f_{0,1}
\end{align*}
using the fact that $d(-d\hol\hol^{-1}) = -\frac{1}{2} [d\hol \hol^{-1}, d\hol \hol^{-1}].$

Therefore we have
\begin{multline*}
k\int_{S^1} f(\nabla\Phi, F)d\theta\\
=  \left(-\frac{1}{2}\right)^{k-1} k \sum_{i=0}^{k-1} \binom{k-1}{i}  \frac{(-1)^i}{k+i} \hol_\Phi^*f(\Theta, [\Theta, \Theta])\\
+ k\sum_{i=0}^{k-i}\binom{k-1}{i} \left(\frac{1}{2}\right)^i  \frac{1}{k+i}   (k-i-1) d c_i.
\end{multline*}

We have seen already that the coefficient above is equal to the coefficient in the definition of the transgression map: 
$$
k \sum_{i=0}^{k-1} \binom{k-1}{i}  \frac{(-1)^i}{k+i} =  \frac{k! (k - 1)!}{(2k-1)!}.
$$
So we see that the pull-back of the transgression of $f$ is cohomologous to the string $(2k-1)$-class.

\end{proof}

Combining Propositions \ref{P:closed}, \ref{P:independence} and \ref{P:s(P)=hol^*f}, we have the following Theorem

\begin{theorem}\label{T:"Chern-Weil"}
The diagram
$$
\xymatrix{ H^{2k}(BG) \ar[r]^-{C\text{-}W_{\widetilde{P}}} \ar_\tau[d] & H^{2k}(M\times S^1) \ar^{\int_{S^1}}[d]\\
		H^{2k-1}(G) \ar^{\hol_\Phi^*}[r] & H^{2k-1}(M)}
$$
commutes. Furthermore, the composition map
$$
H^{2k}(BG) \to H^{2k-1}(M)
$$
associates to any invariant polynomial its string $(2k-1)$-class, which is a characteristic class.




\end{theorem}

\section{The universal string class for $L^{\ssv} G$-bundles}\label{S:universal string class}

We would now like to return to the study of the free loop group. In this section, we shall give a partial generalisation of the results in the previous section. However, we shall be working with a slightly different group than in the rest of this thesis. For the remainder of this chapter we shall be considering the group of smooth maps from the interval $[0,2\pi]$ into $G$ whose endpoints are coincident. This group shall be denoted $L^{\ssv} G.$ Note that $LG \subseteq L^{\ssv}G.$ We also have the based version $\Omega^{\ssv} G$ of this group consisting of maps $[0, 2\pi] \to G$ such that the endpoints are mapped to the identity in $G.$

We will give a classifying theory for $L^{\ssv} G$ bundles and present a calculation for the string class of the universal $L^{\ssv} G$-bundle.

\subsection{Classification of $L^{\ssv} G$-bundles}\label{SS:LG classification}

In order to extend the ideas from the previous section (namely, calculating the string class of the universal $L^{\ssv} G$-bundle) we need a model for $EL^{\ssv} G.$ To construct this we view $L^{\ssv} G$ as the semi-direct product $\Omega^{\ssv} G\rtimes G.$ The group multiplication is given by
$$
(\gamma_1, g_1) (\gamma_2, g_2) = (g_2^{-1} \gamma_1 g_2 \gamma_2 , g_1g_2)
$$
and the isomorphism between $\Omega^{\ssv} G \rtimes G$ and $L^{\ssv} G$ is
\begin{align*}
	\Omega^{\ssv} G \rtimes G &\xrightarrow{\sim} L^{\ssv} G;\quad
	(\gamma, g)		\mapsto g\gamma.
\end{align*}
On the level of Lie algebras, the isomorphism is
\begin{align*}
	\Omega^{\ssv} \fg \rtimes \fg &\xrightarrow{\sim} L^{\ssv} \fg;\quad
	(\xi, X)		\mapsto X + \xi.
\end{align*}
We therefore need a model for the universal $\Omega^{\ssv} G \rtimes G$-bundle. For this, we shall take the product of the universal $\Omega^{\ssv} G$-bundle and the universal $G$-bundle. A model for the universal $\Omega^{\ssv} G$-bundle is given by the space of maps from the interval $[0, 2\pi]$ into $G,$ denoted $P^{\ssv} G.$ 
The based loop group $\Omega^{\ssv} G$ acts on this space by right multiplication and evaluation at the endpoint of a path gives a locally trivial $\Omega^{\ssv} G$-bundle $P^{\ssv} G \to G.$ As our study of $\Omega^{\ssv} G$ will be confined to this section, we shall refer to $P^{\ssv} G$ as the \emph{path fibration} without any risk of confusion. $P^{\ssv} G$ is contractible since any path $p$ can be homotopied to the identity path by the map
$$
h \colon I \times P^{\ssv} G \to P^{\ssv} G; \quad (t, p) \mapsto (\theta \mapsto p(t\theta)).
$$
Therefore the path fibration is a model for the universal $\Omega^{\ssv} G$-bundle. So, for our model for $EL^{\ssv} G$ we shall take the space $P^{\ssv} G \times EG$ which is contractible since $P^{\ssv} G$ and $EG$ are both contractible. This is acted on by $\Omega^{\ssv} G\rtimes G:$
$$
(p, x) (\gamma, g) = (g^{-1} p g \gamma, xg)
$$
where $xg$ is the right action of $G$ on $EG.$ This action is free (since $G$ acts on $EG$ freely) and transitive on fibres (since the action on $EG$ is transitive and the equation $g^{-1}p_1 g\gamma = p_2$ can always be solved) and so $P^{\ssv} G \times EG$ is a model for $EL^{\ssv} G$ and $BL^{\ssv} G$ is equal to $(P^{\ssv} G\times EG)/(\Omega^{\ssv} G \rtimes G).$ In fact, if we consider the map
$$
(P^{\ssv} G\times EG)/(\Omega^{\ssv} G \rtimes G) \to (G\times EG)/G; \quad [p, x] \mapsto [p(2\pi), x],
$$
where $[h, x] = [g^{-1}hg, xg],$ we can see this is well-defined, since
$$
[p, x] = [g^{-1}pg\gamma, xg] \mapsto [g^{-1}p(2\pi)g \gamma(2\pi), xg] = [p(2\pi), x].
$$
Furthermore, this is onto, as the projection $P^{\ssv} G \to G$ is onto, and 1--1, for if we consider two elements $[p,x],\,[q,y] \in(P^{\ssv} G\times EG)/(\Omega^{\ssv} G \rtimes G)$ such that $[p(2\pi), x] = [q(2\pi), y]$ we have $y = xg$ and $q(2\pi) = g^{-1}p(2\pi)g.$ That is, the paths $q$ and $g^{-1}p g$ have the same endpoint. Therefore, the path $g^{-1}p^{-1}gq$ is actually a (based) loop. And since $q = g^{-1}pg (g^{-1}p^{-1}gq),$ we have
\begin{align*}
[q, y]	&= [ g^{-1}pg \gamma, xg]\\
	&= [p, x],
\end{align*}
where $\gamma=  g^{-1}p^{-1}gq \in \Omega^{\ssv} G.$ Thus we have a diffeomorphism between $BL^{\ssv} G$ and $(G\times EG)/G$ (or simply $G\times_G EG$). Note that this allows us to calculate the cohomology of $BL^{\ssv} G$ as the equivariant cohomology of $G$ (with its adjoint action). That is,
$$
H(BL^{\ssv} G) = H_G (G).
$$

Given an $L^{\ssv} G$-bundle $P\to M$ we can write down the classifying map of this bundle as follows. Choose a Higgs field, $\Phi,$ for $P.$ Then define the map $f \colon P \to P^{\ssv} G \times EG$ by 
$$
f(q) = (\hol_\Phi (q) , f_G (q)),
$$
where $\hol_\Phi$ is the Higgs field holonomy and $f_G$ is the classifying map for the $G$-bundle associated to $P$ by the projection $L^{\ssv} G \to G$ given by mapping a loop to its start/endpoint (or equivalently, the projection $\Omega^{\ssv} G\rtimes G \to G$). That is, $f(q) = (p, x)$ where $p^{-1} \partial p  = \Phi(q) $ and $ x $ is $f_G$ applied to the image of $q$ in $P\times_{L^{\ssv} G}G.$ It is easy to see that this is equivariant with respect to the $L^{\ssv} G$ action and hence descends to a map $M\to BL^{\ssv} G$ since if $(\gamma, g) \in \Omega^{\ssv} G \rtimes G$ then
\begin{align*}
f(q(g\gamma))	&= (\hol_\Phi (q(g\gamma)), f_G(q)g)
\end{align*}
and so $f$ is equivariant in the $EG$ slot (by virtue of the fact that $f_G$ is a classifying map) and also in the $P^{\ssv} G$ slot since if $\hol_\Phi(q) = p$ then
\begin{align*}
\Phi(q(g\gamma)) 	&= ad((g\gamma)^{-1})\Phi(q) + (g\gamma)^{-1}\partial (g\gamma)\\
				&= ad((g\gamma)^{-1})\Phi(q) + \gamma^{-1}\partial \gamma
\end{align*}
and
\begin{align*}
(p(\gamma, g))^{-1} \partial (p(\gamma, g)) &= (g^{-1}p g\gamma)^{-1}\partial(g^{-1}p g\gamma)\\
			&= \gamma^{-1}g^{-1}p^{-1}g (g^{-1} \partial p g \gamma + g^{-1}p g \partial \gamma)\\
			&= ad((g\gamma)^{-1}) p^{-1}\partial p + \gamma^{-1}\partial \gamma
\end{align*}
and so $\hol_\Phi (q(g\gamma)) = p (\gamma, g) = \hol_\Phi (q)(g\gamma).$


\subsection{The universal string class}\label{SS:universal string class}

Now that we have a model for the universal $L^{\ssv} G$-bundle we would like to calculate its string class according to Theorem \ref{T:LGstring}. So far everything we have said works on the topological level. In order to use Theorem \ref{T:LGstring} however, the first thing we need is a connection on $P^{\ssv} G \times EG.$ Now, $P^{\ssv} G$ is already a smooth manifold. In order to define a smooth structure and find a connection on $EG$ we use the results in \cite{Narasimhan:1961, Narasimhan:1963}. As long as the dimension of the base of the $G$-bundle $P \to M$ is less than or equal to $n$ this gives a construction of a smooth bundle $EG_n \to BG_n$ with connection which is a model for the universal $G$-bundle. From now on we assume therefore that the dimension of the base of our $L^{\ssv} G$-bundle is fixed (and less than or equal to $n$ for some $n$).

 To define a connection we need to know what a vertical vector looks like. Consider the vector in $T_{(p, x)}(P^{\ssv} G\times EG_n) = T_p P^{\ssv} G \times T_x EG_n$ generated by the Lie algebra element $(\xi, X) \in \Omega^{\ssv} \fg \rtimes \fg:$
\begin{align*}
\iota_{(p,x)} (\xi, X) 	&=  \frac{d}{dt}\bigg|_0 ( (1-tX) p (1+tX)(1+t\xi), x e^{tX})\\
				&=  \frac{d}{dt}\bigg|_0 ( t(-Xp + pX + p\xi), x e^{tX})\\
				&= ( p (X - ad(p^{-1})X + \xi), \iota_x (X)).
\end{align*}
Note that the $P^{\ssv} G$ part of a vertical vector is a vector field along $p$ that ends at $p(2\pi) (X - ad(p(2\pi)^{-1})X)$ (since $\xi$ is a based loop). We will assume that we have a connection in $EG_n$ since this is always possible by the discussion above. Call this connection $a.$ So to find the horizontal part of a vector $(V,W) \in T_p P^{\ssv} G \times T_x EG_n$ we need a vector field along $p$ that ends at $V(2\pi) - p(2\pi) (X - ad(p(2\pi)^{-1})X)$ (since then $V - \{\text{this vector}\}$ will end at the right point to be vertical). Consider the vector field
$$
\left(\frac{\theta}{2\pi}\right) p \left\{ad(p^{-1}) \left(V(2\pi)p(2\pi)^{-1} - ad(p(2\pi))a(W) + a(W)\right) \right\}.
$$
If we define the horizontal projection of $(V,W), h(V,W),$ to be the vector field above together with the horizontal component of $W$ (that is, $hW = W - \iota_x(a(W))$), then we have an invariant splitting of the tangent space at each point in $P^{\ssv} G \times EG_n.$ This is easily verified: Since the $EG_n$ part has a connection, we need only check the $P^{\ssv} G$ part. First calculate the right action on the vector above (which we will call $hV$ even though technically the part of the connection on $P^{\ssv} G$ is not actually a connection itself):
\begin{multline*}
\left(hV(\gamma, g)\right)_{(g^{-1}pg\gamma, xg)}  \\
 = \left(\frac{\theta}{2\pi}\right)g^{-1} p \left\{ad(p^{-1}) \left(V(2\pi)p(2\pi)^{-1} - ad(p(2\pi))a(W) + a(W)\right) \right\}g\gamma.
\end{multline*}
Compare this with the horizontal projection of a vector $V'$ at $(p, x) (\gamma, g)= (g^{-1}pg\gamma, xg):$
\begin{multline*}
hV'_{(g^{-1}pg\gamma, xg)}  \\
\shoveleft{\phantom{hV'_{(g^{-1}pg\gamma}} = \left(\frac{\theta}{2\pi}\right) g^{-1}pg\gamma \left\{ad(g^{-1} p^{-1}g\gamma)^{-1} \left(V'(2\pi)g^{-1}p(2\pi)^{-1}g\right.\right.}\\
\shoveright{ \left.\left.- ad(g^{-1}p(2\pi)g)a(W') + a(W')\right)\right\}}\\
\shoveleft{\phantom{hV'_{(g^{-1}pg\gamma}} = \left(\frac{\theta}{2\pi}\right) g^{-1}p \left\{ ad(p^{-1})g\left(V'(2\pi)g^{-1}p(2\pi)^{-1}g\right.\right.}\\
\shoveright{\left.\left. - ad(g^{-1})ad(p(2\pi))ad(g)a(W') + a(W')\right) g^{-1} \right\}g\gamma}\\
\shoveleft{\phantom{hV'_{(g^{-1}pg\gamma}} = \left(\frac{\theta}{2\pi}\right) g^{-1}p \left\{ ad(p^{-1})g\left(V'(2\pi)g^{-1}p(2\pi)^{-1}g\right.\right.}\\
\left.\left. - ad(g^{-1})ad(p(2\pi))ad(g)ad(g^{-1})a(W) + ad(g^{-1})a(W)\right) g^{-1} \right\}g\gamma
\end{multline*}
(for $W = W' g^{-1}$)
\begin{multline*}
\phantom{hV'_{(g^{-1}pg\gamma}} = \left(\frac{\theta}{2\pi}\right) g^{-1}p \left\{ ad(p^{-1})\left(gV'(2\pi)g^{-1}p(2\pi)^{-1} - ad(p(2\pi))a(W) + a(W)\right)  \right\}g\gamma\\
\shoveleft{\phantom{hV'_{(g^{-1}pg\gamma}} = \left(\frac{\theta}{2\pi}\right) g^{-1}p \left\{ ad(p^{-1})\left(V(2\pi)p(2\pi)^{-1} - ad(p(2\pi))a(W) + a(W)\right)  \right\}g\gamma}\\
\end{multline*}
(for $V = V' (\gamma, g)^{-1},$ so that $V'(2\pi) = g^{-1}V(2\pi)g)$).

So we see that the push forward of the vector $hV$ is horizontal (at $(g^{-1}pg\gamma, xg)$) and conversely the vector $hV'$ is the push forward of a horizontal vector at $(p,x).$ Thus we have defined a horizontal splitting of $T_{(p,x)}(P^{\ssv} G \times EG_n)$ for each $(p,x).$ To find the connection form for this connection we need to recover the Lie algebra element $(\xi, X)$ from the vector $(V,W).$ We know that the vector
\begin{multline*}
v(V,W)\\
= \left(V - \left(\frac{\theta}{2\pi}\right) p \left\{ad(p^{-1}) \left(V(2\pi)p(2\pi)^{-1} - ad(p(2\pi))a(W) + a(W)\right) \right\}, a(W)\right)
\end{multline*}
is the vertical component of $(V,W)$ and that the vertical vector generated by $(\xi, X) \in \Omega^{\ssv} \fg \rtimes \fg$ looks like
$$
( p (X - ad(p^{-1})X + \xi), \iota_x (X)).
$$
Thus to recover $\xi$ from $v(V,W)$ we just subtract $p(a(W) - ad(p^{-1}) a(W))$ and, writing $A$ for the part of the connection on $P^{\ssv} G,$ we have
\begin{multline*}
A(V,W) = \\
p^{-1}V - \left(\frac{\theta}{2\pi}\right)ad(p^{-1}) \left\{V(2\pi)p(2\pi)^{-1} - ad(p(2\pi))a(W) + a(W)\right\}\\
 - (a(W) - ad(p^{-1}) a(W)).
\end{multline*}
Therefore, the connection form $(A,a)$ is given by
$$
(A,a) = \left(\Theta - \left(\frac{\theta}{2\pi}\right) ad(p^{-1})\left\{ \ev_{2\pi}^*\hat{\Theta} - ad(p(2\pi))a + a \right\} - \left(a - ad(p^{-1})a\right), a \right)
$$
where $\Theta$ is the Maurer-Cartan form, $\hat{\Theta}$ is the right Maurer-Cartan form and $\ev_{2\pi} \colon P^{\ssv} G \to G$ is evaluation at the endpoint of a path. It can be easily checked that this form satisfies the conditions for a connection. It will be useful later on to write this as a form valued in $L^{\ssv} \fg.$ To do this we use the isomorphism of Lie algebras given in section \ref{SS:LG classification}. The connection form becomes
$$
A_{L^{\ssv} \fg} = \Theta - \left(\frac{\theta}{2\pi}\right) ad(p^{-1})\left\{ \ev_{2\pi}^*\hat{\Theta} - ad(p(2\pi))a + a \right\}  + ad(p^{-1})a.
$$

To calculate the string class we will need the curvature of this connection and a Higgs field. As usual, the curvature (as an $L^{\ssv} \fg$-valued form) is given by the formula
$$
F_{L^{\ssv} \fg} = DA_{L^{\ssv} \fg}
$$
where $D$ is the covariant exterior derivative. So we have
$$
F_{L^{\ssv} \fg} ((V,W), (V',W')) = \tfrac{1}{2} A_{L^{\ssv} \fg} ([h(V,W), h(V',W')]).
$$
Now,
\begin{multline*}
[h(V,W), h(V',W')] = ([hV, hV'], [hW,hW'])\\
= \left(\left[\left(\frac{\theta}{2\pi}\right) p \left\{ad(p^{-1}) \left(V(2\pi)p(2\pi)^{-1} - ad(p(2\pi))a(W) + a(W)\right) \right\},\right.\right.\\
\left(\left.\frac{\theta}{2\pi}\right) p \left\{ad(p^{-1}) \left(V'(2\pi)p(2\pi)^{-1} - ad(p(2\pi))a(W') + a(W')\right) \right\}\right],\\
\left.\vphantom{\frac{\theta}{2\pi}}[hW , hW']\right)
\end{multline*}
and calculating just the first slot gives
\begin{multline*}
p\left(\frac{\theta}{2\pi}\right)^2 ad(p^{-1}) \left\{ [V(2\pi)p(2\pi)^{-1}, V'(2\pi)p(2\pi)^{-1}] \right.\\
- [V(2\pi)p(2\pi)^{-1},ad(p(2\pi))a(W')] + [V(2\pi)p(2\pi)^{-1},a(W')]\\
-[ad(p(2\pi))a(W),V'(2\pi)p(2\pi)^{-1}] + ad(p(2\pi))[a(W),a(W')]\\
 - [ad(p(2\pi))a(W), a(W')] + [a(W),V'(2\pi)p(2\pi)^{-1}]\\ 
 \left.- [a(W),ad(p(2\pi))a(W')]  + [a(W),a(W')]\right\}.
\end{multline*}
This yields
\begin{multline*}
F_{L^{\ssv} \fg} =\\
 \left(\left(\frac{\theta}{2\pi}\right)^2 - \left(\frac{\theta}{2\pi}\right)\right) ad(p^{-1}) \left\{ \tfrac{1}{2}[\ev_{2\pi}^*\hat{\Theta}, \ev_{2\pi}^*\hat{\Theta}] - [\ev_{2\pi}^*\hat{\Theta}, ad(p(2\pi)^{-1})a] + \tfrac{1}{2} [a,a]\right.\\
 \left.+ [\ev_{2\pi}^*\hat{\Theta}, a] - [ad(p(2\pi))a, a] + [a,a] \right\}  - \left(\frac{\theta}{2\pi}\right) ad(p^{-1})(f - ad(p(2\pi))f) + ad(p^{-1})f
\end{multline*}
where $f$ is the curvature of $a.$

The other piece of data we need to calculate the string class is a Higgs field for $EL^{\ssv} G.$ Define the map $\Phi \colon P^{\ssv} G \times EG_n \to \Omega^{\ssv} \fg \rtimes \fg$ by
$$
\Phi(p,x) = (p^{-1} \partial p, 0).
$$
Or, as a map to $L^{\ssv} \fg,$
$$
\Phi_{L^{\ssv} \fg} (p,x) = p^{-1} \partial p.
$$
Then by the calculation at the end of section \ref{SS:LG classification} we see that $\Phi_{L^{\ssv} \fg}$ is a Higgs field for $P^{\ssv} G \times EG_n.$ Next we need to calculate
$$
\nabla \Phi_{L^{\ssv} \fg} = d \Phi_{L^{\ssv} \fg} + [A_{L^{\ssv} \fg}, \Phi_{L^{\ssv} \fg}] - \partial A_{L^{\ssv} \fg}.
$$
We have
\begin{align*}
d\Phi_{L^{\ssv} \fg} (V,W) 	&= \frac{d}{dt}\bigg|_0 \Phi_{L^{\ssv} \fg}(p e^{t\xi})\\
				&= \frac{d}{dt}\bigg|_0 (e^{-t\xi} p^{-1} \partial (p e^{t\xi}))\\
				&= \frac{d}{dt}\bigg|_0 (e^{-t\xi} p^{-1} \partial p e^{t\xi} + e^{-t\xi}  \partial e^{t\xi}) \\
				&= p^{-1} \partial p \xi - \xi p^{-1} \partial p + \partial \xi,
\end{align*}
for $V = \frac{d}{dt}\big|_0\, p \exp(t\xi).$ That is,
$$
d \Phi_{L^{\ssv} \fg} = [\Phi_{L^{\ssv} \fg}, \Theta] + \partial \Theta.
$$
So
\begin{multline*}
\nabla \Phi_{L^{\ssv} \fg} = [\Phi_{L^{\ssv} \fg}, \Theta] + \partial \Theta\\ 
+ \left[ \Theta - \left(\frac{\theta}{2\pi}\right) ad(p^{-1})\left\{ \ev_{2\pi}^*\hat{\Theta} - ad(p(2\pi))a + a \right\}  + ad(p^{-1})a, \Phi_{L^{\ssv} \fg}\right]\\
 - \partial \left( \Theta - \left(\frac{\theta}{2\pi}\right) ad(p^{-1})\left\{ \ev_{2\pi}^*\hat{\Theta} - ad(p(2\pi))a + a \right\}  + ad(p^{-1})a \right)
\end{multline*}
\vspace{-3ex}
\begin{multline*}
\phantom{\nabla \Phi_{L^{\ssv} \fg}} = \frac{1}{2\pi} ad(p^{-1}) \left\{ \ev_{2\pi}^*\hat{\Theta} - ad(p(2\pi))a + a\right\}.\\
\end{multline*}
So the string class for $P^{\ssv} G \times EG_n$ is
\begin{multline*}
-\frac{1}{4\pi^2}\int_{S^1}\left\< \left( \frac{\theta^2}{4\pi^2} - \frac{\theta}{2\pi}\right) \left( \tfrac{1}{2}[\ev_{2\pi}^*\hat{\Theta}, \ev_{2\pi}^*\hat{\Theta}] - [\ev_{2\pi}^*\hat{\Theta}, ad(p(2\pi)^{-1})a] \right.\right.\\
+ \tfrac{1}{2} [a,a]
 \left.+ [\ev_{2\pi}^*\hat{\Theta}, a] - [ad(p(2\pi))a, a] + [a,a] \right)\\
- \left(\frac{\theta}{2\pi}\right)(f - ad(p(2\pi))f) + f, \left. \frac{1}{2\pi} \left( \ev_{2\pi}^*\hat{\Theta} - ad(p(2\pi))a + a\right)\right\>
\end{multline*}
\begin{multline*}
= -\frac{1}{8\pi^2}\left\< -\tfrac{1}{3}\left(\tfrac{1}{2}[\hat{\Theta}, \hat{\Theta}] - [\hat{\Theta}, ad(p(2\pi)^{-1})a] \right.\right.\\
+ \tfrac{3}{2} [a,a]
 \left.+ [\hat{\Theta}, a] - [ad(p(2\pi))a, a]  \right)\\
+ ad(p(2\pi))f + f, \left.\left(\hat{\Theta} - ad(p(2\pi))a + a\right)\right\>.
\end{multline*}

\chapter{String structures for $LG\rtimes S^1$-bundles}\label{C:LGxS^1}


Thus far we have discussed central extensions of both the loop group (in chapter \ref{C:background}) and the based loop group (in chapter \ref{C:characteristic classes}). The loop group $LG$ has a natural action of the circle given by rotating loops. In this chapter, we shall consider the more general case where we allow rotations of the loops in $LG.$ That is, we shall be working with the semi-direct product $LG \rtimes S^1.$ This group arises when we consider a natural generalisation of the caloron correspondence from section \ref{S:Higgs fields, LG-bundles,...}. There we showed that a $G$-bundle over $M \times S^1$ corresponds to an $LG$-bundle over $M$. If we allow the base space of the $G$-bundle to be a non-trivial $S^1$-bundle (rather than $M \times S^1$) we obtain not an $LG$-bundle but an $\LGS$-bundle. If, further, we consider a non-trivial $S^1$ fibre bundle (instead of a principal bundle), we obtain an $LG \rtimes \Diff(S^1)$-bundle.

In this chapter then, we will calculate the obstruction to lifting a principal $LG \rtimes S^1$-bundle $P$ to a principal $\widehat{LG\rtimes S^1}$-bundle $\widehat{P}.$ In section \ref{S:Higgs fields, LGxS^1-bundles,...} we will construct a correspondence for $\LGS$-bundles in analogy with the caloron correspondence from chapter \ref{C:background}. This will be used to prove a theorem which extends Theorem \ref{T:LGPont} relating the string class and the first Pontrjagyn class. In section \ref{S:Diff(S^1)} we shall consider the lifting problem for the more general case where we allow general (orientation preserving) diffeomorphisms of the loops in $LG,$ that is, principal bundles with structure group $LG\rtimes \Diff(S^1).$

\section{The string class of an $LG\rtimes S^1$-bundle}\label{S:LGxS^1 string class}

In this section we present a formula for the obstruction to lifting a principal $\LGS$-bundle $P$ to a principal $\widehat{\LGS}$-bundle $\widehat{P},$ which we call the \emph{string class} of $P.$ We shall follow the methods of \cite{Murray:2003}, outlined in section \ref{S:LG string class}. In section \ref{SS:reduced splittings} we will give another method for calculating the 3-curvature of a lifting bundle gerbe, first presented in \cite{Gomi:2003}, and apply this to the problem of the string class of an $\LGS$-bundle.

\subsection{The string class via lifting bundle gerbes}\label{SS:LGxS^1 string class}

Let $\LGS$ be the semi-direct product, whose multiplication is given by
$$
(\gamma_1, \phi_1) (\gamma_2, \phi_2) = (\gamma_1 \rho_{\phi_1}(\gamma_2), \phi_1 + \phi_2),
$$
where $\rho_{\phi}(\gamma)(\theta) = \gamma(\theta - \phi).$ For convenience, let us record some facts about the Lie algebra of $\LGS$ here. The bracket on the Lie algebra $\Lgs$ is given by
$$
[ (\xi, x), (\zeta, y) ] = ( [\xi, \zeta] - x \partial \zeta + y \partial \xi ,  0)
$$
and the adjoint action of $\LGS$ on $\Lgs$ is
$$
ad(\gamma, \phi)(\xi, x) = \left( ad(\gamma) \rho_{\phi}(\xi) + x \, \partial \gamma \gamma^{-1}, x \right).
$$

\subsubsection{The central extension of $\LGS$}

Recall from section \ref{S:LG string class} that in order to perform calculations involving the lifting bundle gerbe, we needed an explicit construction of the central extension of $LG.$ This was given following the construction in section \ref{S:central extensions} in terms of a pair of differential forms satisfying a certain compatibility condition. Namely, a pair $(R, \a),$ where $R$ is a closed, integral 2-form on $LG$ and $\a$ is a 1-form on $LG \times LG,$ satisfying the conditions $\d R = d\a$ and $\d \a = 0.$ In a similar manner, for what follows we will require an explicit construction of the central extension of $\LGS.$ Note, however, that the construction in section \ref{S:central extensions} only works for $\cG$ a simply connected Lie group. This is because in order to construct the extension given the pair $(R, \a)$ we used the fact that a flat bundle over a simply connected base has a section satisfying certain conditions. This allowed us to find a $U(1)$-bundle $P$ over $\cG$ such that $\d P \to \cG \times \cG$ was trivial and had a section which defined the multiplication on the central extension.\footnote{See the discussion in section \ref{SS:slb's and ce's}.} However, even though the semi-direct product $\LGS$ is not simply connected we can modify the construction from section \ref{S:central extensions} slightly to cover this case \cite{Murray:2003}. This involves replacing the 2-form $R$ with a differential character \cite{Cheeger:1985} for the bundle $\widehat{\cG} \to \cG$. That is, we add to our pair $(R, \a)$ a homomorphism $h\colon Z_1(\cG) \to U(1)$ satisfying
$$
h(\partial \sigma) = \exp \left( \int_\sigma R \right)
$$
for every two-cycle $\sigma$ in $\cG.$ We also require the compatibility condition
$$
(\d h) (\gamma) = \exp \left( \int_\gamma \a \right)
$$
for every closed one-cycle $\gamma$ in $\cG \times \cG.$

Therefore, we need to find a triple of objects $(R, \a, h)$ as above. Note first that
$$
H^2(\LGS) \simeq H^2(LG).
$$
To see this, we observe that as $\LGS = LG \times S^1$ as a space, the K\"unneth formula (see \cite{Bott-Tu}) gives
$$
H^2(\LGS) \simeq H^2(LG) \otimes H^0(S^1) \oplus H^1(LG) \otimes H^1(S^1),
$$
since $H^2(S^1) = 0.$ Now, $H^0(S^1) \simeq H^1(S^1) \simeq \RR,$ so we have
$$
H^2(\LGS) \simeq H^2(LG) \oplus H^1(LG).
$$
Recall, however, that $LG \simeq \Omega G \times G$ as a space, and so $\pi_1(LG) = \pi_2(G) \times \pi_1(G).$ Therefore, as $G$ is simply connected, so is $LG,$ and thus $H_1(LG, \RR) = 0$ by the Hurewicz Theorem (see for example \cite{Hatcher}). Therefore, by the Universal Coefficient Theorem (see for example \cite{Massey:1991}) $H^1(LG) =0,$ and so $H^2(\LGS) \simeq H^2(LG).$ Thus, we take as the 2-form $R,$ the pull-back of the form from section \ref{S:LG string class} to $\LGS.$ That is,
$$
R = \frac{i}{4\pi} \int_{S^1} \< \Theta, \partial \Theta \>\, d\theta.
$$
Note that since we are integrating over the circle, this expression is unchanged when each term is rotated by a fixed angle. That is,
$$
\frac{i}{4\pi} \int_{S^1} \< \rho_\phi (\Theta), \partial \rho_\phi(\Theta) \>\, d\theta = \frac{i}{4\pi} \int_{S^1} \< \Theta, \partial \Theta \>\, d\theta
$$
Now, to find $\a$ we need to calculate $\d R = \pi_1^*R - m^*R + \pi_2^*R,$ where as before, $\pi_i$ is the projection $\LGS \times \LGS \to \LGS$ which omits the $i^\text{th}$ factor and $m$ is the multiplication defined above. As in chapter \ref{C:background}, $\pi_i^*R$ is given by
$$
\frac{i}{4\pi} \int_{S^1} \< \pi_i^*\Theta, \partial \pi_i^*\Theta \>\, d\theta
$$
and so it remains to calculate $m^*R.$ For this, note that a tangent vector to $\LGS$ at the point $(\gamma, \phi)$ can be written as $(\gamma, \phi)(\xi, x) = (\gamma \rho_\phi (\xi), x_\phi)$ for some $(\xi, x) \in \Lgs$ by using the left multiplication to transport elements of the Lie algebra to the point $(\gamma, \phi).$ Therefore, we can calculate $m^*R$ by noting that
$$
m^*R((\gamma_1 \rho_{\phi_1} (\xi_1), x_{1\phi_1}), (\gamma_2 \rho_{\phi_2} (\xi_2), x_{2\phi_2})) = R(m_*((\gamma_1 \rho_{\phi_1} (\xi_1), x_{1\phi_1}), (\gamma_2 \rho_{\phi_2} (\xi_2), x_{2\phi_2})))
$$
and calculating the push-forward of $m.$ We have
\begin{multline*}
	m_*((\gamma_1 \rho_{\phi_1} (\xi_1), x_{1\phi_1}), (\gamma_2 \rho_{\phi_2} (\xi_2), x_{2\phi_2}))\\
	=\frac{d}{dt}\bigg|_0 \left( \gamma_1(1+t\xi_1^{\rho_1})\rho_{(\phi_1+tx_1)}(\gamma_2)\rho_{(\phi_1+tx_1)}(1+t\xi_2^{\rho_2})), \phi_1+\phi_2+t(x_1+x_2)\right),
\end{multline*}
where we have written (for example) $\xi_1^{\rho_1}$ for $\rho_{\phi_1}(\xi_1)$. As the multiplication on the $S^1$ factor is not twisted, the second slot above will give $x_1 + x_2.$ Thus it suffices to calculate the first slot only. Using the fact that
$$
\frac{d}{dt}\bigg|_0 \rho_{(\phi_1 + tx_1)} (\gamma_2) = -x_1\rho_{\phi_1}(\partial \gamma_2),
$$
we have
\begin{align*}
	\frac{d}{dt}\bigg|_0 &\left( \gamma_1(1+t\xi_1^{\rho_1})\rho_{(\phi_1+tx_1)}(\gamma_2)\rho_{(\phi_1+tx_1)}(1+t\xi_2^{\rho_2})\right)\\
	&= \gamma_1 \xi_1^{\rho_1} \gamma_2^{\rho_1} + \gamma_1 \gamma_2^{\rho_1} \xi_2^{\rho_2} -x_1 \gamma_1\partial \gamma_2^{\rho_1}\\
	&= \gamma_1 \gamma_2^{\rho_1} \rho_{(\phi_1+\phi_2)} \left((\gamma_2^{-1}\xi_1\gamma_2^{\vphantom{-1}})^{\rho_2^{-1}} + \xi_2 - x_1(\gamma_2^{-1}\partial \gamma_2^{\vphantom{-1}})^{\rho_2^{-1}}\right).
\end{align*}
Therefore, $m^*R$ evaluated on the pairs of tangent vectors $\left((\gamma_1, \phi_1)(\xi_1, x_1), (\gamma_2, \phi_2)(\xi_2, x_2)\right)$ and $\left((\gamma_1, \phi_1)(\zeta_1, y_1), (\gamma_2, \phi_2)(\zeta_2, y_2)\right)$ is given by
$$
\frac{i}{4\pi} \int_{S^1} \left\< (ad(\gamma_2^{-1})\xi_1)^{\rho_2^{-1}} + \xi_2 - x_1(\gamma_2^{-1}\partial \gamma_2^{\vphantom{-1}})^{\rho_2^{-1}}, \partial \left( (ad(\gamma_2^{-1})\zeta_1)^{\rho_2^{-1}} + \zeta_2 - y_1(\gamma_2^{-1}\partial \gamma_2^{\vphantom{-1}})^{\rho_2^{-1}}\right) \right\> \, d\theta,
$$
where we have used the fact that the integral is unchanged by rotating everything by $\rho_{(\phi_1 + \phi_2)}^{-1}.$ Expanding this, we have
\begin{multline*}
\frac{i}{4\pi} \int_{S^1} \left\< (ad(\gamma_2^{-1})\xi_1), \partial (ad(\gamma_2^{-1})\zeta_1) \right\>
+ \left\< \xi_2, \partial \zeta_2 \right\>\\
+ x_1 y_1 \left\< (ad(\gamma_2^{-1}) Z_2), \partial (ad(\gamma_2^{-1}) Z_2) \right\>\\
+ \left\<(ad(\gamma_2^{-1})\xi_1), \partial \zeta_2^{\rho_2} \right\>
+ \left\< \xi_2^{\rho_2}, \partial (ad(\gamma_2^{-1}) \zeta_1) \right\>\\
- y_1 \left\< \xi_2^{\rho_2}, \partial (ad(\gamma_2^{-1})Z_2) \right\>
- x_1 \left\< (ad(\gamma_2^{-1})Z_2), \partial \zeta_2^{\rho_2} \right\>\\
- y_1 \left\< (ad(\gamma_2^{-1}) \xi_1), \partial (ad(\gamma_2^{-1}) Z_2) \right\>\\
- x_1 \left\< (ad(\gamma_2^{-1}) Z_2), \partial (ad(\gamma_2^{-1}) \zeta_2) \right\>  d \theta,
\end{multline*}
where as before $Z$ is the function $\gamma \mapsto \partial \gamma \gamma^{-1}$ and, again, we have used the rotation invariance of the integral. Using the $ad$-invariance of the Killing form and integration by parts, along with the identity from section \ref{S:LG string class},
$$
\partial \left(ad(\gamma^{-1}) X \right) = ad(\gamma^{-1})[X, Z] + ad(\gamma^{-1}) \partial X
$$
for a vector $X \in L\fg,$ this simplifies to
\begin{multline*}
\frac{i}{4\pi} \int_{S^1} \left\< [\xi_1, \zeta_1], Z_2 \right\>
+ \left\< \xi_1, \partial \zeta_1 \right\>
+ \left\< \xi_2, \partial \zeta_2 \right\>\\
+ \left\< ad(\gamma_2^{-1}) \xi_1, \partial \zeta_2^{\rho_2} \right\>
- \left\< \partial \xi_2^{\rho_2}, ad(\gamma_2^{-1}) \zeta_1 \right\>
- x_1 \left< Z_2, \partial \zeta_1 \right\>
 + y_1 \left\< \partial \xi_1, Z_2 \right\>\\
- x_1 \left\<  ad(\gamma_2^{-1}) Z_2, \partial \zeta_2^{\rho_2} \right\>
+ y_1 \left\< \partial \xi_2^{\rho_2}, ad(\gamma_2^{-1}) Z_2 \right\> d\theta,
\end{multline*}
or simply
\begin{multline*}
m^*R = \frac{i}{4\pi} \int_{S^1} \left\< [\Theta_1, \Theta_1], Z_2 \right\>
+ \left\< \Theta_1, \partial \Theta_1 \right\>
+ \left\< \Theta_2, \partial \Theta_2 \right\>\\
+ 2 \left\< ad(\gamma_2^{-1}) \Theta_1, \partial \Theta_2^{\rho_2} \right\>
- 2\left\< \mu_1 ad(\gamma_2^{-1}) Z_2, \partial \Theta_2^{\rho_2} \right\>
- 2\left\< \mu_1 Z_2, \partial \Theta_1 \right\> d\theta,
\end{multline*}
where $\mu$ represents the Maurer-Cartan form on $S^1$. Therefore, we have
\begin{multline*}
\d R = \frac{i}{2\pi} \int_{S^1} -\tfrac{1}{2} \left\< [\Theta_1, \Theta_1], Z_2 \right\>
- \left\< ad(\gamma_2^{-1}) \Theta_1, \partial \Theta_2^{\rho_2} \right\>\\
+ \left\< \mu_1 ad(\gamma_2^{-1}) Z_2, \partial \Theta_2^{\rho_2} \right\>
+ \left\< \mu_1 Z_2, \partial \Theta_1 \right\> d\theta.
\end{multline*}
Recall from section \ref{S:LG string class} that for the loop group case, the form $\a$ such that $d\a = \d R$ is given by
$$
\a = \frac{i}{2\pi} \int_{S^1} \left\< \Theta_1, Z_2 \right\> d\theta.
$$
When evaluated on the vector $(\gamma_1 \xi_1, \gamma_2 \xi_2)$ tangent to the point $(\gamma_1, \gamma_2) \in LG \times LG,$ $\a$ is given by
$$
\frac{i}{2\pi} \int_{S^1} \left\< \xi_1, \partial \gamma_2^{\vphantom{-1}}  \gamma_2^{-1}\right\> d\theta.
$$
Consider the generalisation of this form to $\LGS \times \LGS.$ That is, define $\a_1$ as
$$
\a_1 (\gamma_1 \xi^{\rho_1}, x_{1\phi_1}, \gamma_2 \xi^{\rho_2}, x_{2\phi_2}) = \frac{i}{2\pi} \int_{S^1} \left\< \xi_1, \partial \gamma_2^{\vphantom{-1}} \gamma_2^{-1}\right\>,
$$
or
$$
\a_1 = \frac{i}{2\pi} \int_{S^1} \big\< \Theta_1^{\rho_1^{-1}}, Z_2 \big\> \, d\theta.
$$
We can calculate the derivative of this form via
$$
d\a_1(X, Y) = \tfrac{1}{2} \left\{ X(\a_1 (Y)) - Y(\a_1(X)) - \a_1 ([X, Y]) \right\},
$$
for tangent vectors $X$ and $Y.$ Thus we need to calculate
\begin{align*}
(\gamma_1\xi_1^{\rho_1}, & x_{1\phi_1},\gamma_2 \xi_2^{\rho_2}, x_{2\phi_2})\left(\a_1(\gamma_1\zeta_1^{\rho_1}, y_{1\phi_1}, \gamma_2 \zeta_2^{\rho_2}, y_{2\phi_2})\right)\\
		& = \frac{d}{dt}\bigg|_0 \frac{i}{2\pi} \int_{S^1} \left\< \zeta_1, \partial (\gamma_2(1 + t \xi_2^{\rho_2})) (1 - t \xi_2^{\rho_2}) \gamma_2^{-1} \right\> d\theta\\
		& =  \frac{i}{2\pi} \int_{S^1} \left\< \zeta_1, ad(\gamma_2) \partial \xi_2^{\rho_2}  \right\> d\theta,
\end{align*}
and
\begin{align*}
\a_1([(\gamma_1\xi_1^{\rho_1}, x_{1\phi_1}), & (\gamma_1\zeta_1^{\rho_1}, y_{1\phi_1})], [(\gamma_2 \xi_2^{\rho_2}, x_{2\phi_2}), (\gamma_2 \zeta_2^{\rho_2}, y_{2\phi_2})])\\
		& = \frac{i}{2\pi} \int_{S^1} \left\< [ (\xi_1, x_1), (\zeta_1, y_1) ], \partial \gamma_2^{\vphantom{-1}} \gamma_2^{-1} \right\> d\theta\\
		& = \frac{i}{2\pi} \int_{S^1} \left\< [ \xi_1, \zeta_1 ], \partial \gamma_2^{\vphantom{-1}} \gamma_2^{-1} \right\> - \left\< x_1 \partial \zeta_1 - y_1 \partial \xi_1 , \partial \gamma_2^{\vphantom{-1}} \gamma_2^{-1} \right\> d\theta.
\end{align*}
Therefore, we have
$$
d \a_1 = \frac{1}{2\pi} \int_{S^1} -\tfrac{1}{2} \left\< [ \Theta, \Theta], Z_2^{\vphantom{-1}} \right\>
- \left\< ad(\gamma_2^{-1}) \Theta_1, \partial \Theta_2^{\rho_2} \right\>
+ \left\< \mu_1 Z_2, \partial \Theta_1^{\vphantom{-1}} \right\> d\theta.
$$
Note that $\d R$ does not equal $d\a_1.$ However,
$$
\d R - d \a_1 = \frac{i}{2\pi} \int_{S^1} \left\< \mu_1 ad(\gamma_2^{-1}) Z_2, \partial \Theta_2^{\rho_2} \right\> d\theta.
$$
Using the identity
$$
ad(\gamma) \partial \Theta^\rho = d Z,
$$
we see that
$$
\d R - d \a_1 = \frac{i}{2\pi} \int_{S^1} \left\< \mu_1  Z_2, dZ_2 \right\> d\theta.
$$
Now, if we define
$$
\a_2 = -\frac{i}{4\pi} \int_{S^1} \left\< \mu_1 Z_2, Z_2 \right\> d\theta,
$$
then
\begin{align*}
d\a_2	&= \frac{i}{4\pi} \int_{S^1} \left\< \mu_1 d Z_2, Z_2 \right\> + \left\< \mu_1 Z_2, d Z_2 \right\> d\theta\\
		&= \frac{i}{2\pi} \int_{S^1} \left\< \mu_1 Z_2, d Z_2 \right\> d\theta\\
		&= \d R - d\a_1.
\end{align*}
Thus $\a$ is given by
$$
\a = \frac{i}{2\pi} \int_{S^1} \big\< \pi_2^*\Theta^{\rho^{-1}} \! - \tfrac{1}{2} \pi_2^*\mu \, \pi_1^*Z, \pi_1^*Z \big\> \, d\theta,
$$
and $\d R = d\a.$ One can also easily check that $\d \a = 0.$ Notice that the $2$-form $R$ is left invariant and the $1$-form $\a$ is left invariant in the first slot. To find the homomorphism $h\colon Z_1(\LGS) \to U(1)$ we note that since $\pi_1(\LGS) = \ZZ$ any cycle $a \in Z_1(\LGS)$ can be written as $n\gamma + \partial \sigma,$ for some two-cycle $\sigma,$ where $\gamma$ is the generator of $H_1(\LGS),$ a loop around the $S^1$ factor. It is easy to see that the integral of $\a$ over the generators of $H_1(\LGS \times \LGS)$ vanishes, that is,
$$
\int_{\gamma_1} \a = 0 = \int_{\gamma_2} \a 
$$
for $\gamma_1, \gamma_2$ loops around the first and second $S^1$ factors respectively. This suggests that we define
$$
h(a) = h(\partial \sigma) = \exp \left( \int_\sigma R \right).
$$
This is well defined since if $a = n\gamma + \partial \sigma = n \gamma + \partial \sigma'$ then $\partial (\sigma - \sigma') = 0$ and so $\ds\int_{\sigma - \sigma'} R \in 2\pi i \ZZ$ (since $R$ is integral). Because the integral of $\a$ over the generators of $H_1(\LGS \times \LGS)$ vanishes, it is easy to check that for any one-cycle $\gamma$ we have
$$
(\d h) (\gamma) = \exp \left( \int_\gamma \a \right).
$$
Thus we have proven
\begin{proposition}
The triple $(R, \a, h)$ as above determines a central extension of the semi-direct product $\LGS.$
\end{proposition}

\subsubsection{A connection for the lifting bundle gerbe}

Now that we have a construction of the central extension of $\LGS,$ the next step is to write down a bundle gerbe connection for the lifting bundle gerbe. Recall from section \ref{S:LG string class} that if $P$ is an $\LGS$-bundle and $\nu$ is a connection on the central extension $\widehat{\LGS}$ thought of as a bundle over $\LGS$ then a bundle gerbe connection is given by $\tau^*\nu - \epsilon,$ where $\epsilon$ is some $1$-form on $P^{[2]}$ satisfying $\d \epsilon = \tau^*\a.$ In the $LG$ case, this form was given by
$$
\epsilon = \frac{i}{2\pi} \int_{S^1} \< \pi_2^* A , \tau^*Z\>\, d\theta,
$$
where $A$ is a connection on $P.$ As mentioned in section \ref{S:LG string class}, it is possible to write $\epsilon$ in general in terms of $\a$ \cite{Stevenson:comm}. We shall now demonstrate how to do this. Let $P$ be a $\cG$-bundle with connection $A.$ Recall that $A$ satisfies
$$
\pi_1^*A = ad(\tau_{12}^{-1}) \pi_2^* A + \tau_{12}^*\Theta.
$$
For tangent vectors $(X_1, X_2, X_3)$ at $(p_1, p_2, p_3) \in P^{[3]},$ we can calculate
\begin{multline*}
(\d\a)_{(1, \tau_{12}, \tau_{23})}( A(X_1), \tau_{12}(X_1, X_2), \tau_{23}(X_2,X_3)) =\\
	\shoveleft{\phantom{(\d\a)}\a_{(\tau_{12},\tau_{23})}( \tau_{12}(X_1, X_2), \tau_{23}(X_2,X_3))}\\
	 - \a_{(\tau_{12},\tau_{23})}( m_*(A(X_1),\tau_{12}(X_1, X_2)), \tau_{23}(X_2,X_3))\\
	+ \a_{(1,\tau_{12}\tau_{23})}( A(X_1),m_*(\tau_{12}(X_1, X_2), \tau_{23}(X_2,X_3)))\\
	 - \a_{(1,\tau_{12})}(A(X_1),\tau_{12}(X_1, X_2)).
\end{multline*}
Notice that the first term above is actually $\tau^* \a.$ Since $\d\a = 0,$ we have
\begin{multline*}
(\tau^*\a)_{(p_1,p_2,p_3)} (X_1,X_2,X_3) =\\
	\a_{(\tau_{12},\tau_{23})}(m_*(A(X_1),\tau_{12}(X_1, X_2)), \tau_{23}(X_2,X_3))\\
	- \a_{(1,\tau_{12}\tau_{23})}(A(X_1),m_*(\tau_{12}(X_1, X_2), \tau_{23}(X_2,X_3)))\\
	+ \a_{(1,\tau_{12})}(A(X_1),\tau_{12}(X_1, X_2)).
\end{multline*}
Now, if we define $\epsilon$ in terms of $\a$ and $A$ as
$$
\epsilon_{(p_1,p_2)}(X_1, X_2) =\a_{(1, \tau_{12})}(A(X_1), \tau_{12}(X_1,X_2))
$$
then we have
\begin{multline*}
(\d\epsilon)_{(p_1, p_2,p_3)}(X_1,X_2,X_3) =\\
	\a_{(1, \tau_{23})}(A(X_2), \tau_{23}(X_2, X_3)) - \a_{(1, \tau_{13})}(A(X_1), \tau_{13}(X_1,X_3))\\
	+\a_{(1,\tau_{12})}(A(X_1), \tau_{12}(X_1,X_2)).
\end{multline*}
Using the fact that $\tau_{13} = \tau_{12}\tau_{23},$ we see
$$
\a_{(1, \tau_{13})}( A(X_1), \tau_{13}(X_1,X_2))=
\a_{(1,\tau_{12}\tau_{23})}( A(X_1),m_*(\tau_{12}(X_1, X_2), \tau_{23}(X_2,X_3)))
$$
and since $\a$ is left invariant in the first slot, and using the equation above relating $A(X_1)$ and $A(X_2),$ we have
\begin{align*}
\a_{(1, \tau_{23})}( & A(X_2), \tau_{23}(X_2, X_3))\\
			&= \a_{(\tau_{12}, \tau_{23})}( \tau_{12}A(X_2), \tau_{23}(X_2, X_3))\\
			&= \a_{(\tau_{12}, \tau_{23})}( \tau_{12}(ad(\tau_{12}^{-1})A(X_1) + \tau_{12}^{-1}(\tau_{12}(X_1,X_2))), \tau_{23}(X_2, X_3))\\
			&= \a_{(\tau_{12}, \tau_{23})}( \tau_{12}ad(\tau_{12}^{-1})A(X_1) + \tau_{12}(X_1,X_2), \tau_{23}(X_2, X_3)),
\end{align*}
which equals
$$
\a_{(\tau_{12},\tau_{23})}( m_*(A(X_1),\tau_{12}(X_1, X_2)), \tau_{23}(X_2,X_3)).
$$
Thus we have $\d \epsilon = \tau^*\a.$

Consider now the $\LGS$-bundle $P.$ Choose a connection $(A, a)$ for $P,$ where $A$ and $a$ are $1$-forms on $P$ with values in $L\fg$ and $i\RR$ respectively. Note that $a$ is a connection for the associated $S^1$-bundle $P/LG$ whereas $A$ is not a connection form. In fact, if $X$ is a tangent vector to $P,$ we have
\begin{align*}
(A, a)(X(\gamma, \phi)) &= ad(\gamma, \phi)^{-1}(A(X), a(X))\\
				&= \left( \rho_{\phi^{-1}}\left( ad(\gamma^{-1}) A(X) - a(X) \gamma^{-1} \partial \gamma \right) , a(X) \right),
\end{align*}
and so $A$ does not have the correct transformation properties to be a connection.\footnote{Notice the similarity with the treatment of connections for the universal $\Omega G \rtimes G$-bundle in section \ref{S:universal string class}.} Given $(A, a)$ then, we can write down the $1$-form $\epsilon \in \Omega^1(P^{[2]})$ as above:
$$
\epsilon = \frac{i}{2\pi} \int_{S^1} \left\< \pi_2^* A - \tfrac{1}{2} \pi_2^*a \, \tau^*Z, \tau^*Z \right\> d\theta.
$$
It is easy to check that $\d \epsilon = \tau^* \a$ and so we have that $\tau^*\nu - \epsilon$ is a connection for the lifting bundle gerbe. Of course, as in section \ref{S:LG string class} we are concerned with finding a curving for this bundle gerbe and so we are really interested in calculating the curvature of this connection, given by $\tau^*R - d\epsilon.$ Recall that for a connection $A$ on a $\cG$-bundle, we have the formula
$$
\pi_1^*A = ad(\tau^{-1}) \pi_2^* A + \tau^*\Theta.
$$
In the case where $\cG = \LGS,$ the formula relating $\pi_1^*(A, a) = (A_2, a_2)$ and $\pi_2^*(A, a) = (A_1, a_1)$ is
$$
(A_2,a_2) =\left( \rho_{\tau_{S^1}}^{-1}\left(ad(\tau_{LG}^{-1})A_1 - a_1 \tau_{LG}^{-1} \partial\tau_{LG}^{\vphantom{-1}} \right) + \tau_{LG}^* (\rho_{\tau_{S^1}}^{-1}(\Theta) ), a_1 + \tau_{S^1}^*\mu \right)
$$
where we have written the difference map $\tau$ as $(\tau_{LG}, \tau_{S^1}).$ That is, $\tau_{LG}$ is the $LG$ part of $\tau$ and $\tau_{S^1}$ is the circle part. From now on, we will simply write $\tau$ and assume that it is clear from the context which part we mean. In particular, then, we have
$$
\tau^* \rho_{\tau}^{-1}(\Theta) =A_2 - \rho_{\tau}^{-1}\left(ad(\tau^{-1})A_1 + a_1\tau^{-1} \partial\tau \right).
$$
Note that here we have used the fact that the Maurer-Cartan form on $\LGS$ is not the pair $(\Theta, \mu)$ but in fact includes a rotation of $\Theta.$ So at the point $(\gamma, \phi),$ it is given by $(\rho_{\phi^{-1}}(\Theta), \mu).$

We can use this to calculate $\tau^*R - d\epsilon.$ Writing $A^\rho$ for $\rho(A)$ and so on, as before, we have
\begin{multline*}
\tau^*R = \frac{i}{4\pi} \int_{S^1} \<\tau^*\Theta, \partial \tau^* \Theta \> d\theta\\
\shoveleft{\phantom{\tau^*R} = \frac{i}{4\pi} \int_{S^1} \<A_2^\rho - ad(\tau^{-1})A_1 + a_1 \tau^{-1} \partial\tau, \partial (A_2^\rho - ad(\tau^{-1})A_1 + a_1 \tau^{-1} \partial\tau)\> d\theta}\\
\shoveleft{\phantom{\tau^*R} = \frac{i}{4\pi} \int_{S^1} \<A_2, \partial A_2 \> - 2\<A_2^\rho , \partial  (ad(\tau^{-1})A_1) \> + 2\<A_2^\rho, a_1 \partial(\tau^{-1} \partial\tau)\>}\\
+ \< ad(\tau^{-1})A_1 , \partial  (ad(\tau^{-1})A_1)\>
- 2\< ad(\tau^{-1})A_1 , a_1 \partial  (\tau^{-1} \partial\tau)\>\\
+ \< a_1 \tau^{-1} \partial\tau, a_1 \partial (\tau^{-1} \partial\tau)\> d\theta\\
\shoveleft{\phantom{\tau^*R} = \frac{i}{4\pi} \int_{S^1} \<A_2, \partial A_2 \> - 2\<A_2^\rho , \partial  (ad(\tau^{-1})A_1) \> 
+ 2\<A_2^\rho , a_1\partial (\tau^{-1} \partial\tau)\>}\\
+ \< ad(\tau^{-1})A_1 , \partial  (ad(\tau^{-1})A_1)\> - 2\< ad(\tau^{-1})A_1 , a_1 \partial  (\tau^{-1} \partial\tau)\> d\theta.
\end{multline*}
The last term vanishes since $a_1 \wedge a_1 = 0.$ For $d\epsilon$ we have:
\begin{multline*}
d\epsilon = \frac{i}{2\pi} d \int_{S^1}\< A_1 - \tfrac{1}{2} a_1 \tau^*Z, \tau^*Z \>d\theta\\
\shoveleft{\phantom{d\epsilon} = \frac{i}{2\pi} \int_{S^1} \< dA_1, \tau^*Z\>
 - \< A_1, d(\tau^*Z)\> 
 - \tfrac{1}{2}\< da_1\tau^*Z, \tau^*Z\>
 + \< a_1 \tau^*Z, d(\tau^*Z)\> d\theta}\\
%
\end{multline*}
and using the fact that $d(\tau^*Z) = ad(\tau)\partial (\tau^*\Theta^\rho),$
\begin{multline*}
\phantom{d\epsilon} = \frac{i}{2\pi} \int_{S^1}\< dA_1, \tau^*Z\>
 - \< A_1, ad(\tau)\partial (\tau^*\Theta^\rho)\>\\
\shoveright{- \tfrac{1}{2} \<da_1\tau^*Z, \tau^*Z\>
+ \< a_1 \tau^*Z, ad(\tau)\partial (\tau^*\Theta^\rho)\> d\theta}\\
\shoveleft{\phantom{d\epsilon} = \frac{i}{2\pi} \int_{S^1} \< dA_1, \tau^*Z\>}\\
- \< A_1, ad(\tau) \partial (A_2^\rho - ad(\tau^{-1})A_1
 + a_1 \tau^{-1} \partial\tau)\>
 - \tfrac{1}{2} \<da_1\tau^*Z, \tau^*Z\>\\
\shoveright{ + \< a_1 \tau^*Z, ad(\tau)\partial (A_2^\rho - ad(\tau^{-1})A_1 + a_1 \tau^{-1} \partial\tau)\> d\theta}\\
 %
%
\shoveleft{\phantom{d\epsilon} = \frac{i}{2\pi} \int_{S^1} \< dA_1, \tau^*Z\>
 - \< A_1, ad(\tau) \partial A_2^\rho \>
  + \< A_1, ad(\tau) \partial (ad(\tau^{-1})A_1) \>}\\
- \< A_1, a_1 ad(\tau) \partial (\tau^{-1} \partial\tau)\>
 - \tfrac{1}{2} \<da_1\tau^*Z, \tau^*Z\>\\
+ \< a_1 \tau^*Z, ad(\tau)\partial A_2^\rho \>
 - \< a_1 \tau^*Z, ad(\tau)\partial (ad(\tau^{-1})A_1) \> d\theta.
\end{multline*}
Therefore,
\begin{multline*}
\tau^*R - d\epsilon = \frac{i}{4\pi} \int_{S^1} \<A_2, \partial A_2 \>
-  2\< dA_1, \tau^*Z\>
- \< A_1, ad(\tau)\partial (ad(\tau^{-1})A_1) \>\\
+ 2\< a_1 \tau^{-1} \partial\tau, \partial(ad(\tau^{-1})A_1)\>
+ \<da_1\tau^*Z, \tau^*Z\> d\theta,
\end{multline*}
using the $ad$ invariance of the Killing form and integration by parts. Then, using the identity from before,
$$
\partial(ad(\tau^{-1})A) = ad(\tau^{-1})[A,\tau^*Z] + ad(\tau^{-1})\partial A,
$$
yields
\begin{multline*}
\tau^*R - d\epsilon = \frac{i}{4\pi} \int_{S^1} \<A_2, \partial A_2 \>
-  2\< dA_1, \tau^*Z\>
- \< A_1, [A_1,\tau^*Z] \>
- \< A_1 , \partial A_1\> \\
\shoveright{+ 2\< \tau^*Z a_1, [A_1,\tau^*Z] \>  
+ 2 \< a_1 \tau^*Z, \partial A_1\> 
+ \<da_1\tau^*Z, \tau^*Z\> d\theta}\\
\shoveleft{\phantom{\tau^*R - d\epsilon} = \frac{i}{4\pi} \int_{S^1} \<A_2, \partial A_2 \> 
- \< A_1 , \partial A_1\> 
-  2\< dA_1, \tau^*Z\> 
- \< [A_1, A_1],\tau^*Z \>} \\
+ 2 \< \tau^*Z a_1, \partial A_1\> 
+ \<da_1\tau^*Z, \tau^*Z\> d\theta.
\end{multline*}
Note now that if $(F,f)$ is the curvature of the connection $(A,a)$ then we have
\begin{align*}
(F,f)(X,Y)	&= (dA(X,Y) + \tfrac{1}{2}[(A,a)(X), (A,a)(Y)], da(X,Y))\\
		&= (dA(X,Y) + \tfrac{1}{2}([A(X), A(Y)] -a(X)\partial A(Y) + a(Y) \partial A(X)), da(X,Y)).
\end{align*}
That is,
$$
(F,f) = (dA + \tfrac{1}{2} [A,A] - a\wedge \partial A, da).
$$
Therefore, the formula above for $\tau^*R - d\epsilon$ reads
$$
\tau^*R - d\epsilon = \frac{i}{4\pi} \int_{S^1} \left\<\pi_1^*A, \partial \pi_1^*A \right\> - \left\< \pi_2^*A , \partial \pi_2^*A \right\> - 2\left\< \pi_2^* F - \tfrac{1}{2} \pi_2^* f \, \tau^*Z, \tau^*Z\right\> d\theta.
$$

\subsubsection{A curving for the lifting bundle gerbe}

Recall that in order to find the 3-curvature of the lifting bundle gerbe, and hence a representative for the image in real cohomology of the Dixmier-Douady class, we need a curving for $\tau^*\widehat{\LGS}.$ That is, some 2-form $B$ on $P$ such that $\d B = \tau^* R - d\epsilon.$ Note that $\d = \pi_1^* - \pi_2^*$ and
$$
\tau^*R - d\epsilon = \d \left(\frac{i}{4\pi} \int_{S^1} \left\< A, \partial A \right\> d\theta \right)
- \frac{i}{2\pi} \int_{S^1} \left\< \pi_2^* F - \tfrac{1}{2} \pi_2^* f \, \tau^*Z, \tau^*Z\right\> d\theta.
$$
To deal with the second term above, we use a similar method to the one in section \ref{S:LG string class}. Namely, we will need a Higgs field for the $\LGS$-bundle $P.$

\begin{definition}
A \emph{Higgs field} for $P$ is a map $\Phi \colon P \to L\fg$ satisfying
$$
\Phi(p (\gamma, \phi)) = \rho_{\phi}^{-1} \left(ad(\gamma^{-1}) \Phi(p) + \gamma^{-1} \partial \gamma \right).
$$
\end{definition}
We shall explain the geometric significance of this map in section \ref{S:Higgs fields, LGxS^1-bundles,...}. As in the $LG$ case, Higgs fields exist for $\LGS$-bundles. Note that the condition above implies that a Higgs field $\Phi$ satisfies
$$
\pi_1^*\Phi = \rho_{\tau}^{-1} \left( ad(\tau^{-1}) \pi_2^* \Phi + \tau^{-1} \partial \tau \right)
$$
or simply,
$$
ad(\tau) \Phi_2^\rho = \Phi_1 + \tau^*Z.
$$
Using this, the second term in $\tau^*R - d\epsilon$ becomes
$$
\frac{i}{2\pi} \int_{S^1} \left\< F_1 - \tfrac{1}{2} f_1 \, \tau^*Z, ad(\tau) \Phi_2^\rho - \Phi_1 \right\> d\theta.
$$
Since $(F,f)$ is a curvature, it satisfies
$$
\pi_1^* (F,f) = ad(\tau^{-1}) \pi_2^*(F,f).
$$
That is, $f_2 = f_1$ and
$$
F_2 = \rho_{\tau}^{-1} \left(ad(\tau^{-1}) F_1 - f_1 \tau^{-1}\partial \tau \right),
$$
or
$$
ad(\tau) F_2^\rho = F_1 -  f_1 \tau^*Z.
$$
Using this, we have
\begin{align*}
\frac{i}{2\pi} \int_{S^1} 	&\left\< F_1 - \tfrac{1}{2} f_1 \, \tau^*Z, ad(\tau) \Phi_2^\rho - \Phi_1 \right\> d\theta\\
	& = \frac{i}{4\pi} \int_{S^1}  \left\< F_1 + ad(\tau) F_2^\rho , ad(\tau) \Phi_2^\rho - \Phi_1 \right\> d\theta\\
	& = \frac{i}{4\pi} \int_{S^1}  \left\< F_1  , ad(\tau) \Phi_2^\rho  \right\> 
- \left\< F_1  ,  \Phi_1 \right\>
+ \left\<   F_2 ,  \Phi_2  \right\>
- \left\<  ad(\tau) F_2^\rho ,  \Phi_1 \right\> d\theta\\
	& = \d \left( \frac{i}{4\pi} \int_{S^1} \left< F, \Phi \right\> d\theta \right)
+ \frac{i}{4\pi} \int_{S^1} \left\< F_1  , ad(\tau) \Phi_2^\rho  \right\>
- \left\<  ad(\tau) F_2^\rho ,  \Phi_1 \right\> d\theta.
\end{align*}
Note, however, that the second integral above simplifies further
\begin{align*}
 \frac{i}{4\pi} \int_{S^1}	&\left\< F_1  , ad(\tau) \Phi_2^\rho  \right\>
- \left\<  ad(\tau) F_2^\rho ,  \Phi_1 \right\> d\theta\\
	& = \frac{i}{4\pi} \int_{S^1} \left\< ad(\tau) F_2 +  f_1 \tau^*Z, ad(\tau) \Phi_2^\rho  \right\>  - \left\< F_1 -  f_1 \tau^*Z,  \Phi_1 \right\> d\theta\\
	& = \frac{i}{4\pi}  \int_{S^1} \< F_2, \Phi_2 \> - \< F_1 ,  \Phi_1\> + \< f_1 \tau^*Z, ad(\tau) \Phi_2^\rho + \Phi_1\> d\theta\\
	& = \d \left( \frac{i}{4\pi} \int_{S^1} \left\< F, \Phi \right\> d\theta \right)  +  \frac{i}{4\pi} \int_{S^1} \left\< f_1 \tau^*Z,  2\Phi_1 + \tau^*Z \right\> d\theta\\
	& =\d \left( \frac{i}{4\pi} \int_{S^1} \left\< F, \Phi \right\> d\theta \right)  +  \frac{i}{4\pi} \int_{S^1} 2\left\< f_1\tau^*Z, \Phi_1 \right\> + \left\< f_1\tau^*Z, \tau^*Z \right\> d\theta.
\end{align*}
Therefore, $\tau^* R - d\epsilon$ is equal to
$$
\d \left( \frac{i}{4\pi} \int_{S^1} \left\< A, \partial A \right\> - 2 \left\< F, \Phi \right\> d\theta \right) - \frac{i}{4\pi} \int_{S^1} 2\left\< f_1\tau^*Z, \Phi_1 \right\> + \left\< f_1\tau^*Z, \tau^*Z \right\> d\theta.
$$
So it is enough to find a $B_2 \in \Omega^2(P)$ such that
$$
\d B_2 = \frac{i}{4\pi} \int_{S^1} 2\left\< f_1\tau^*Z, \Phi_1 \right\> + \left\< f_1\tau^*Z, \tau^*Z \right\> d\theta.
$$
Consider, then, the form
$$
\frac{i}{4\pi} \int_{S^1} \left\< \Phi, f \Phi \right\> d\theta.
$$
We have
\begin{multline*}
\d \left( \frac{i}{4\pi} \int_{S^1} \left\< \Phi, f\Phi \right\> d\theta \right)\\
\shoveleft{\phantom{\d (  \int_{S^1} \< \Phi}
= \frac{i}{4\pi} \int_{S^1} \left\< \Phi_2, f_2 \Phi_2 \right\> 
 - \left\< \Phi_1, f_1 \Phi_1 \right\> d\theta}\\
\shoveleft{\phantom{\d (  \int_{S^1} \< \Phi} 
= \frac{i}{4\pi} \int_{S^1} \left\< ad(\tau^{-1}) (\Phi_1 + \tau^*Z), f_1 ad(\tau^{-1}) (\Phi_1 + \tau^*Z) \right\> 
- \left\< \Phi_1, f_1 \Phi_1 \right\> d\theta}\\
\shoveleft{\phantom{\d(  \int_{S^1} \< \Phi} 
= \frac{i}{4\pi} \int_{S^1} \left\< \Phi_1, f_1 \Phi_1 \right\> 
+ \left\< \tau^*Z, f_1 \Phi_1 \right\> 
+ \left\< \Phi_1 , f_1 \tau^*Z \right\>}\\
\shoveright{ + \< \tau^*Z, f_1  \tau^*Z\> 
- \left\< \Phi_1, f_1 \Phi_1 \right\> d\theta}\\
\shoveleft{\phantom{\d(  \int_{S^1} \< \Phi} 
= \frac{i}{4\pi} \int_{S^1} 2\left\< f_1\tau^*Z, \Phi_1 \right\> 
+ \left\< f_1\tau^*Z, \tau^*Z \right\> d\theta.}\\
\end{multline*}
Therefore, a curving for the lifting bundle gerbe is given by
$$
B = \frac{i}{4\pi} \int_{S^1} \left\<A, \partial A \right\> - 2\< F + \tfrac{1}{2}f\Phi, \Phi\> \, d\theta.
$$

\subsubsection{The string class of an $\LGS$-bundle}

The last step now that we have found a curving for the lifting bundle gerbe is to calculate the $3$-curvature $H = dB.$ Then $H/2\pi i$ is integral and represents the real image of the Dixmier-Douady class of $\tau^*\widehat{\LGS}$ (and hence the obstruction to lifting $P$). We have
\begin{align*}
dB	&= \frac{i}{4\pi} \int_{S^1} \left\< dA, \partial A\right\> - \left\< A, d(\partial A)\right\> - 2\left\< dF, \Phi \right\> - 2 \left\< F, d\Phi \right\> - \left\< d\Phi, f\Phi\right\>  - \left\< \Phi, f d\Phi\right\> d\theta\\
	&= \frac{i}{4\pi} \int_{S^1} \left\< dA, \partial A\right\> + \left\< \partial A, dA\right\> - 2\left\< dF, \Phi\right\> - 2 \left\< F, d\Phi\right\> - 2 \left\< d\Phi, f\Phi\right\> d\theta\\
	&= \frac{i}{2\pi} \int_{S^1} \left\< dA, \partial A\right\> - \left\< dF, \Phi\right\> - \left\< F, d\Phi\right\> - \left\< d\Phi, f\Phi \right\> d\theta.
\end{align*}
To proceed further, we require the Bianchi identity for $(F, f).$ Note that
$$
d(F, f) = ( [dA, A] - f \wedge \partial A + a \wedge \partial (dA) , d^2 a).
$$
In particular, this means that
$$
dF = [F, A] - f \wedge \partial A + a \wedge \partial F,
$$
since
\begin{align*}
[F,A] &= [dA, A] + \tfrac{1}{2}[[A,A],A] - [a\wedge \partial A,A]\\
	&= [dA, A] -[a\wedge \partial A,A],
\end{align*}
and
\begin{align*}
a\wedge \partial F &= a \wedge \partial (dA) + \tfrac{1}{2} a\wedge \partial [A, A] - a\wedge\partial (a\wedge \partial A)\\
			& = a \wedge \partial (dA) + [a\wedge \partial A, A].
\end{align*}
Using this, and the fact that $\int_{S^1}\<[A, A], \partial A\> d\theta$ and  $\< a\wedge \partial A, \partial A\>$ both vanish (so that $\int_{S^1} \<dA, \partial A\> d \theta = \int_{S^1} \<F, \partial A\> d \theta$), the expression for $dB$ becomes
\begin{multline*}
dB = \frac{i}{2\pi} \int_{S^1} \left\< F, \partial A \right\> 
- \left\< [F,A] - f \wedge \partial A + a \wedge \partial F, \Phi\right\> 
-  \left\< F, d\Phi\right\> 
- \left\< d\Phi, f\Phi\right\> d\theta\\
\shoveleft{\phantom{dB} = \frac{i}{2\pi} \int_{S^1} \left\< F, \partial A\right\> 
- \left\< [F,A], \Phi\right\> 
+\left\< f \wedge \partial A, \Phi\right\> 
-\left\< a \wedge \partial F, \Phi\right\> 
-  \left\< F, d\Phi\right\> 
- \left\< d\Phi, f\Phi\right\> d\theta}\\
\shoveleft{\phantom{dB} = \frac{i}{2\pi} \int_{S^1} \left\< F+ f\Phi, \partial A\right\> 
- \left\< F,[A, \Phi]\right\>  
- \left\< a \wedge \partial F, \Phi\right\> 
-  \left\< F+ f\Phi, d\Phi\right\> d\theta}\\
\shoveleft{\phantom{dB} = \frac{i}{2\pi} \int_{S^1} \left\< F+ f\Phi, \partial A\right\> 
- \left\< F,[A, \Phi]\right\>  
+ \left\< F,a \partial \Phi\right\> 
-  \left\< F+ f\Phi, d\Phi\right\> d\theta}\\
\shoveleft{\phantom{dB} = \frac{i}{2\pi} \int_{S^1} \left\< F+ f\Phi, \partial A - [A, \Phi]  + a \partial \Phi -   d\Phi\right\> d\theta.}\\
\end{multline*}
Where the last line follows from the fact that $\int_{S^1} \< f\Phi, a \partial \Phi\> d\theta$ and $\< f\Phi , [A, \Phi]\>$ both vanish. If we define the covariant derivative of $\Phi$ by
$$
\nabla \Phi = d\Phi + [A, \Phi] - \partial A - a \partial \Phi,
$$
then one can easily check that it is (twisted) equivariant for the adjoint action. That is,
$$
\nabla \Phi (X (\gamma, \phi)) = \rho_{\phi}^{-1} \left(ad(\gamma^{-1}) \nabla \Phi (X) \right),
$$
for any tangent vector $X.$ The same is true for the quantity $F + f \Phi,$ and so using the $ad$-invariance of the Killing form and the rotation invariance of the integral, Lemma \ref{L:dpsi=Dpsi} implies that $H= dB$ descends to a form on $M.$ Thus we have proven

\begin{theorem}\label{T:LGxS^1string}
Let $P \to M$ be a principal $LG\rtimes S^1$-bundle and let $\Phi$ be a Higgs field for $P$ and $(A,a)$ be a connection for $P$ with curvature $(F,f).$ Then the string class of $P,$ that is, the obstruction to lifting $P$ to an $\widehat{\LGS}$-bundle, is represented in de Rham cohomology by
$$
-\frac{1}{4\pi^2}\int_{S^1} \langle F+ f\Phi, \nabla\Phi \rangle d\theta,
$$
where
$$
\nabla\Phi = d\Phi + [A,\Phi] - \partial A - a\partial \Phi.
$$
\end{theorem}

\subsection{Reduced splittings for lifting bundle gerbes}\label{SS:reduced splittings}

In this section we shall present an alternative method for finding the curving of a lifting bundle gerbe and show how to apply this to the problem above. This method uses \emph{reduced splittings} and was first introduced by Gomi \cite{Gomi:2003}.

In \cite{Brylinski:1993} Brylinski considers the problem of lifting a principal $\cG$-bundle $P$ to a $\widehat{\cG}$-bundle $\widehat{P},$ for which he uses a \emph{bundle splitting}. He relates the obstruction class to the \emph{scalar curvature} of a certain connection on $\widehat{P}.$ In \cite{Gomi:2003} Gomi phrases this in such a way that he can use the theory of lifting bundle gerbes in order to calculate the obstruction class. We shall begin by briefly outlining Brylinski's results before describing the reduced splittings of Gomi.

Let $\cG$ be a Lie group with central extension $\widehat{\cG}.$ If $\mathfrak{G}$ and $\widehat{\mathfrak{G}}$ are the Lie algebras of $\cG$ and $\widehat{\cG}$ respectively then we have an extension of Lie algebras
$$
0 \to i\RR \to \widehat{\mathfrak{G}} \to \mathfrak{G} \to 0.
$$
We can define an action of $\cG$ on $\widehat{\mathfrak{G}}$ by lifting the adjoint action of $\cG$ on its Lie algebra. That is, we define $ad \colon \cG \times \widehat{\mathfrak{G}} \to \widehat{\mathfrak{G}}$ by
$$
ad(g) \hat{\xi} = ad(\hat{g}) \hat{\xi},
$$
where $\hat{\xi} \in \widehat{\mathfrak{G}}$ and $\hat{g} \in\widehat{\cG}$ projects to $g \in \cG.$ This is well-defined since $U(1)$ acts trivially on $\widehat{\mathfrak{G}}$ and any two lifts of $g$ differ by an element of $U(1).$ Consider now a principal $\cG$-bundle $P.$ We can write down an exact sequence of vector bundles associated to $P$ as follows. Let $\Ad_{\fg}(P)$ denote the adjoint bundle of $P$ where $\cG$ acts on the Lie algebra $\fg.$ For example, $\Ad_{\mathfrak{G}}(P)$ is the usual adjoint bundle of $P$ and $\Ad_{i\RR}(P) = P \times_{ad} i\RR \simeq M \times i\RR.$ Since $\cG$ acts via the adjoint action on the exact sequence above, we have an exact sequence of vector bundles
$$
0 \to \Ad_{i\RR}(P) \to \Ad_{\mathfrak{G}}(P) \to \Ad_{\widehat{\mathfrak{G}}}(P) \to 0.
$$
This means that $\Ad_{\mathfrak{G}}(P)$ is isomorphic to the direct sum of $M\times i \RR$ and $\Ad_{\widehat{\mathfrak{G}}}(P).$ A choice of isomorphism is called a \emph{bundle splitting}. That is,

\begin{definition}[\cite{Brylinski:1993}]
A \emph{bundle splitting} of $P$ is a vector bundle map
$$
L \colon \Ad_{\widehat{\mathfrak{G}}}(P) \to \Ad_{i\RR} (P)
$$
which is the identity on the (trivial) subbundle $\Ad_{i\RR} (P).$
\end{definition}

As mentioned above, Brylinski uses the notion of scalar curvature to calculate the obstruction to the existence of a lift of $P.$ This is essentially the $i\RR$ part of the curvature of a connection on $\widehat{P}.$ More precisely,

\begin{definition}[\cite{Brylinski:1993}]
Let $\hat{A}$ be a connection on $\widehat{P}$ with curvature $\hat{F},$ viewed as a $2$-form on $M$ with values in $\Ad_{\widehat{\mathfrak{G}}}(\widehat{P}) \simeq \Ad_{\widehat{\mathfrak{G}}}(P).$ Let $L$ be a bundle splitting of $P.$ The \emph{scalar curvature} of $\hat{A}$ is the $i\RR$-valued $2$-form
$$
K = L \circ \hat{F}.
$$

\end{definition}

To see how this is related to the obstruction class, let $\{U_\a\}$ be a good cover of $M$ over which $P$ is trivial. Then there exists a lift $\widehat{P}_\a$ of $P |_{U_\a} \to U_\a.$ Choose a connection $A_\a$ on $P |_{U_\a}$ and let $K_\a$ be the scalar curvature of a connection $\hat{A}_\a$ on $\widehat{P}_\a$ which is compatible with $A_\a$ in the sense that the pull-back of $A_\a$ to $\widehat{P}_\a$ coincides with the image of $\hat{A}_\a$ in $\mathfrak{G}.$ That is,
$$
f^* A_\a = p ( \hat{A}_\a),
$$
where $f$ is the bundle map $\widehat{P}_\a \to P |_{U_\a}$ and $p$ is the projection $\widehat{\cG} \to \cG.$ Brylinski's result, then, is that the (real image of the) obstruction class restricted to $U_\a$ coincides with the derivative of the scalar curvature, $dK_\a.$

As mentioned, Gomi's results interpolate between the method described above and the theory of lifting bundle gerbes which we have used extensively. He utilises so-called reduced splittings to write down a formula for the curving of the lifting bundle gerbe associated to a lifting problem and relates the curving to the scalar curvature. In the case where a splitting of the Lie algebra of $\widehat{\mathfrak{G}}$ has been specified, reduced splittings are equivalent to bundle splittings. To describe Gomi's results, let us assume we have chosen a splitting of the Lie algebra $\widehat{\mathfrak{G}}$ as $\mathfrak{G} \oplus i\RR.$

\begin{definition}[\cite{Gomi:2003}]
The \emph{group cocycle} for the central extension $\widehat{\cG}$ is the map $\sigma \colon \cG \times \mathfrak{G} \to i\RR$ defined by
$$
\sigma(g, \xi) = ad(g) (\xi, 0) - (ad(g) \xi, 0),
$$
where $ad(g)$ acts on $\widehat{\mathfrak{G}}$ as described above.
\end{definition}

The group cocycle gives information about the multiplication in $\widehat{\cG}$ in the same way as the $1$-form $\a$ which we used. In fact, as we shall see, to apply Gomi's results to the case where $\cG$ is either the loop group $LG$ or the semi-direct product $\LGS,$ we shall give $\sigma$ in terms of $\a.$

\begin{definition}[\cite{Gomi:2003}]
A \emph{reduced splitting} for a principal $\cG$-bundle $P$ is a map $\ell \colon P \times \mathfrak{G} \to i\RR$ which is linear in the second factor and satisfies
$$
\ell (p, \xi) = \ell (p g, ad(g^{-1})\xi) + \sigma(g^{-1}, \xi).
$$
\end{definition}

The relation between reduced splittings and bundle gerbe curvings is given by the following theorem.

\begin{theorem}[\cite{Gomi:2003}]\label{T:Gomi}
Let $F$ be the curvature of a connection $A$ on $P$ and $\ell$ be a reduced splitting for $P.$ Define a $2$-form $\kappa$ on $P$ by $\kappa_p = \ell (p, F).$ Then a curving for the lifting bundle gerbe associated to the lifting problem for $P$ is given by
$$
B = \frac{1}{2} \omega (A, A) + \kappa,
$$
where $\omega(\xi, \zeta) = [(\xi,0), (\zeta,0)]_{\widehat{\mathfrak{G}}} - ([\xi, \zeta]_{\mathfrak{G}}, 0)$ is the cocycle classifying the central extension.
\end{theorem}

To connect this with Brylinski's work, Gomi proves the following theorem relating the curving and the scalar curvature.

\begin{theorem}[\cite{Gomi:2003}]
Let $P$ be a principal $\cG$-bundle and $\widehat{P}$ be a lift of $P.$ Let $A$ be a connection on $P$ and $\hat{A}$ be a compatible connection on $\widehat{P}.$ Then the curving can be written as
$$
B = \pi^*K - \tilde{F},
$$
where $\tilde{F}$ is the curvature of the connection $\hat{A} - f^*A$ on $\widehat{P}$ (for $f \colon \widehat{P} \to P$ the bundle map defining the lift) and $K$ is the scalar curvature of $\hat{A}.$
\end{theorem}

We would now like to consider the case where $\cG = \LGS.$ We shall define a reduced splitting for $P$ so we can use Theorem \ref{T:Gomi} to calculate a curving and show that it is in agreement with the results from section \ref{SS:LGxS^1 string class}. The group cocycle in this case is given by
\begin{align*}
\sigma ((\gamma, \phi)^{-1}, (\xi, x))	& = \a_{((1,1), (\gamma, \phi))} ((\xi, x), (0,0))\\
			& = \frac{i}{2\pi} \int_{S^1} \left\< \xi - \tfrac{1}{2} x \partial \gamma \gamma^{-1}, \partial \gamma \gamma^{-1} \right\> d\theta,
\end{align*}
and we have

\begin{proposition}
A reduced splitting for the $\LGS$-bundle $P$ is given by
$$
\ell (p, (\xi, x)) = - \frac{i}{2\pi} \int_{S^1} \left\< \xi + \tfrac{1}{2}x \, \Phi(p), \Phi(p) \right\> d\theta,
$$
where $\Phi$ is a Higgs field for $P.$
\end{proposition}

\begin{proof}
We need only show that it satisfies the transformation property above. We can calculate
\begin{multline*}
\ell (p (\gamma, \phi), ad(\gamma, \phi)^{-1} (\xi, x))\\
\shoveleft{\phantom{\ell (p } = - \frac{i}{2\pi} \int_{S^1} \left\< ad(\gamma^{-1}) (\xi - x Z) + \tfrac{1}{2} x\, ad(\gamma^{-1})(\Phi(p) +  Z), ad(\gamma^{-1})(\Phi(p) +  Z) \right\> d\theta}\\
\shoveleft{\phantom{\ell (p } = - \frac{i}{2\pi} \int_{S^1} \left\< \xi, \Phi(p) \vphantom{\tfrac{1}{2}} \right\> 
- \left\< x Z , \Phi(p) \vphantom{\tfrac{1}{2}} \right\>
+ \left\< \tfrac{1}{2} x\, \Phi(p), \Phi(p) \right\> 
+ \left\< \tfrac{1}{2} x Z , \Phi(p) \right\> }\\
\shoveright{+ \left\< \xi, Z \vphantom{\tfrac{1}{2}} \right\> 
- \left\< xZ , Z \vphantom{\tfrac{1}{2}} \right\>
+ \left\< \tfrac{1}{2} x\, \Phi(p) , Z \right\>
+ \left\< \tfrac{1}{2} x\, Z, Z \right\> d\theta}\\
\shoveleft{\phantom{\ell (p } = - \frac{i}{2\pi} \int_{S^1} \left\< \xi + \tfrac{1}{2} x \, \Phi(p) , \Phi (p) \right\>
+ \left\< X - \tfrac{1}{2}x Z, Z \right\> d\theta}\\
\shoveleft{\phantom{\ell (p } = \ell(p, (\xi,x)) - \sigma((\gamma, \phi)^{-1} , (\xi, x))}\\
\end{multline*}
as required.
\end{proof}

Note that in order to use Theorem \ref{T:Gomi}, we need the cocycle $\omega.$ This is simply given by the form $R$ which defines the central extension. In particular,
$$
\omega((\xi, x), (\zeta, y)) = \frac{i}{2\pi} \int_{S^1} \left\< \xi, \partial \zeta \right\> d\theta.
$$
Therefore, for the curving of the lifting bundle gerbe, Theorem \ref{T:Gomi} gives
\begin{align*}
B	&= \frac{1}{2} R((A, a), (A, a)) - \frac{i}{2\pi} \int_{S^1} \left\< F + \tfrac{1}{2}f \, \Phi, \Phi \right\> d\theta\\
	&= \frac{i}{4\pi} \int_{S^1} \left\< A, \partial A \right\> - 2 \< F + \tfrac{1}{2}f \, \Phi, \Phi\> \, d\theta,
\end{align*}
where as before, $(A,a)$ is a connection on $P$ and $(F, f)$ is its curvature.

\section{Higgs fields, $LG\rtimes S^1$-bundles and the string class}\label{S:Higgs fields, LGxS^1-bundles,...}

Now that we have an explicit formula for the string class of an $\LGS$-bundle $P,$ it is natural to ask whether there is some relation with the Pontrjagyn class of a $G$-bundle related to $P$ in some way, as was the case with the string class of an $LG$-bundle presented in chapter \ref{C:background}. In particular, in section \ref{S:Higgs fields, LG-bundles,...}, following \cite{Murray:2003}, we saw that there was a correspondence between $LG$-bundles over $M$ (with connection and Higgs field) and $G$-bundles over $M\times S^1$ (with connection) (Propositions \ref{P:LG correspondences} and \ref{P:LG connection correspondences}) and we used this to prove that the string class of $P$ is given by integrating over the circle the first Pontrjagyn class of the corresponding $G$-bundle (Theorem \ref{T:LGPont}). In this section, we shall show there is a correspondence between $\LGS$-bundles over $M$ and $G$-bundles over $S^1$-bundles over $M,$ which holds on the level of connections as well. As in section \ref{S:Higgs fields, LG-bundles,...} we shall use this correspondence to prove that the string class of $P$ is given in terms of the Pontrjagyn class of some $G$-bundle.

\subsection{Higgs fields and $LG\rtimes S^1$-bundles}

The following correspondence first appeared in \cite{Bergman:2005}. We will present it here in detail and also extend it to the level of connections.

Suppose that we have a principal $G$-bundle over a principal $S^1$-bundle:
$$
\xymatrix{ \widetilde{P}\ar^{G}[d] \\
		Y\ar^{S^1}[d]\\
		M}
$$
We would like to mimic the construction of the $LG$-bundle in section \ref{S:Higgs fields, LG-bundles,...} where we essentially took loops in $\widetilde{P}$ such that their image in $M\times S^1$ commuted with the obvious $S^1$ action on this space. That is, for a loop $f \in L\widetilde{P}_m$ in the fibre above $\{m\}\times S^1$ we required that $\tilde{\pi}(f(\theta)) = (m,\theta).$ The difference here is that we cannot choose a global section $M\to Y$ and thus there is no way of choosing a `starting point' for the loop $\tilde{\pi}(f) \colon S^1 \to Y.$ We can, however, still require that the map $\tilde{\pi}(f)$ commutes with the $S^1$ action on $Y$ (which we will write as addition). That is, we can define
$$
P = \{ f \colon S^1 \to \widetilde{P} \mid \tilde{\pi}(f(\theta + \phi)) = \tilde{\pi}(f(\theta)) + \phi \}
$$
and there is a canonical map $P \to M.$ $P$ is acted on by $LG\rtimes S^1:$
$$
(f (\gamma, \phi))(\theta) = f(\theta +\phi)\gamma(\theta +\phi),
$$
i.e.
$$
f(\gamma, \phi) = \rho_{\phi}^{-1}(f\gamma).
$$
It is a right action since
\begin{align*}
f (\gamma_1, \phi_1)(\gamma_2, \phi_2) &=  \rho_{(\phi_1 +\phi_2)}^{-1}f \rho_{(\phi_1+\phi_2)}^{-1}\gamma_1 \rho_{\phi_2}^{-1}( \gamma_2)\\
			&= \rho_{(\phi_1 +\phi_2)}^{-1}(f \gamma_1 \rho_{\phi_1}( \gamma_2))\\
			&= f (\gamma_1\rho_{\phi_1}(\gamma_2), \phi_1 + \phi_2).
\end{align*}
It preserves the fibres of $P$ since the $G$ action on $\widetilde{P}$ preserves fibres and the $S^1$ action on $Y$ preserves fibres. It is also free and transitive on fibres and therefore $P\to M$ is a principal $LG\rtimes S^1$-bundle. Note that local triviality of this bundle follows from the local triviality of $Y$ as follows: Choose a good cover of $M$ and let $U$ be an open set such that we can find a local section $s \colon U \to Y_{|_U}.$ There is a map $P \to Y$ given by $f\mapsto \tilde{\pi}(f(0)).$ If we pull-back $P$ by $s$ then $s^* P \to U$ is trivial (since $U$ is contractible).

Conversely, suppose we are given a principal $LG\rtimes S^1$-bundle $P \to M.$ Following the construction in section \ref{S:Higgs fields, LG-bundles,...}, define
$$
\widetilde{P} = (P\times G\times S^1)/LG\rtimes S^1,
$$
where $[p,g,\theta] = [p (\gamma, \phi), \gamma(\theta)^{-1}g, \theta - \phi].$ A $G$ action on $\widetilde{P}$ is given by $[p, g, \theta] h = [p, gh, \theta].$ There is a natural projection from $\widetilde{P}$ to the $S^1$-bundle associated to $P$ via the homomorphism $LG\rtimes S^1\to S^1,$ that is, $\widetilde{P} \to (P\times S^1)/LG\rtimes S^1 \simeq P/ LG,$ given by $\tilde{\pi}([p, g, \theta]) = [p, \theta].$ This makes $\widetilde{P}$ into a principal $G$-bundle. Thus, given the $LG\rtimes S^1$-bundle $P \to M$ we can construct a $G$-bundle over an $S^1$-bundle:
$$
\xymatrix{ \dfrac{P\times G \times S^1}{LG\rtimes S^1}\ar^{G}[d] \\
		\dfrac{P\times S^1}{LG\rtimes S^1}\ar^{S^1}[d]\\
		M}
$$
We would like to show that both constructions above are invertible (as we did for the constructions in the $LG$ case). Assume, then, that we are given an $LG\rtimes S^1$-bundle $P\to M$ and have constructed the $G$-bundle $\widetilde{P}$ over the $S^1$-bundle $P/LG \to M$ as above. Then use the first correspondence above to form the $LG\rtimes S^1$-bundle $P' \to M$ (by taking certain loops in $\widetilde{P}$). So we have
$$
P' = \{f \colon S^1 \to (P\times G\times S^1)/LG\rtimes S^1 \mid \tilde{\pi}(f(\theta + \phi)) = \tilde{\pi}(f(\theta)) + \phi \}
$$
and a bundle isomorphism is given by
$$
P \to P' ; \quad p \mapsto f_p = (\theta \mapsto [p, 1, \theta]).
$$
It is easily checked that this map commutes with the $LG\rtimes S^1$ action, for
\begin{align*}
p (\gamma, \phi) \mapsto &~ f_{p(\gamma, \phi)}\\
					= &~ [p(\gamma, \phi), 1, \theta]
\intertext{and on the other hand,}
f_p (\gamma, \phi) 	= &~ f_p (\theta + \phi) \gamma ( \theta + \phi)\\
					= &~ [p, 1, \theta+\phi] \gamma(\theta + \phi)\\
					= &~ [p, \gamma(\theta+\phi), \theta+\phi]\\
					= &~ [p(\gamma, \phi), 1, \theta].
\end{align*}
So we have that $P \simeq P'.$ If, on the other hand, we are given the $G$-bundle over the $S^1$-bundle $\widetilde{P} \to Y\to M$ and have constructed $P \to M,$ then we can construct $\widetilde{P}' \to P/LG \to M$ and we would like for these bundles to be isomorphic. That is, $\widetilde{P}' \simeq \widetilde{P}$ and $P/LG \simeq Y.$ Firstly, consider the map $P/LG \simeq P\times_{LG\rtimes S^1} S^1 \to Y$ defined by $[f, \theta] \mapsto \tilde{\pi}(f(\theta)).$ This is well-defined on equivalence classes:
$$
[f, \theta] = [\rho_{\phi}^{-1}(f\gamma), \theta -\phi] \mapsto \tilde{\pi}(f(\theta -\phi + \phi)\gamma(\theta - \phi +\phi )) = \tilde{\pi}(f(\theta)).
$$
It commutes with the $S^1$ action on $P\times_{LG\rtimes S^1} S^1:$
$$
[f, \theta+ \a] \mapsto \tilde{\pi}(f(\theta + \a)) = \tilde{\pi}(f(\theta)) + \a
$$
by the definition of $P$ in terms of $\widetilde{P}.$ Thus $P\times_{LG\rtimes S^1} S^1 \simeq Y.$ For $\widetilde{P}'$ and $\widetilde{P}$ consider the bundle map
$$
\widetilde{P}' \to \widetilde{P}; \quad [f, g, \theta] \mapsto f(\theta)g.
$$
This is well-defined:
$$
[f, g, \theta] = [f(\gamma, \phi), \gamma(\theta)^{-1}g, \theta - \phi] \mapsto f(\theta - \phi- \phi)\gamma(\theta - \phi + \phi)\gamma(\theta)^{-1} g = f(\theta)g,
$$
and commutes with the $G$ action:
$$
[f, g,\theta]h = [f, gh, \theta] \mapsto f(\theta)gh = (f(\theta)g)h.
$$
Therefore, it is a bundle isomorphism and $\widetilde{P}' \simeq \widetilde{P}.$ Thus we have proven

\begin{proposition}[\cite{Bergman:2005}]\label{P:LGxS^1 correspondences}
There is a bijective correspondence between isomorphism classes of principal $\LGS$-bundles over $M$ and isomorphism classes of principal $G$-bundles over principal $S^1$-bundles over $M.$
\end{proposition}

As in section \ref{S:Higgs fields, LG-bundles,...} the correspondences here hold on the level of connections as well. We shall now describe how to derive the connections corresponding to one another.

Suppose we are given a connection $\tilde{A}$ on $\widetilde{P}\to Y$ and a connection $\tilde{a}$ on $Y\to M.$ This amounts to a splitting of the tangent spaces $T_{\tilde{p}}\widetilde{P} \simeq V_{\tilde{p}}\widetilde{P} \oplus H_{\tilde{p}}\widetilde{P}$ at each point $\tilde{p} \in \widetilde{P}$ and also $T_yY \simeq  V_y Y \oplus H_y Y$ at each point $y\in Y.$ Since $P$ is given by certain loops in $\widetilde{P},$ a vector $X\in T_f P$ is really a vector field along $f$ in $\widetilde{P}.$ So, $X_\theta \in T_{f(\theta)} \widetilde{P}.$ Thus we can use the splittings of the tangent spaces of $\widetilde{P}$ and $Y$ to define a splitting for the tangent space to $P$ at $f$ for each $\theta.$ So we have
\begin{align*}
T_{f(\theta)} \widetilde{P}	&\simeq V_{f(\theta)} \widetilde{P} \oplus H_{f(\theta)} \widetilde{P}\\
					& \simeq V_{f(\theta)} \widetilde{P} \oplus V_{\tilde{\pi}(f(\theta))}Y \oplus H_{\tilde{\pi}(f(\theta))} Y\\
					& \simeq V_{f(\theta)} \widetilde{P} \oplus V_{\tilde{\pi}(f(\theta))}Y \oplus T_{(\pi_Y\circ\tilde{\pi})(f(\theta))} M,
\end{align*}
using the isomorphisms $H_{f(\theta)}\widetilde{P} \simeq T_{\tilde{\pi}(f(\theta))} Y$ and $H_{\tilde{\pi}(f(\theta))} Y \simeq T_{(\pi_Y\circ\tilde{\pi})(f(\theta))} M.$ We can find the 1-form for this connection by calculating
$$
X_\theta - \widehat{\widehat{\pi_* X_\theta}}
$$
which equals $\iota_{f(\theta)}(A_f (X)_\theta),$ where $\pi = \pi_Y \circ \tilde{\pi}$ and $\widehat{\widehat{V}}$ is the horizontal lift of a vector on $M$ first to $Y,$ then to $\widetilde{P}.$ Note that using the connection on $Y$ we have
$$
\iota_{\tilde{\pi}(f(\theta))}(\tilde{a} (\tilde{\pi}_* X_\theta)) = \tilde{\pi}_* X_\theta - \widehat{\pi_* X_\theta},
$$
and so
$$
\widehat{\pi_* X_\theta} = \tilde{\pi}_* X_\theta -\iota_{\tilde{\pi}(f(\theta))}(\tilde{a} (\tilde{\pi}_* X_\theta)).
$$
Lifting everything, we have
$$
\widehat{\widehat{{\pi_* X_\theta}}} = \widehat{\tilde{\pi}_* X_\theta} - \widehat{\iota_{\tilde{\pi}(f(\theta))}(\tilde{a} (\tilde{\pi}_* X_\theta))},
$$
and thus
$$
\iota_{f(\theta)}(A_f (X)_\theta) = X_\theta - \widehat{\tilde{\pi}_* X_\theta} + \widehat{\iota_{\tilde{\pi}(f(\theta))}(\tilde{a} (\tilde{\pi}_* X_\theta))}.
$$
But $X_\theta - \widehat{\tilde{\pi}_* X_\theta} = \iota_{f(\theta)} (\tilde{A}(X_\theta))$ and so we have
$$
\iota_{f(\theta)}(A_f (X)_\theta) =  \iota_{f(\theta)} (\tilde{A}(X_\theta)) + \widehat{\iota_{\tilde{\pi}(f(\theta))}(\tilde{a} (\tilde{\pi}_* X_\theta))}.
$$
To make use of this we need to be able to write $A$ as an $L\fg$-valued 1-form and an $i\RR$-valued 1-form. That is, $A(X)_\theta = (\xi(\theta), x)$ for $\xi \in L\fg$ and $x \in i\RR.$ To that end, consider the vertical vector $V$ in $T_f P$ generated by the Lie algebra element $(\xi, x):$
\begin{align*}
V_\theta	&= \frac{d}{dt}_{|_0} f (\exp(t\xi), tx) (\theta)\\
		&= \frac{d}{dt}_{|_0} f(\theta + tx)\exp(t\xi(\theta + tx))\\
		&= \frac{d}{dt}_{|_0} \left(f(\theta) +f'(\theta) tx \vphantom{^2}\right) \left(1 + t\xi(\theta) + O(t^2)\right)\\
		&= \iota_{f(\theta)}(\xi(\theta)) + x f'(\theta).
\end{align*}
Since $A$ is a connection, it returns the Lie algebra element corresponding to the vertical part of a vector $X.$ Therefore, we must solve the following equation for $\xi$ and $x:$
$$
\iota_{f(\theta)} (\tilde{A}(X_\theta)) + \widehat{\iota_{\tilde{\pi}(f(\theta))}(\tilde{a} (\tilde{\pi}_* X_\theta))} = \iota_{f(\theta)}(\xi(\theta)) + x f'(\theta).
$$
Applying $\tilde{A}$ to both sides gives
$$
\tilde{A}(X_\theta)  = \xi(\theta) + x \tilde{A}(f'(\theta)),
$$
since $\widehat{\iota_{\tilde{\pi}(f(\theta))}(\tilde{a} (\tilde{\pi}_* X_\theta))}$ is horizontal with respect to $\tilde{A}.$ Thus, we have
$$
\xi(\theta) = \tilde{A} (X_\theta - x f'(\theta)).
$$
Taking instead, $\tilde{\pi}_*$ of both sides gives
$$
\iota_{\tilde{\pi}(f(\theta))}(\tilde{a} (\tilde{\pi}_* X_\theta)) = x\, \tilde{\pi}_* f'(\theta),
$$
since the vectors $\iota_{f(\theta)} (\tilde{A}(X_\theta))$ and $\iota_{f(\theta)}(\xi(\theta))$ are vertical in $\widetilde{P}$. Then applying $\tilde{a}$ to both sides yields
$$
\tilde{a} (\tilde{\pi}_* X_\theta) = x\, \tilde{a}(\tilde{\pi}_* f'(\theta)).
$$
So (with a slight abuse of notation) we can write the connection form on $P$ as
$$
(A, a)_f(X)_\theta = (\tilde{A} (X_\theta - a(X) f'(\theta)), a(X)),
$$
where $\tilde{A}$ and $\tilde{a}$ are connection forms on $\widetilde{P}$ and $Y$ respectively and $a(X)$ is given by the formula for $x$ above. Now that we have the connection on $P$ in this form we can check explicitly that it satisfies the conditions for a connection. By construction, it satisfies $(A, a)(\iota_f(\xi, x)) = (\xi ,x)$ and so we just need to check that $(A, a)(X(\gamma, \phi)) = ad(\gamma, \phi)^{-1} (A, a)(X).$ Recall that the adjoint action of $LG\rtimes S^1$ on its Lie algebra is given by
$$
ad(\gamma, \phi)^{-1} (\xi , x) = \left( \rho_{\phi}^{-1} \left(ad(\gamma^{-1})\xi - \gamma^{-1}\partial \gamma\, x\right), x\right)
$$
and so
$$
ad(\gamma, \phi)^{-1} (A, a)(X)_\theta = \left( \rho_{\phi}^{-1} (ad(\gamma^{-1}) \tilde{A} (X_\theta - a(X) f'(\theta)) - \gamma^{-1}\partial \gamma\, a(X)), a(X)\right).
$$
On the other hand, the action of $LG\rtimes S^1$ on the tangent vector $X$ is
$$
X(\gamma, \phi) = \rho_{\phi}^{-1}(X\gamma).
$$
Therefore,
\begin{align*}
(A, a)(X(\gamma, \phi))_\theta	&= \left(\tilde{A} (X(\gamma, \phi)_\theta - a(X(\gamma, \phi))\partial (f(\theta + \phi)\gamma(\theta + \phi))), a(X(\gamma, \phi))\right)\\
						&= \left(\tilde{A} (\rho_{\phi}^{-1}(X\gamma)_\theta - a(\rho_{\phi}^{-1}(X\gamma)_\theta)\partial (f(\theta + \phi)\gamma(\theta + \phi))), a(\rho_{\phi}^{-1}(X\gamma)_\theta)\right)\\
						\begin{split}
						&= \left(\tilde{A} (\rho_{\phi}(X\gamma)_\theta - a(\rho_{\phi}^{-1}(X\gamma)_\theta)\{ \partial f(\theta + \phi)\gamma(\theta + \phi)\right.\\
						&\left. \phantom{(\tilde{A} (\rho_{\phi}^{-1}(X\gamma)_\theta - a(\rho_{\phi}^{-1}(X\gamma} + f(\theta + \phi) \partial \gamma(\theta + \phi) \}), a(\rho_{\phi}^{-1}(X\gamma)_\theta)\right).
						\end{split}
\end{align*}
Since $\tilde{A}$ is a connection, we have $\tilde{A} (\rho_{\phi}^{-1}(X\gamma))_\theta = \rho_{\phi}^{-1} (ad(\gamma^{-1})\tilde{A}(X))_\theta$ and $\tilde{A} (\partial f(\theta + \phi)\gamma(\theta + \phi)) = \rho_{\phi}^{-1} (ad(\gamma^{-1})\tilde{A}(\partial f(\theta))).$ Also, since $a$ is $i\RR$-valued, we have $a(\rho_{\phi}^{-1}(X\gamma)_\theta) = a(X).$ Therefore,
\begin{multline*}
(A, a)(X(\gamma, \phi))_\theta	= \left(\rho_{\phi}^{-1} (ad(\gamma^{-1})\tilde{A}(X_\theta - a(X)\partial f(\theta)))\right.\\
\left. - a(X)\tilde{A}(f(\theta + \phi) \partial \gamma(\theta + \phi))), a(X)\vphantom{\tilde{A}}\right).
\end{multline*}
But, $f(\theta + \phi) \partial \gamma(\theta + \phi)$ is really shorthand for $\iota_{f(\theta + \phi)}(\rho_{\phi}^{-1}(\gamma^{-1} \partial \gamma))$ and so
$$
\tilde{A}(f(\theta + \phi) \partial \gamma(\theta + \phi)) = \tilde{A}(\iota_{f(\theta + \phi)}(\rho_{\phi}^{-1}(\gamma^{-1} \partial \gamma))) = \rho_{\phi}^{-1}(\gamma^{-1} \partial \gamma).
$$
Thus, we have
$$
(A, a)(X(\gamma, \phi))_\theta	= \left(\rho_{\phi}^{-1} (ad(\gamma^{-1})\tilde{A}(X_\theta - a(X)\partial f(\theta)) - a(X)\gamma^{-1} \partial \gamma), a(X)\vphantom{\tilde{A}}\right),
$$
as required.

As for the $LG$-bundle case in section \ref{S:Higgs fields, LG-bundles,...}, to define a connection\footnote{Of course, here we need to define a connection on $Y$ as well as on $\widetilde{P}.$} on $\widetilde{P}$ given the bundle $P$ we need a connection and Higgs field for $P$. Unlike the case in the previous section, however, in order to define a connection we require a Higgs field to satisfy a slightly different condition. Recall that a Higgs field for an $LG\rtimes S^1$-bundle $P$ satisfies
$$
\Phi(p (\gamma, \phi) ) = \rho_{\phi}^{-1} \left(ad(\gamma^{-1}) \Phi(p) + \gamma^{-1} \partial \gamma\right).
$$
It will be instructive to define a Higgs field for $P$ given the bundles $\widetilde{P} \to Y \to M$ now since we will need this later to show that the construction is invertible. Define then, the map $\Phi \colon P \to L\fg$ by
$$
\Phi(f) = \tilde{A}(\partial f).
$$
This is a Higgs field since
\begin{align*}
\Phi(f (\gamma, \phi) )	&= \tilde{A}(\rho_{\phi}^{-1}(\partial f \gamma) + \rho_{\phi}^{-1}(\gamma^{-1}\partial \gamma )\\ 
					&= \tilde{A}(\rho_{\phi}^{-1}(\partial f \gamma)) + \iota_{\rho_\phi (f)}( \rho_{\phi}^{-1}(\gamma^{-1}\partial \gamma )\\
					&= ad(\rho_{\phi}^{-1} (\gamma^{-1})) \tilde{A}(\rho_{\phi}^{-1}(\partial f )) + \rho_{\phi}^{-1}(\gamma^{-1}\partial \gamma )\\
					&= \rho_{\phi}^{-1} \left(ad(\gamma^{-1}) \tilde{A}(\partial f ) + \gamma^{-1}\partial \gamma \right)\\
					&= \rho_{\phi}^{-1} \left(ad(\gamma^{-1}) \Phi(f) + \gamma^{-1} \partial \gamma\right).
\end{align*}

To define a connection on $\widetilde{P} = (P\times G \times S^1)/LG \rtimes S^1$ we need to be able to write a form on $P\times G \times S^1$ which is zero on vertical vectors (with respect to the $LG\rtimes S^1$ action) and invariant under the $LG\rtimes S^1$ action (so as to ensure that it is well-defined). Thus we need to calculate the action of $LG\rtimes S^1$ on: the connection, $(A,a),$ on $P,$ the Higgs field, $\Phi,$ on $P$ and the Maurer-Cartan forms $\Theta$ and $d\theta$ on $G$ and $S^1$ respectively. Then we can combine these in an invariant way. We can calculate the action of $(\gamma, \phi)$ on the connection on $P$:
\begin{align*}
(\gamma, \phi)^*(A,a)(X)	&= (A,a)(X(\gamma, \phi))\\
						&= ad(\gamma, \phi)^{-1}(A, a)(X)\\
						&= \left( \rho_{\phi}^{-1} \left(ad(\gamma^{-1})A(X) - \gamma^{-1}\partial \gamma\, a(X)\right), a(X)\right),
\end{align*}
and we know that the Higgs field satisfies
$$
\Phi(p (\gamma, \phi) ) = \rho_{\phi}^{-1} \left(ad(\gamma^{-1}) \Phi(p) + \gamma^{-1} \partial \gamma\right),
$$
and the Maurer-Cartan form on $S^1$ is unchanged. To calculate the action on the Maurer-Cartan form on $G,$ consider a vector $(X, g\zeta, x_\theta) \in T_{(p, g, \theta)}(P\times G\times S^1).$ We have:
\begin{align*}
(X, g\zeta, x_\theta) (\gamma, \phi)	&=\frac{d}{dt}\bigg|_0 ( \gamma_X (t) (\gamma, \phi), \gamma(\theta +tx)^{-1}g \exp(t\zeta), \theta + tx - \phi)\\
							&= \frac{d}{dt}\bigg|_0  ( \gamma_X (t) (\gamma, \phi), (\gamma(\theta)^{-1} - \gamma(\theta)^{-1}\partial \gamma(\theta)\gamma(\theta)^{-1}tx )g (1+t\zeta), \theta + tx - \phi)\\
							&= (X(\gamma, \phi), \gamma(\theta)^{-1}g\zeta - \gamma(\theta)^{-1}\partial \gamma(\theta)\gamma(\theta)^{-1}gx, x)\\
							&= (X(\gamma, \phi), \gamma(\theta)^{-1}g \{\zeta - x\, ad(g^{-1}) \partial \gamma(\theta)\gamma(\theta)^{-1} \}, x),
\end{align*}
and so
\begin{align*}
(\gamma, \phi)^* \Theta (g\zeta) 	&= \Theta_{\gamma(\theta)^{-1}g}(\gamma(\theta)^{-1}g \{\zeta - x\, ad(g^{-1}) \partial \gamma(\theta)\gamma(\theta)^{-1} \})\\
							&= \zeta - x\, ad(g^{-1}) \partial \gamma(\theta)\gamma(\theta)^{-1}.
\end{align*}
Now consider the form on $P \times G\times S^1$ given by
$$
\tilde{A} = ad(g^{-1})A + \Theta + ad(g^{-1})\Phi ( a + d\theta).
$$
This is invariant under the $LG\rtimes S^1$ action, for
\begin{align*}
(\gamma, \phi)^*\tilde{A}_{(p, g, \theta)}&(X, g\zeta, x_\theta)\\
				&= \tilde{A}_{(p(\gamma, \phi), \gamma(\theta)^{-1}g, \theta + \phi)}(X(\gamma, \phi),\gamma(\theta)^{-1}g \{\zeta - x\, ad(g^{-1}) \partial \gamma(\theta)\gamma(\theta)^{-1} \}, x)\\
				&= ad(g^{-1}\gamma(\theta))\rho_{\phi}^{-1} \left(ad(\gamma^{-1})A(X)_{\theta - \phi} - \gamma^{-1}\partial \gamma_{\theta - \phi}\, a(X)\right)\\
				&\qquad\qquad + \zeta - x\, ad(g^{-1}) \partial \gamma(\theta)\gamma(\theta)^{-1}\\
				&\qquad\qquad+ ad(g^{-1}\gamma(\theta))\rho_{\phi}^{-1} \left(ad(\gamma^{-1}) \Phi(p)_{\theta - \phi} + \gamma^{-1} \partial \gamma_{\theta - \phi}\right) \left(a(X) + x\right)\\
				&= ad(g^{-1}\gamma(\theta)) \left(ad(\gamma^{-1})A(X)_{\theta} - \gamma^{-1}\partial \gamma_{\theta}\, a(X)\right)\\
				&\qquad\qquad + \zeta - x\, ad(g^{-1}) \partial \gamma(\theta)\gamma(\theta)^{-1}\\
				&\qquad\qquad+ ad(g^{-1}\gamma(\theta))\left(ad(\gamma^{-1}) \Phi(p)_{\theta} + \gamma^{-1} \partial \gamma_{\theta}\right) \left(a(X) + x\right)\\
				&= ad(g^{-1}) A(X)_{\theta} + \zeta + ad(g^{-1})\Phi(p)_{\theta} \left(a(X) + x\right)\\
				&= \tilde{A}_{(p, g, \theta)}(X, g\zeta, x_\theta),
\end{align*}
by the calculations above. So for it to be well-defined on the quotient space we just need to check that it vanishes on vertical vectors. The vertical vector at the point $(p, g, \theta)$ generated by the vector $(\xi , x)$ is
\begin{align*}
V	&= \frac{d}{dt}\bigg|_0 (p, g, \theta)(\exp(t\xi), tx)\\
	&= \frac{d}{dt}\bigg|_0 (p(\exp(t\xi), tx), (1 - t\xi(\theta))g, \phi- tx)\\
	&= (\iota_{p}(\xi, x), - g\, ad(g^{-1})\xi(\theta), -x),
\end{align*}
and so
\begin{align*}
\tilde{A}(V)	&= ad(g^{-1})A(\iota_{p}(\xi, x))_\theta - ad(g^{-1})\xi(\theta) + ad(g^{-1})\Phi(p)(a(\iota_p (\xi. x)) - x)\\
			&= ad(g^{-1})\xi(\theta) - ad(g^{-1})\xi(\theta) + ad(g^{-1})\Phi(p)(x - x)\\
			&= 0.
\end{align*}
Thus we have defined a $G$-valued 1-form on $\widetilde{P}.$ $\tilde{A}$ is in fact a connection form, since if we evaluate it on the vertical vector generated by $\zeta\in \fg,$ that is, $\iota_{[p, g,\theta]}(\zeta) = (0, g\zeta, 0),$ we get $\tilde{A}(g\zeta) = \zeta$ and further,
\begin{align*}
\tilde{A}((X, g\zeta, x_\theta)h)	&= \tilde{A}(X, gh h^{-1}\zeta h, x_\theta)\\
						&= \left(ad(gh)^{-1}A + ad(h^{-1})\Theta + ad(gh)^{-1}\Phi ( a + d\theta)\right)(X, g\zeta, x_\theta)\\
						&= ad(h^{-1})\tilde{A}(X, g\zeta, x_\theta).
\end{align*}
To define a connection on the $S^1$-bundle $P/ LG$ we just take the projection of the $i\RR$-valued 1-form $a$ which is a connection form.

What remains to be shown now is that the constructions presented here for connections on $P,$ $\widetilde{P}$ and $Y$ are invertible. In particular, suppose we have the $LG\rtimes S^1$-bundle $P\to M$ with connection $(A,a)$ and Higgs field $\Phi$ and have constructed $\widetilde{P}\to Y\to M$ with connections $\tilde{A}$ and $\tilde{a}.$ Then if we construct the corresponding $LG\rtimes S^1$-bundle $P'$ (which is isomorphic to $P$ via the map $f \colon P \to P' ; \, p \mapsto f_p = (\theta \mapsto [p, 1, \theta])$) and the connection $(A', a')$ for $P',$ we would like to show that $f^*(A',a') = (A,a).$ Note that for the vector $X \in T_{p} P$ we have
$$
f_*X = (X, 0,0) \in T_{f_p}P'.
$$
Therefore,
\begin{align*}
f^*(A',a') (X) 	&= (A',a') (X,0,0)\\
			&= (\tilde{A}(X), a'(X))\\
			&= (A(X), a'(X))
\end{align*}
by the definition of $A'$ in terms of $\tilde{A}$ and $\tilde{A}$ in terms of $A$ and also $a'(X) = \tilde{a}(\tilde{\pi}_*X) = a(X).$ On the other hand, suppose we had the bundles $\widetilde{P} \to Y \to M$ with connections $\tilde{A}$ and $\tilde{a}$ and constructed $P\to M$ with connection $(A, a)$ and Higgs field $\Phi(f) =\tilde{A}(\partial f).$ Then we would like to show that if we construct the bundles $\widetilde{P}\to Y \to M$ with connections $\tilde{A}'$ and $\tilde{a}',$ we have $\tilde{A}' = f^*\tilde{A}$  where $f \colon \widetilde{P}' \xrightarrow{\sim} \widetilde{P}$ is the isomorphism given by $[f,g,\theta] \mapsto f(\theta)g.$ Note that at the point $[p, g, \theta]$ we have
\begin{align*}
f_*(X, g\zeta, x_\theta)	&= \frac{d}{dt}\bigg|_0 \gamma_{X(\theta + tx)} (t)g\exp(t\zeta)\\
				&= \frac{d}{dt}\bigg|_0 \gamma_{X(\theta)}(t)g +\partial\gamma_{X(\theta)} (0) x g + \gamma_{X(\theta)}(0) g\zeta\\
				&= X(\theta) g +\partial p(\theta) x g + p(\theta)g\zeta
\end{align*}
and therefore,
\begin{align*}
 f^*\tilde{A} (X, g\zeta, x_\theta)	&= \tilde{A} (X(\theta) g + \partial p(\theta) x g + p(\theta)g\zeta)\\
 						&= \tilde{A}(X(\theta) g) + x\tilde{A} (\partial p(\theta)  g) + \zeta\\
						&= ad(g^{-1})\tilde{A}(X(\theta)) + x\,ad(g^{-1})\tilde{A} (\partial p(\theta)) + \zeta
\end{align*}
while for $\tilde{A}'$ we have
\begin{align*}
\tilde{A}'(X, g\zeta, x_\theta)&= ad(g^{-1})A(X) + \zeta + ad(g^{-1})\Phi(p) ( a(X) + x)\\
					&= ad(g^{-1})\left(\tilde{A}(X) -a(X)\tilde{A}(\partial p)\right) + \zeta + ad(g^{-1})\tilde{A}(\partial p) ( a(X) + x)\\
					&=  f^*\tilde{A} (X, g\zeta, x_\theta).
\end{align*}
Thus we have proven the analogue of Proposition \ref{P:LG connection correspondences}

\begin{proposition}\label{P:LGxS^1 connection correspondences}
The correspondence from Proposition \ref{P:LGxS^1 correspondences} extends to a bijection between $G$-bundles with connection over $S^1$-bundles with connection and $\LGS$-bundles with connection and Higgs field.
\end{proposition}

\subsection{The string class and the first Pontrjagyn class}\label{SS:string class and p1}

Now that we have extended the correspondence from section \ref{S:Higgs fields, LG-bundles,...}, we are in a position to extend the result concerning the string class and the Pontrjagyn class (Theorem \ref{T:LGPont}). Recall that Theorem \ref{T:LGPont} extended Killingback's result to a general $LG$-bundle $P \to M$ by relating the string class of $P$ to the first Pontrjagyn class of the corresponding $G$-bundle $\widetilde{P} \to M\times S^1.$ In particular, the string class of $P$ is given by integrating $p_1(\widetilde{P})$ over the circle. We would like now to extend this further to the case where $P \to M$ is an $\LGS$-bundle and $\widetilde{P}$ is the corresponding $G$-bundle over a circle bundle $Y$ over $M.$ In this case we find that the string class is given by integrating the first Pontrjagyn class of $\widetilde{P}$ over the fibre of the circle bundle $Y$. In particular,
we have the following theorem

\begin{theorem}\label{T:LGxS^1Pont}
Let $P \to M$ be a principal $LG\rtimes S^1$-bundle and $\widetilde{P} \to Y \to M$ be the corresponding $G$-bundle over an $S^1$-bundle. Then the string class of $P$ is given by the integration over the fibre of the first Pontrjagyn class of $\widetilde{P}.$ That is,
$$
s(P) = \int_{S^1} p_1 (\widetilde{P}).
$$
\end{theorem}

\begin{proof}

We prove this in analogy with the proof of Theorem \ref{T:LGPont}, that is, by calculating the integral of the first Pontrjagyn class of $\widetilde{P}.$

Recall that the first Pontrjagyn class is given by
$$
p_1 = -\frac{1}{8\pi^2} \< \tilde{F}, \tilde{F}\>,
$$
where $\tilde{F} = d\tilde{A} + \tfrac{1}{2} [\tilde{A}, \tilde{A}]$ is the curvature of the connection $\tilde{A}$ corresponding to the pair $(A, \Phi)$ on $P.$ We have
\begin{multline*}
\tilde{F} = d (ad(g^{-1})A + \Theta + ad(g^{-1})\Phi ( a + d\theta))\\
\shoveright{+ \tfrac{1}{2}[ad(g^{-1})A + \Theta + ad(g^{-1})\Phi ( a + d\theta), ad(g^{-1})A + \Theta + ad(g^{-1})\Phi ( a + d\theta)]}\\
\shoveleft{\phantom{\tilde{F}} = d (ad(g^{-1})A + \Theta + ad(g^{-1})\Phi ( a + d\theta))}\\
+ \tfrac{1}{2}[ad(g^{-1})A, ad(g^{-1})A ] 
+ [ad(g^{-1})A ,  \Theta]
+ [ad(g^{-1})A , ad(g^{-1})\Phi ( a + d\theta)]\\ 
+ \tfrac{1}{2}[ \Theta ,  \Theta]
+ [ \Theta ,  ad(g^{-1})\Phi ( a + d\theta)] 
+ \tfrac{1}{2}[ad(g^{-1})\Phi ( a + d\theta), ad(g^{-1})\Phi ( a + d\theta)].
\end{multline*}
To calculate $d(ad(g^{-1})A + \Theta + ad(g^{-1})\Phi ( a + d\theta))$ we use
\begin{multline*}
d(ad(g^{-1})A + \Theta + ad(g^{-1})\Phi ( a + d\theta)) (( X, g \xi, x_\theta),( Y, g\zeta, y_\theta))\\
\shoveleft{\phantom{d(ad} = \tfrac{1}{2} \left\{( X, g \xi, x_\theta )(ad(g^{-1})A(Y)_\theta) 
- ( Y, g\zeta, y_\theta)(ad(g^{-1})A(Y)_\theta)\right.}\\
\shoveright{\left.- ad(g^{-1})A([X, Y])_\theta\right\}}\\
+  d\Theta\\
\shoveleft{\phantom{d(ad}+ \tfrac{1}{2} \left\{( X, g \xi, x_\theta)(ad(g^{-1})(a(Y) + y) \Phi(p)_\theta)\right.}\\
\shoveright{\left.- ( Y, g\zeta, y_\theta)(ad(g^{-1})(a(X) + x) \Phi(p)_\theta)
- ad(g^{-1})[x,y]\Phi(p)_\theta\right\},}\\
\end{multline*}
for tangent vectors $(X, g\xi, x_\theta)$ and $(Y, g\zeta, y_\theta)$ at the point $[p, g, \theta] \in \widetilde{P}.$ For the first term, calculate
\begin{align*}
( X, g & \xi, x_\theta )(ad(g^{-1})A(Y)_\theta)\\
& = \frac{d}{dt}\bigg|_0 (1 - t \xi)g^{-1} A_{\gamma_{p}(t)}(Y)_{(\theta + tx)}g(1 + t \xi)\\
& = \frac{d}{dt}\bigg|_0 g^{-1} A_{\gamma_{p}(t)}(Y)_{(\theta + tx)}g -  t \xi g^{-1} A_{\gamma_{p}(t)}(Y)_{(\theta + tx)}g + g^{-1} A_{\gamma_{p}(t)}(Y)_{(\theta + tx)}g t \xi\\
& = \frac{d}{dt}\bigg|_0 g^{-1} A_{\gamma_{p}(t)}(Y)_{\theta}g + g^{-1}\partial A(Y)_{\theta}xg -  \xi g^{-1} A(Y)_{\theta}g + g^{-1} A(Y)_{\theta}g \xi.\\
\end{align*}
Combining this with the other terms for the first derivative above, we have
$$
d(ad(g^{-1}) A) = ad(g^{-1}) dA - ad(g^{-1}) \partial A \wedge d\theta - [\Theta, ad(g^{-1}) A].
$$
For the last term, calculate
\begin{multline*}
( X, g  \xi, x_{ 1 \theta})(ad(g^{-1})(a(Y) + y) \Phi(p)_\theta)\\
\shoveleft{ = \frac{d}{dt}\bigg|_0 (1 - t \xi)g^{-1} (a_{\gamma_{p}(t)}(Y) + y)\Phi(\gamma_p (t))_{(\theta + tx)}g(1 + t \xi)}\\
\shoveleft{ = \frac{d}{dt}\bigg|_0 g^{-1} (a_{\gamma_{p}(t)}(Y) + y)\Phi(\gamma_p (t))_{(\theta + tx)}g}\\
-  t \xi g^{-1} (a_{\gamma_{p}(t)}(Y) + y)\Phi(\gamma_p (t))_{(\theta + tx)}g\\
\shoveright{+ g^{-1}(a_{\gamma_{p}(t)}(Y) + y)\Phi(\gamma_p (t))_{(\theta + tx)}g t \xi}\\
\shoveleft{= \frac{d}{dt}\bigg|_0 g^{-1} (a_{\gamma_{p}(t)}(Y) + y)\Phi(p)_{\theta }g + \frac{d}{dt}_{|_0} g^{-1} (a(Y) + y)\Phi(\gamma_p (t))_{\theta}g}\\
+ g^{-1} (a(Y)+ y)\partial\Phi(p)_{\theta}xg - \xi g^{-1} (a(Y) + y)\Phi(p)_{\theta}g\\
\shoveright{+ g^{-1}(a(Y) + y)\Phi(p)_{\theta }g \xi.}\\
\end{multline*}
Subtracting $ ( Y, g\zeta, y_\theta)(ad(g^{-1})(a(X) + x) \Phi(p)_\theta)$ from this gives
\begin{multline*}
d(ad(g^{-1})\Phi ( a + d\theta)) = ad(g^{-1})f \Phi + ad(g^{-1})d\Phi\wedge (a+ d\theta)\\
- [\Theta, ad(g^{-1})(a+d\theta)\Phi] - ad(g^{-1})a\partial \Phi \wedge d\theta.
\end{multline*}
We also have $d\Theta = -\tfrac{1}{2} [\Theta, \Theta].$ Therefore, the curvature of $\widetilde{P}$ is given by
\begin{multline*}
\tilde{F} = d (ad(g^{-1})A + \Theta + ad(g^{-1})\Phi ( a + d\theta))\\
+ \tfrac{1}{2}[ad(g^{-1})A, ad(g^{-1})A ] 
+ [ad(g^{-1})A ,  \Theta]
+ [ad(g^{-1})A , ad(g^{-1})\Phi ( a + d\theta)]\\ 
+ \tfrac{1}{2}[ \Theta ,  \Theta]
+ [ \Theta ,  ad(g^{-1})\Phi ( a + d\theta)] 
+ \tfrac{1}{2}[ad(g^{-1})\Phi ( a + d\theta), ad(g^{-1})\Phi ( a + d\theta)]\\
\shoveleft{\phantom{\tilde{F}} = ad(g^{-1}) dA 
- ad(g^{-1}) \partial A \wedge d\theta 
+ ad(g^{-1})f \Phi}\\ 
+ ad(g^{-1})d\Phi\wedge (a+ d\theta)
- ad(g^{-1})a\partial \Phi \wedge d\theta \\
\shoveright{+ \tfrac{1}{2}[ad(g^{-1})A, ad(g^{-1})A ]
+ [ad(g^{-1})A , ad(g^{-1})\Phi ( a + d\theta)].}\\
\end{multline*}
That is,
$$
\tilde{F} = ad(g^{-1})\left( F + f\Phi + \nabla \Phi \wedge (a + d\theta)\right).
$$
So the first Pontrjagyn class is
\begin{multline*}
p_1 = - \frac{1}{8\pi^2} \<\tilde{F}, \tilde{F}\>\\
%
%
\shoveleft{\phantom{p_1} = - \frac{1}{8\pi^2} \left\<  F + f\Phi + \nabla \Phi \wedge (a + d\theta),  F + f\Phi + \nabla \Phi \wedge (a + d\theta)\right\> }\\
\shoveleft{\phantom{p_1} = - \frac{1}{8\pi^2} \Big(\left\<  F + f\Phi ,  F + f\Phi \right\> 
-2 \left\<  F + f\Phi , \nabla \Phi \wedge (a + d\theta)\right\>}\\
\shoveright{- \left\<  \nabla \Phi \wedge (a + d\theta), \nabla \Phi \wedge (a + d\theta)\right\> \Big) }\\
%
%
%
%
%
\shoveleft{\phantom{p_1} = -\frac{1}{8\pi^2} \Big( \left\<  F + f\Phi ,  F + f\Phi \right\> 
-2 \left\<  F + f\Phi , \nabla \Phi \wedge a \right\> 
-2 \left\<  F + f\Phi , \nabla \Phi \right\> d\theta \Big).}\\
\end{multline*}
Thus, integrating $p_1$ over the fibre, we get
$$
-\frac{1}{4\pi^2}\int_{S^1} \langle F+ f\Phi, \nabla\Phi \rangle d\theta,
$$
which is the expression from Theorem \ref{T:LGxS^1string}.

\end{proof}

\section{String structures for $LG\rtimes \Diff(S^1)$-bundles}\label{S:Diff(S^1)}

So far in this chapter we have generalised the results from \cite{Murray:2003} to include the possibility of rotating loops. That is, we have worked with the semi-direct product $\LGS.$ We would like to conclude now with a brief outline of one way in which the results we have seen regarding $\LGS$ lead us to information about a more general situation. Namely, we shall consider the problem of lifting a bundle whose structure group is the semi-direct product $LG \rtimes \Diff(S^1).$ That is, we shall allow an action of the orientation preserving diffeomorphisms of the circle on the loops in $LG.$

The group $\Diff(S^1)$ has a well known central extension. In particular, the Lie algebra of this extension is the Virasoro algebra (see for example \cite{Mickelsson:1989}). In this section, we would like to consider the central extension of the semi-direct product above
$$
U(1) \to \widehat{LG \rtimes \Diff(S^1)} \to LG \rtimes \Diff(S^1).
$$

Thus far, we have seen that principal $LG$-bundles over $M$ correspond to principal $G$-bundles over $M\times S^1$ (via the caloron correspondence) and in the previous section we showed that isomorphism classes of principal $\LGS$-bundles are in bijective correspondence with isomorphism classes of principal $G$-bundles over principal $S^1$-bundles. If instead we considered a principal $G$-bundle over a general $S^1$ fibre bundle\footnote{Such bundles have structure group $\Diff(S^1)$ and give rise to principal $\Diff(S^1)$-bundles in a natural way.} we would find that these bundles correspond to principal $LG \rtimes \Diff(S^1)$-bundles.

Now let $R \to M$ be a principal $LG \rtimes \Diff(S^1)$-bundle. We are interested in finding the obstruction to lifting this bundle to an $\widehat{LG \rtimes \Diff(S^1)}$-bundle $\widehat{R}.$ The following result, due to Smale, gives us a way of using our previous results to solve this problem. Namely, we have
\begin{theorem}[\cite{Smale:1959}]\label{T:Diff(S^1)=S^1}
$\Diff(S^1)$ is homotopy equivalent to $S^1.$
\end{theorem}
This means that if $Y \to M$ is a $\Diff(S^1)$-bundle then its transition functions can be chosen to be valued in $S^1$ and so $Y$ actually admits an action of the circle (by identifying $Y$ locally with $S^1 \times U$ (for some open subset $U \subseteq M$) and rotating the $S^1$ factor). 
This makes $Y$ into a principal $S^1$-bundle. In particular, then, if we have a $G$-bundle $\widetilde{P}$ over an $S^1$ fibre bundle $Y \to M$ we can replace the $LG \rtimes \Diff(S^1)$-bundle in question with an $\LGS$-bundle. That is, $R$ has a reduction to a principal $\LGS$-bundle $P,$ so $R = P \times_{\LGS} LG \rtimes \Diff(S^1).$ We can thus give the lift of $R$ in terms of the central extension of $LG \rtimes \Diff(S^1)$ and the lift $\widehat{P}$ of $P.$ In particular, we have a bundle map
$$
\widehat{P} \times_{\widehat{\LGS}} \widehat{LG \rtimes \Diff(S^1)} \to P \times_{\LGS} LG \rtimes \Diff(S^1)
$$
given by
$$
[\hat{p}, \widehat{(\gamma, \varphi)}] \mapsto [p, (\gamma, \varphi)],
$$
where $\hat{p}$ is a lift of $p$ to $\widehat{P}$ and $\widehat{(\gamma, \varphi)}$ is a lift of $(\gamma, \varphi)$ to the central extension of $LG \rtimes \Diff(S^1).$ This map commutes with the homomorphism $\widehat{LG \rtimes \Diff(S^1)} \to LG \rtimes \Diff(S^1)$ and so $\widehat{P} \times_{\widehat{\LGS}} \widehat{LG \rtimes \Diff(S^1)}$ is a lift of $R.$


\renewcommand{\theequation}{\Alph{chapter}.\arabic{section}.\arabic{equation}}

\appendix
\chapter{Infinite-dimensional manifolds and Lie groups}\label{A:infinite-dimensional}


In this thesis we have largely been concerned with the loop group of a compact group. This is an example of an infinite-dimensional Lie group -- specifically, it is a \emph{Fr\'echet} Lie group. In this Appendix we collect some of the basic results on Fr\'echet manifolds and Lie groups. We follow closely the expositions presented in \cite{Hamilton:1982}, \cite{Milnor:1984} and \cite{Pressley-Segal}

\section{Fr\'echet spaces}\label{S:Frechet}

We will begin with some basic definitions and examples of the sorts of spaces we shall be dealing with. An infinite-dimensional manifold, like any manifold, is a topological space modelled on some sort of Euclidean space. In the case we are considering, this is a locally convex topological vector space called a \emph{Fr\'echet} space.

\begin{definition}
A \emph{Fr\'echet} space is a complete metrisable Hausdorff locally convex topological vector space, where by \emph{locally convex} we mean a space whose topology is generated from some family of seminorms.\footnote{An equivalent definition of local convexity for a topological vector space is that every neighbourhood of $0$ contains a neighbourhood which is convex. This is the definition used in \cite{Milnor:1984}.}
\end{definition}

Perhaps the most immediate example of a Fr\'echet space is given by any Banach space. In general, however, there are examples of Fr\'echet spaces which are not Banach spaces. The particular example we will consider is the space of all smooth maps\footnote{More generally, the space of smooth sections of a vector bundle over a compact manifold is also a Fr\'echet space. We shall restrict our interest however, to the case of a trivial bundle (that is, the space of all maps as above) since this covers the case we are really interested in -- the Lie algebra of the loop group, $\Map(S^1, \fg).$} from a compact manifold $X$ into a vector space $V,$ that is, the space $\Map(X,V).$ We define the topology on this space in terms of a collection of neighbourhoods of the zero map. (Since this is a topological vector space this will give the topology completely.) To do this, choose a small neighbourhood $E$ of $0 \in V.$ Then consider an open coordinate chart $U \subseteq X$ with local coordinates $x_1, \ldots, x_m$ and a compact set $K \subseteq U.$ We define a family of sub-basic neighbourhoods (for each choice of coordinate chart, compact set, neighbourhood of $0 \in V$ and non-negative integer $n$)
$$
N =\{ f \colon X \to V \mid \partial^k f / \partial x_{i_1} \ldots x_{i_k} \in E \ \forall \, x \in K, 0 \leq k \leq n, i_j \in \{1, \ldots, m\} \}.
$$
Finite intersections of sets of this form give the basic neighbourhoods for the topology on $\Map(X,V).$

The above example is important for our purposes since the special case of maps from the circle into the Lie algebra of a compact group $G$ will be the Fr\'echet space on which the loop group $LG$ is modelled.

\section{Groups of maps}\label{S:groups of maps}

Now that we have seen an example of a Fr\'echet space, we can give an example of an infinite-dimensional manifold modelled on this space. This is the space $\Map(X, G)$ of smooth maps from a compact manifold $X$ into a compact Lie group $G$ and it is in fact an example of an infinite-dimensional Lie group.

To define the coordinate charts for this manifold consider an open neighbourhood $U$ of the identity in $G$. Using the exponential map, this is homeomorphic to an open neighbourhood of the identity in $\fg,$ say $\tilde{U}.$ The set $\tilde{\mathcal{U}} := \Map(X, \tilde{U})$ is then an open neighbourhood of the identity in $\Map(X, \fg)$ and an atlas for $\Map(X, G)$ is given by the open sets $\mathcal{U} f$ (where $\mathcal{U} := \Map(X, U)$), which are also homeomorphic to $\tilde{\mathcal{U}}.$ The case where $X$ is the circle is the loop group $LG.$

Note that there is a slightly more general example given by taking sections of a fibre bundle over $X.$ Recall from the previous section that sections of a vector bundle form a Fr\'echet space. Given a fibre bundle $Y \xrightarrow{\pi} X$ we can associate to any section $f\colon X \to Y$ a vector bundle over $X,$ called the vertical tangent bundle to $f$ and denoted $T_{\text{vert}}Y_f,$ whose fibre at $x \in X$ is given by all vertical tangent vectors to $Y$ at $f(x).$ That is, $T_{\text{vert}}Y_f = \{ V \in T_{f(x)} Y \mid \pi_* V =0 \}.$ Then the sections of $T_{\text{vert}}Y_f \to X$ form a Fr\'echet space and there is a diffeomorphism from a neighbourhood of the zero section to a neighbourhood of the image of $f$ in $Y$ which serves as a coordinate chart.

\section{The path fibration}\label{S:path fibration}

In chapter \ref{C:characteristic classes} we made extensive use of a particular $\Omega G$-bundle called the \emph{path fibration}. This is a model for the universal $\Omega G$-bundle. In this section we shall explain why this is in fact a locally trivial $\Omega G$-bundle. Recall that the total space of the path fibration is defined as
$$
PG = \{ p: \RR \to G \mid p(0) = 1 \text{ and } p^{-1}dp \text{ is periodic} \}.
$$
We can equivalently view this as the space of connections on the trivial $G$-bundle over the circle, since if $p$ is a path in $G$ as above then $p^{-1} dp$ is a $\fg$-valued 1-form on $S^1$ and conversely, each connection form $A$ on the trivial $G$-bundle over $S^1$ uniquely determines a periodic path by solving the ordinary differential equation $A = p^{-1} dp$ subject to the initial condition $p(0) = 1.$ This means that $PG$ is contractible. Note that when viewed as the space of connections $\Omega G$ acts freely on the right of this space by gauge transformations. Notice also that if $p$ and $q$ are two paths in the same fibre of the projection $PG \xrightarrow{\pi} G$ (so $p(2\pi) = q(2\pi)$) then $p^{-1}q$ is a smooth based loop, since if $f(t) = (p^{-1}q)(t+2\pi)$ then $f$ satisfies the same differential equation as $p^{-1}q$ and $f(0) = 1$ so $f= p^{-1}q$ and thus $p^{-1}q$ is periodic. This means that $q = p\gamma$ for some $\gamma \in \Omega G$ and so $PG/\Omega G = G.$

For the local triviality of this bundle, consider an open neighbourhood $U$ of the identity in $G.$ We can define a map
$$
U \times \Omega G \xrightarrow{\sim} \pi^{-1}(U); \quad (g, \gamma) \mapsto p,
$$
where $p(t) = \exp (t \xi) \gamma(t)$ and $\exp(2\pi \xi) = g.$ The inverse of this map is given by
$$
p \mapsto (\pi(p) , \exp(t \pi(p))^{-1} p ).
$$
This gives us a trivialisation near the identity. To extend this to a local trivialisation for the entire bundle we consider the open cover $\{Uh\}$ for $h \in G.$ Let $\tilde{h}$ be a path ending at $h$ (that is, $\pi(\tilde{h}) = h$). Then the maps
$$
Uh \times \Omega G \xrightarrow{\sim} \pi^{-1}(Uh); \quad (g, \gamma) \mapsto p_h,
$$
for $p_h (t) = \tilde{h}(t) \exp (t \xi ) \gamma(t),$ give a local trivialisation for the path fibration. So we have that the path fibration is a model for the universal $\Omega G$-bundle.

\chapter{Classification of semi-direct product bundles}\label{A:classification}

\section{Classification of semi-direct product bundles}\label{S:KxH-bundles}


In section \ref{S:universal string class} we gave a model for the universal $L^\vee G$-bundle (where $L^\vee G$ is the group of smooth maps $[0, 2\pi] \to G$ with coincident endpoints) by utilising its description as the semi-direct product $\Omega^\vee G \rtimes G.$ Following those ideas we can actually give a classification theory for general $K \rtimes H$-bundles.

Suppose $K$ and $H$ are Lie groups and we have an action $\varphi \colon H\to \Aut (K).$ Then we can form the semi-direct product $K\rtimes H,$ where the multiplication is defined by
$$
(k_1, h_1) (k_2, h_2) = (k_1 \varphi_{h_1} (k_2), h_1 h_2).
$$
We can give a model for the classifying space $E(K\rtimes H)$ as follows. Consider the space $EK \times EH.$ This is contractible, since both $EK$ and $EH$ are. Suppose we can find a left action of $H$ on $EK.$ That is, some $\tilde{\varphi} \colon H \to \Diff (EK)$ such that $\tilde{\varphi}_{h_1} \tilde{\varphi}_{h_2} = \tilde{\varphi}_{h_1 h_2}.$ Suppose also that this action satisfies
$$
\tilde{\varphi}_h (xk) = \tilde{\varphi}_h (x) \varphi_h (k)
$$
for all $x \in EK.$ Then we can define a right action of $K\rtimes H$ on $EK\times EH$ by
$$
(x, y) (k, h) = (\tilde{\varphi}_{h^{-1}}(x k), yh),
$$
where $(x, y) \in  EK \times EH.$ This is clearly a right action since
\begin{align*}
(\tilde{\varphi}_{h_{1}^{-1}}(x k_1), y h_1) (k_2, h_2) &= (\tilde{\varphi}_{h_2^{-1}}(\tilde{\varphi}_{h_1^{-1}}(x k_1) k_2) , y h_1 h_2)\\
	&= (\tilde{\varphi}_{(h_1 h_2)^{-1}} ( x k_1 \varphi_{h_1}(k_2)), y h_1 h_2)\\
	&= (x, y) (k_1 \varphi_{h_1} (k_2), h_1 h_2).
\end{align*}
It is also free and transitive on fibres and so
$$
\xymatrix{EK \times EH \ar[d]\\ (EK\times EH)/(K\rtimes H)}
$$
is a model for the universal bundle. To see that $\tilde{\varphi}$ exists, consider the following construction of $EK$ \cite{Segal:1968} (see also \cite{Dupont:1978}). Let $\Delta^n$ be the standard $n$-simplex in $\RR^{n+1}.$ That is,
$$
\Delta^n = \{ (t_0, \ldots, t_n) \mid t_i \geq 0, \textstyle\sum t_i = 1\}.
$$
Then
$$
EK = \bigsqcup_{n \geq 0} \Delta^n \times K^{n+1} / \sim,
$$
where we make the identifications
$$
\left( (t_0, \ldots, t_{i-1}, 0, t_{i+1}, \ldots, t_n), (k_0, \ldots, k_n) \right) \sim \left( (t_0, \ldots, t_n), (k_0, \ldots, k_{i-1}, 1, k_{i+1}, \ldots, k_n) \right).
$$
Equivalently, we can think of $EK$ as the set of formal linear combinations of elements of $K:$
$$
EK = \left\{ \textstyle\sum t_i k_i \mid t_i \geq 0, \textstyle\sum t_i = 1, k_i \in K \right\}
$$
where in any given sum, only finitely many of the $t_i$'s are non-zero. Then $\tilde{\varphi}$ is given by
$$
\tilde{\varphi}_h \left( \textstyle\sum t_i k_i \right) =  \textstyle\sum t_i \varphi_h (k_i).
$$

Using this construction, we can also write down a classifying map for any $K\rtimes H$-bundle $P \xrightarrow{\pi} M.$ For this we will need a correspondence between these bundles and certain pairs of $K$-bundles and $H$-bundles. Let us briefly outline this correspondence now. First note that there is a homomorphism $K\rtimes H \to H$ and so we can form the associated $H$-bundle $P\times_{K\rtimes H} H \xrightarrow{\pi_H} M,$ where $[p,h] = [p(k',h'), h'^{-1}h],\,$ $[p, h]h' = [phh']$ and $\pi_H([p,h]) = \pi(p).$ Further, there's a free action of $K$ on $P$ that identifies $P\times_{K\rtimes H} H$ with $P/K.$ Namely, $pk = p(k, 1).$ Then we have that $P \xrightarrow{\pi_K} P\times_{K\rtimes H}H$ is a principal $K$-bundle.\footnote{For the proof of the local triviality of this bundle, see \cite{Kobayashi:1963}, Proposition 5.5, p 57.} Thus, we have constructed a $K$-bundle over an $H$-bundle out of the $K\rtimes H$-bundle $P$ that we started with. In addition, we have an action of $H$ on $P$ that covers the $H$ action on $P/K.$ That is, define $ph= p(1,h)$ and then $\pi_K(ph) = [p(1,h), 1] = [p, h] = [p,1]h = \pi_K(p) h.$ This $H$ action also has the property that $(ph)k = p(1, h)(k, 1) = p(\varphi_h(k), h) = (p\varphi_h(k))h.$ Therefore, we have constructed a $K$-bundle with a twisted $H$-equivariant action as above over an $H$-bundle:
$$
\xymatrix{ P\ar^{K, H}[d] \\
		P/K\ar^{H}[d]\\
		M}
$$

In fact, this construction is invertible. That is, given a $K$-bundle over an $H$-bundle that satisfies the properties above, we can construct a $K\rtimes H$-bundle. Suppose, then, that we have two Lie groups $K$ and $H$ with an action $\varphi \colon H\to \Aut(K)$ as above. Suppose also that we have a principal $K$-bundle $P\xrightarrow{\pi_K} P/K$ and a principal $H$-bundle $P/K \xrightarrow{\pi_H} M$ and that there's an $H$ action on $P$ covering that on $P/K$ and such that $(ph)k = (p\varphi_h(k))h.$ We can define an action of $K\rtimes H$ on $P$ by $p(k,h) = (pk)h.$ This is a right action since
\begin{align*}
p(k_1, h_1)(k_2,h_2)	&= (((pk_1)h_1)k_2)h_2\\
					&= (((pk_1)\varphi_{h_1}(k_2))h_1)h_2\\
					&= p(k_1 \varphi_{h_1}(k_2), h_1 h_2).
\end{align*}
It is a free action, for suppose that $p(k, h) =p.$ Then $(pk)h =p$ and so $\pi_K((pk)h) =\pi_K(p).$ But $\pi_K((pk)h) = \pi_K(pk)h$ and $\pi_K(pk) = \pi_K(p),$ so we have $\pi_K(p)h = \pi_K(p)$ and therefore $h=1$ since the $H$ action is free. But if $h=1$ we have that $pk =p$ and so $k=1.$ We also have that $P/(K\rtimes H) = (P/K)/H = M.$ To see that $P \to M$ is locally trivial, consider an open set $U\subset M$ over which $P/K$ is trivial. Then there exists a section $s \colon U\to P/K.$ Since $U$ is contractible, the pull-back $s^*P$ over $U$ is trivial and so there exists a section $s' \colon U \to s^*P.$ But a choice of section $s' \colon U\to s^*P$ is equivalent to a map $\sigma \colon U \to P$ such that $\pi_K(\sigma(x)) = s(x).$ That is, such that $\pi(\sigma(x)) = x.$ So $\sigma$ is a local section of $P\to M.$ Therefore, we have that $P\to M$ is a principal $K\rtimes H$-bundle.

Using this correspondence, we can write down a classifying map for $P.$ That is, a map $f \colon P \to EK \times EH$ such that $f(p(k, h)) = f(p) (k,h).$ Firstly, note that if $P \xrightarrow{\pi} M$ is a $\cG$-bundle then we can write the classifying map as follows: Let $\{U_\a\}$ be an open cover of $M$ over which $P$ is trivial. Then $\pi^{-1}(U_\a)$ is isomorphic to $U_\a \times \cG.$ Now choose local sections $s_\a \colon U_\a \to \pi^{-1}(U_\a)$ and define the functions $g_\a \colon \pi^{-1}(U_\a) \to \cG$ by $s_\a(m) = (m, g_\a(s_\a(m))),$ where we have used the isomorphism to identify $\pi^{-1}(U_\a)$ with $U_\a \times \cG.$ Now, let $\{\psi_\a\}$ be a partition of unity subordinate to $\{U_\a\}.$ Then define the map $f_{\cG} \colon P \to E\cG$ by
$$
f_{\cG} ( p) = \sum \psi_\a ( \pi(p)) g_\a (p).
$$
This is clearly $\cG$-equivariant and so defines the classifying map for $P.$

Now consider again the case where $\cG = K\rtimes H.$ Write the classifying map $f$ as a pair of functions $(f_K, f_H).$ Then we require that
$$
(f_K(p (k,h)), f_H(p (k,h))) = (\tilde{\varphi}_{h^{-1}}(f_K (p) k), f_H (p)h).
$$
Using the correspondence above, we can construct a pair of bundles $P \xrightarrow{\pi_K} P/K \xrightarrow{\pi_H} M.$ Define $f_H$ to be the classifying map of the $H$-bundle $P/K.$ To define $f_K,$ consider an open cover $\{U_\a\}$ of $M$ as above. Consider the cover $\{V_\a\}$ of $P/K$ where $V_\a = \pi_H^{-1}(U_\a).$ $P$ is trivial over $V_\a$ since we can construct a local section as follows. Identify $V_\a$ with $U_\a \times H.$ Then over the subset $U_\a \times \{ 1\},$ $P$ has a section, say $\sigma_\a.$ We can define a section of $P$ over $U_\a \times H$ by forcing $H$-equivariance. That is, by defining $\chi_\a (s_\a (m) h) \colon = \sigma_\a(m)h,$ where $s_\a$ is a local section of $P/K.$ So $\pi_K^{-1}(V_\a) \simeq U_\a \times H \times K.$ Then we can define the functions $k_\a$ as above and we see that $k_\a (p h) = \varphi_{h^{-1}}(k_\a (p))$ (which follows from the fact that $(pk)h = (ph)\varphi_{h^{-1}}(k)$). Therefore, if we choose partitions of unity $\{ \psi_a\}$ subordinate to $\{U_\a\}$ and $\{\chi_\a\}$ subordinate to $\{V_\a\},$ we can define
$$
f(p) = \left( \sum \chi_\a (\pi(p)) k_\a (p), \sum \psi_\a( \pi(p)) h_\a (\pi_K(p)) \right),
$$
which is $K\rtimes H$-equivariant because
\begin{align*}
f(p(k,h))	&= \left( \sum \chi_\a (\pi(p)) k_\a ((pk)h), \sum \psi_\a( \pi(p)) h_\a (\pi_K(ph)) \right)\\
		&= \left( \sum \chi_\a (\pi(p)) \varphi_{h^{-1}}(k_\a (p)k), \sum \psi_\a( \pi(p)) h_\a (\pi_K(p))h \right)\\
		&= f(p)(k,h).
\end{align*}
Thus $f$ is a classifying map for $P.$

\section{$\LGS$-bundles}\label{S:LGxS^1}

We have shown in the previous section that principal $K\rtimes H$-bundles are equivalent to $K$-bundles with a twisted equivariant $H$ action over $H$-bundles. Consider now the case where $K=LG$ and $H= S^1,$ as in chapter \ref{C:LGxS^1}. We have already seen (see section \ref{S:Higgs fields, LGxS^1-bundles,...}) that there is a bijective correspondence between isomorphism classes of principal $LG\rtimes S^1$-bundles and isomorphism classes of principal $G$-bundles over $S^1$-bundles. The result from section \ref{S:KxH-bundles}, however, implies that we could construct a principal $LG$-bundle over a circle bundle. Namely, the bundle $P \to P/LG = (P\times S^1)/LG\rtimes S^1$ is a principal $LG$ bundle. We would like to understand the relationship between the $LG$-bundle we have constructed and the $G$-bundle we have constructed in section \ref{S:Higgs fields, LGxS^1-bundles,...}. Consider the map
$$
\xymatrix{ P\ar^{f}[r] \ar[d]		&\widetilde{P}\ar[d] \\
		P/LG\ar[r]			&\widetilde{P}/G}
$$
given by $f(p) = [p,1,1]$ (and where the induced map $LG \to G$ is the homomorphism $\gamma \mapsto \gamma(1)$). This is a bundle map since
\begin{align*}
f(p\gamma)	&= [p (\gamma, 1), 1, 1]\\
			&= [p, \gamma(1), 1]\\
			&= [p, 1, 1]\gamma(1)\\
			&= f(p)\gamma(1).
\end{align*}
Therefore, we see that $\widetilde{P} \simeq P \times_{LG}G$ (via the isomorphism $[p,g] \mapsto [p,g,1]$). So $\widetilde{P}$ is given by extending the structure group of $P$ from $LG$ to $G.$

\addcontentsline{toc}{chapter}{Bibliography}

\bibliography{bibliography/BibTeX_References}

\begin{thebibliography}{10}

\bibitem{Bergman:2005}
A.~Bergman and U.~Varadarajan.
\newblock Loop groups, {K}aluza-{K}lein reduction and {M}-theory.
\newblock {\em J. High Energy Phys.}, (6):043, 28 pp. (electronic), 2005.

\bibitem{Bott-Tu}
R.~Bott and L.~W. Tu.
\newblock {\em Differential forms in algebraic topology}, volume~82 of {\em
  Graduate Texts in Mathematics}.
\newblock Springer-Verlag, New York, 1982.

\bibitem{Bouwknegt:2002}
P.~Bouwknegt, A.~L. Carey, V.~Mathai, M.~K. Murray, and D.~Stevenson.
\newblock Twisted {$K$}-theory and {$K$}-theory of bundle gerbes.
\newblock {\em Comm. Math. Phys.}, 228(1):17--45, 2002.

\bibitem{Brylinski:1993}
J.-L. Brylinski.
\newblock {\em Loop spaces, characteristic classes and geometric quantization},
  volume 107 of {\em Progress in Mathematics}.
\newblock Birkh\"auser Boston Inc., Boston, MA, 1993.

\bibitem{Brylinski:1994}
J.-L. Brylinski and D.~A. McLaughlin.
\newblock The geometry of degree four characteristic classes of line bundles on
  loop spaces {I}.
\newblock {\em Duke Math. J.}, 75(3):603--638, 1994.

\bibitem{Carey:2000}
A.~Carey and J.~Mickelsson.
\newblock A gerbe obstruction to quantization of fermions on odd-dimensional
  manifolds with boundary.
\newblock {\em Lett. Math. Phys.}, 51(2):145--160, 2000.

\bibitem{Carey:1997}
A.~Carey, J.~Mickelsson, and M.~Murray.
\newblock Index theory, gerbes, and {H}amiltonian quantization.
\newblock {\em Comm. Math. Phys.}, 183(3):707--722, 1997.

\bibitem{Carey:2005}
A.~L. Carey, S.~Johnson, M.~K. Murray, D.~Stevenson, and B.-L. Wang.
\newblock Bundle gerbes for {C}hern-{S}imons and {W}ess-{Z}umino-{W}itten
  theories.
\newblock {\em Comm. Math. Phys.}, 259(3):577--613, 2005.

\bibitem{Carey:2002}
A.~L. Carey and J.~Mickelsson.
\newblock The universal gerbe, {D}ixmier-{D}ouady class, and gauge theory.
\newblock {\em Lett. Math. Phys.}, 59(1):47--60, 2002.

\bibitem{Carey:2000a}
A.~L. Carey, J.~Mickelsson, and M.~K. Murray.
\newblock Bundle gerbes applied to quantum field theory.
\newblock {\em Rev. Math. Phys.}, 12(1):65--90, 2000.

\bibitem{Carey:1991}
A.~L. Carey and M.~K. Murray.
\newblock String structures and the path fibration of a group.
\newblock {\em Comm. Math. Phys.}, 141(3):441--452, 1991.

\bibitem{Cheeger:1985}
J.~Cheeger and J.~Simons.
\newblock Differential characters and geometric invariants.
\newblock In {\em Geometry and topology ({C}ollege {P}ark, {M}d., 1983/84)},
  volume 1167 of {\em Lecture Notes in Math.}, pages 50--80. Springer, Berlin,
  1985.

\bibitem{Chern:1974}
S.~S. Chern and J.~Simons.
\newblock Characteristic forms and geometric invariants.
\newblock {\em Ann. of Math. (2)}, 99:48--69, 1974.

\bibitem{Dupont:1978}
J.~L. Dupont.
\newblock {\em Curvature and characteristic classes}.
\newblock Lecture Notes in Mathematics, Vol. 640. Springer-Verlag, Berlin,
  1978.

\bibitem{Ekstrand:2000}
C.~Ekstrand and J.~Mickelsson.
\newblock Gravitational anomalies, gerbes, and {H}amiltonian quantization.
\newblock {\em Comm. Math. Phys.}, 212(3):613--624, 2000.

\bibitem{Garland:1988}
H.~Garland and M.~K. Murray.
\newblock Kac-{M}oody monopoles and periodic instantons.
\newblock {\em Comm. Math. Phys.}, 120(2):335--351, 1988.

\bibitem{Giraud:1971}
J.~Giraud.
\newblock {\em Cohomologie non ab\'elienne}.
\newblock Springer-Verlag, Berlin, 1971.
\newblock Die Grundlehren der mathematischen Wissenschaften, Band 179.

\bibitem{Gomi:2003}
K.~Gomi.
\newblock Connections and curvings on lifting bundle gerbes.
\newblock {\em J. London Math. Soc. (2)}, 67(2):510--526, 2003.

\bibitem{Hamilton:1982}
R.~S. Hamilton.
\newblock The inverse function theorem of {N}ash and {M}oser.
\newblock {\em Bull. Amer. Math. Soc. (N.S.)}, 7(1):65--222, 1982.

\bibitem{Hatcher}
A.~Hatcher.
\newblock {\em Algebraic topology}.
\newblock Cambridge University Press, Cambridge, 2002.

\bibitem{Heitsch:1974}
J.~L. Heitsch and H.~B. Lawson, Jr.
\newblock Transgressions, {C}hern-{S}imons invariants and the classical groups.
\newblock {\em J. Differential Geometry}, 9:423--434, 1974.

\bibitem{Killingback:1987}
T.~P. Killingback.
\newblock World-sheet anomalies and loop geometry.
\newblock {\em Nuclear Phys. B}, 288(3-4):578--588, 1987.

\bibitem{Kobayashi:1963}
S.~Kobayashi and K.~Nomizu.
\newblock {\em Foundations of differential geometry. {V}ol {I}}.
\newblock Interscience Publishers, a division of John Wiley \& Sons, New
  York-Lond on, 1963.

\bibitem{Kobayashi:1969}
S.~Kobayashi and K.~Nomizu.
\newblock {\em Foundations of differential geometry. {V}ol. {II}}.
\newblock Interscience Tracts in Pure and Applied Mathematics, No. 15 Vol. II.
  Interscience Publishers John Wiley \& Sons, Inc., New York-London-Sydney,
  1969.

\bibitem{Lawson:1989}
H.~B. Lawson, Jr. and M.-L. Michelsohn.
\newblock {\em Spin geometry}, volume~38 of {\em Princeton Mathematical
  Series}.
\newblock Princeton University Press, Princeton, NJ, 1989.

\bibitem{MacLane:1971}
S.~MacLane.
\newblock {\em Categories for the working mathematician}.
\newblock Springer-Verlag, New York, 1971.
\newblock Graduate Texts in Mathematics, Vol. 5.

\bibitem{Massey:1991}
W.~S. Massey.
\newblock {\em A basic course in algebraic topology}, volume 127 of {\em
  Graduate Texts in Mathematics}.
\newblock Springer-Verlag, New York, 1991.

\bibitem{Meinrenken:2003}
E.~Meinrenken.
\newblock The basic gerbe over a compact simple {L}ie group.
\newblock {\em Enseign. Math. (2)}, 49(3-4):307--333, 2003.

\bibitem{Mickelsson:1989}
J.~Mickelsson.
\newblock {\em Current algebras and groups}.
\newblock Plenum Monographs in Nonlinear Physics. Plenum Press, New York, 1989.

\bibitem{Mickelsson:2008}
J.~Mickelsson.
\newblock From gauge anomalies to gerbes and gerbal actions.
\newblock math-ph/0812.1640, 2008.

\bibitem{Milnor:1984}
J.~Milnor.
\newblock Remarks on infinite-dimensional {L}ie groups.
\newblock In {\em Relativity, groups and topology, {II} ({L}es {H}ouches,
  1983)}, pages 1007--1057. North-Holland, Amsterdam, 1984.

\bibitem{Murray:1996}
M.~K. Murray.
\newblock Bundle gerbes.
\newblock {\em J. London Math. Soc. (2)}, 54(2):403--416, 1996.

\bibitem{Murray:2000}
M.~K. Murray and D.~Stevenson.
\newblock Bundle gerbes: stable isomorphism and local theory.
\newblock {\em J. London Math. Soc. (2)}, 62(3):925--937, 2000.

\bibitem{Murray:2001}
M.~K. Murray and D.~Stevenson.
\newblock Yet another construction of the central extension of the loop group.
\newblock In {\em Geometric analysis and applications (Canberra, 2000)},
  volume~39 of {\em Proc. Centre Math. Appl. Austral. Nat. Univ.}, pages
  194--200. Austral. Nat. Univ., Canberra, 2001.

\bibitem{Murray:2003}
M.~K. Murray and D.~Stevenson.
\newblock Higgs fields, bundle gerbes and string structures.
\newblock {\em Comm. Math. Phys.}, 243(3):541--555, 2003.

\bibitem{Narasimhan:1961}
M.~S. Narasimhan and S.~Ramanan.
\newblock Existence of universal connections.
\newblock {\em Amer. J. Math.}, 83:563--572, 1961.

\bibitem{Narasimhan:1963}
M.~S. Narasimhan and S.~Ramanan.
\newblock Existence of universal connections. {II}.
\newblock {\em Amer. J. Math.}, 85:223--231, 1963.

\bibitem{Paycha:2004}
S.~Paycha and S.~Rosenberg.
\newblock Traces and characteristic classes on loop spaces.
\newblock In {\em Infinite dimensional groups and manifolds}, volume~5 of {\em
  IRMA Lect. Math. Theor. Phys.}, pages 185--212. de Gruyter, Berlin, 2004.

\bibitem{Pressley-Segal}
A.~Pressley and G.~Segal.
\newblock {\em Loop groups}.
\newblock Oxford Mathematical Monographs. The Clarendon Press Oxford University
  Press, New York, 1986.
\newblock Oxford Science Publications.

\bibitem{Schreiber:2007}
U.~Schreiber, C.~Schweigert, and K.~Waldorf.
\newblock Unoriented {WZW} models and holonomy of bundle gerbes.
\newblock {\em Comm. Math. Phys.}, 274(1):31--64, 2007.

\bibitem{Segal:1968}
G.~Segal.
\newblock Classifying spaces and spectral sequences.
\newblock {\em Inst. Hautes \'Etudes Sci. Publ. Math.}, (34):105--112, 1968.

\bibitem{Smale:1959}
S.~Smale.
\newblock Diffeomorphisms of the {$2$}-sphere.
\newblock {\em Proc. Amer. Math. Soc.}, 10:621--626, 1959.

\bibitem{Stevenson:comm}
D.~Stevenson.
\newblock Private communication.

\bibitem{Sury:2004}
B.~Sury, T.~Wang, and F.-Z. Zhao.
\newblock Identities involving reciprocals of binomial coefficients.
\newblock {\em J. Integer Seq.}, 7(2):Article 04.2.8, 12 pp. (electronic),
  2004.

\bibitem{Witten:1985}
E.~Witten.
\newblock Global anomalies in string theory.
\newblock In {\em Symposium on anomalies, geometry, topology (Chicago, Ill.,
  1985)}, pages 61--99. World Sci. Publishing, Singapore, 1985.

\end{thebibliography}


\end{document}